\documentclass[11pt,a4paper]{article}

\usepackage[bbgreekl]{mathbbol}
\usepackage{amsfonts}

\DeclareSymbolFontAlphabet{\mathbb}{AMSb}
\DeclareSymbolFontAlphabet{\mathbbl}{bbold}

\usepackage{amsmath,amssymb,nccmath}
\usepackage{theorem}
\usepackage{graphicx}
\usepackage{tcolorbox}
\usepackage{epstopdf}
\usepackage{subfig}
\usepackage{array}
\usepackage{subfig}
\usepackage{xspace}
\usepackage{footmisc}
\usepackage{lipsum}
\usepackage{afterpage}
\usepackage{dcolumn,booktabs}
\newcolumntype{d}[1]{D{.}{.}{#1}}
\allowdisplaybreaks

\usepackage{epstopdf}

\usepackage[margin=1in]{geometry} 

\usepackage{setspace}
\usepackage{pdflscape}
\usepackage{multirow}
\usepackage{dsfont}
\usepackage{hhline}
\usepackage{rotating}
\usepackage{float}
\usepackage{natbib}

\usepackage{enumerate}
\usepackage{color,colortbl}
\usepackage{subeqnarray}

\usepackage{textcomp}
\usepackage{float}
\usepackage{lscape}
\usepackage{epsfig}
\usepackage{tablefootnote}

\newtheorem{assumption}{Assumption}[section]
\newtheorem{example}{Example}[section]
{\theorembodyfont{\upshape} }

\newtheorem{proposition}{Proposition}
\newtheorem{lemma}{Lemma}
\newtheorem{definition}{Definition}

\newtheorem{corollary}{Corollary}
\newtheorem{theorem}{Theorem}

\providecommand{\keywords}[1]
{
  \small	
  \textbf{\textit{Keywords:}} #1
}

\newcommand{\proof}{\noindent \textbf{Proof:}\, }

\usepackage{color}
\usepackage{multirow, makecell,nccmath,amsmath,amsfonts}
\usepackage{soul,verbatim}
\usepackage{hyperref}
\usepackage{thmtools}
\usepackage[justification=centering]{caption}

\DeclareMathAlphabet\mathbfcal{OMS}{cmsy}{b}{n}

\newcommand{\V}{\overline{V}}
\newcommand{\bvarphi}{\pmb{\varphi}}

\newcommand{\bnu}{\pmb{\nu}}

\newcommand{\bu}{\mathbf{u}}
\newcommand{\bs}{\mathbf{s}}

\newcommand{\tzeta}{\tilde{\zeta}}

\newcommand{\mA}{\mathcal{A}}

\newcommand{\mX}{\mathcal{X}}
\newcommand{\mS}{\mathcal{S}}
\newcommand{\omS}{\overline{\mathcal{S}}}
\newcommand{\mbS}{\mathbfcal{S}}
\newcommand{\mU}{\mathcal{U}}
\newcommand{\omU}{\overline{\mathcal{U}}}
\newcommand{\mbU}{\mathbfcal{U}}
\newcommand{\mB}{\mathcal{B}}
\newcommand{\mD}{\mathcal{D}}
\newcommand{\mC}{\mathcal{C}}
\newcommand{\mM}{\mathcal{M}}

\newcommand{\mMS}{\mathcal{M}(\mathcal{S})}

\newcommand{\mplusMS}{\mathcal{M}_+(\mathcal{S})}
\newcommand{\mMpS}{\mathcal{M}_{\mathbb{P}}(\mathcal{S})}
\newcommand{\mMp}{\mathcal{M}_{\mathbb{P}}}

\newcommand{\mBS}{\mathcal{B}(\mathcal{S})}

\newcommand{\mBU}{\mathcal{B}(\mathcal{U})}

\newcommand{\oV}{\overline{V}}
\newcommand{\of}{\overline{f}}

\newcommand{\oq}{\overline{q}}
\newcommand{\oQ}{\overline{Q}}
\newcommand{\ovr}{\overline{r}}
\newcommand{\ovR}{\overline{R}}
\newcommand{\oJ}{\overline{J}}

\newcommand{\oh}{\overline{h}}
\newcommand{\ovh}{\overline{h}}
\newcommand{\orho}{\overline{\rho}}

\begin{document}

\title{\bfseries {\Large Measurized Markov Decision Processes }
}

\author{\small Daniel Adelman, Alba V. Olivares-Nadal\thanks{The University of Chicago, Booth School of Business, \textit{Dan.Adelman@chicagobooth.edu}, UNSW Business School, \textit{a.olivares\_nadal@unsw.edu.au}}}

\date{}
\maketitle
\vspace{-2em}
\thispagestyle{empty}


\abstract{
In this paper, we explore lifting Markov Decision Processes (MDPs) to the space of probability measures and consider the so-called {\it measurized} MDPs: deterministic processes where states are probability measures on the original state space, and actions are stochastic kernels on the original action space. We show that measurized MDPs are a generalization of stochastic MDPs, thus the measurized framework can be deployed without loss of fidelity. Bertsekas and Shreve studied similar deterministic MDPs under the discounted infinite-horizon criterion in the context of universally measurable policies. Here, we also consider the long-run average reward case, but we cast lifted MDPs within the semicontinuous-semicompact framework of Hernández-Lerma and Lasserre. This makes the lifted framework more accessible as it entails (i) optimal Borel-measurable value functions and policies, (ii) reasonably mild assumptions that are easier to verify than those in the universally-measurable framework, and (iii) simpler proofs. In addition, we showcase the untapped potential of lifted MDPs by demonstrating how the measurized framework enables the incorporation of constraints and value function approximations that are not available from the standard MDP setting. Furthermore, we introduce a novel algebraic lifting procedure for any MDP, showing that non-deterministic measure-valued MDPs can emerge from lifting MDPs impacted by external random shocks.
}


\maketitle

\section{Introduction}

Markov decision processes (MDPs) are stochastic systems that can be influenced by actions, taken by the user. Every time an action is implemented, we accrue a revenue $r(s,u)$ that depends on the system's current state $s \in \mS$ and the action $u \in \mathcal{U}(s)$, where $\mathcal{U}(s)$ denotes the set of implementable actions from state $s$. Following this, the system transitions to the next state according to a stochastic kernel $Q$. The problem of determining an optimal control policy according to some performance criterion is the core of dynamic programming \citep{bellman}. The versatility of MDPs translates into diverse and copious applications. We find, for instance, prolific literature in the fields of revenue management \citep{zhang,IJOC}, supply chain management \citep{inventoryrouting,klabjan}, and healthcare \citep{schaefer,shechter}.

Traditionally, there are two criteria to optimize the policy within an infinite horizon in an MDP: the Discounted Infinite-Horizon criterion (DIH), and the Long-Run Average Reward (AR).  In both cases, their optimal {\it value functions} assign a value to each state $s$. 

\paragraph{\textcolor{black}{ The DIH criterion.}} According to the DIH criterion, the optimal value function is defined as

\begin{equation}\label{infinite_horizon}
V^*(s_0):=\sup_{\Pi}  \ \mathbb{E}_{s_0}^\pi \left[ \sum_{t=0}^\infty \alpha^t r(s_t,u_t) \right]\qquad \forall s_{0} \in \mathcal{S}. 
\end{equation}
where $\Pi$ is a set of policies $\pi$ and $\alpha \in (0,1)$ is the discount factor. Under certain assumptions, the optimal DIH value function $V^*(\cdot)$ and the optimal action to take at any given time is deterministic and can be obtained by solving the DIH {\it optimality equations} \cite{bellman}:
   \begin{equation}
   V^*(s) = \sup_{u \in \mathcal{U}(s)} \left\{ r(s, u) + \alpha \mathbb{E}_Q \left[ V^*(s')| s, u \right] \right\}, \qquad \forall s \in \mS, \label{DCOE} \tag{$\alpha$-DCOE}
   \end{equation}
which simultaneously obtain the optimal value function $V^*(\cdot)$ and the optimal decision to take from each state. 

%


%

Instead of working in $\mS$, we propose to generalize the state space of this {\it baseline} MDP so it encompasses all probability measures on $\mS$. 
Specifically, we derive and analyze a lifted MDP whose states are distributions $\nu$ over $\mS$  and whose controls $\varphi$ are stochastic kernels, with $\varphi(\cdot|s)$ being a probability distribution over $\mU(s)$. We show that such an MDP, which we call a {\it measurized} MDP,  generalizes the {\it original} MDP associated with \eqref{DCOE}, and thus one can work in the measurized framework without loss of fidelity. We also show that the measurized MDP is deterministic with optimality equations

\begin{align} \label{M-DCOE-simple}
\overline{V}^*(\nu)&=\sup_{\varphi \in \Phi} \ \left\{ \ovr(\nu,\varphi) + \alpha   \overline{V}^*(\nu')\right\}, \qquad \forall \nu \in \mMpS,
\end{align}
where $\mMpS$ denotes the space of probability measures on $\mS$, $ \Phi$ is the set of implementable Markov decision rules, $\ovr(\nu,\varphi):=\mathbb{E}_\varphi \mathbb{E}_\nu[r(s,u)]$ is the expected reward and $\nu'$ is the future state distribution, which is {\it known} and depends deterministically on $\nu$, $\varphi$ and $Q$. 

Lifting problems to higher dimensional spaces in order to endow them with properties and constraints that are not reachable from their original spaces has been explored before in other contexts. For example, in integer programming, lifting is employed to enhance minimal cover inequalities, a common component in contemporary branch-and-bound solvers \citep{IP}. Support Vector Machines may transform the covariates in a classification problem through a high-dimensional application \citep{smola_learning}. This procedure may enable the construction of a linear separating hyperplane in the lifted space, even when no such hyperplane exists in the original space. In a dynamic setting, a similar concept has been studied in the uncontrolled literature \citep{Tweedie} through measure-valued Markov chains. Lifting MDPs to the space of probability measures has rarely been explored in the literature and is typically confined to specific settings. For instance, \cite{bauerle} and \cite{carmona} leverage the structure of Mean-Field MDPs (MFMDPs) to derive deterministic optimality equations similar to \eqref{M-DCOE-simple} in the absence of common noise. To the best of the authors' knowledge, the only exceptions to these setting-specific approaches are \cite{BSdeterministic,bertsekas} and this paper, where we develop a general framework applicable to any standard MDP.

Although the foundational work in \cite{BSdeterministic} and  \citet[Chapter 9]{bertsekas} laid important groundwork, it has had limited practical impact, and, to our knowledge, little research has been conducted in this direction since its publication. Our work aims to resurface, simplify and extend the theory of deterministic MDPs introduced therein by making the lifted framework more accessible and practical in several ways. First, we use simpler controls to make the model more interpretable. More specifically, a key difference between the deterministic MDPs in \cite{BSdeterministic, bertsekas}, which we refer to as {\it BS-deterministic} MDPs, and the lifted MDPs in this paper, termed {\it measurized} MDPs, lies in their control structures. BS-deterministic MDPs use joint probability distributions \( \gamma \) over the coupled state-action space, whereas measurized MDPs employ stochastic kernels \( \varphi \) of the original action space. The use of \( \gamma \) requires an additional constraint to ensure that its marginal distribution matches the current state \( \nu \), adding complexity to the model. In contrast, using \( \varphi \) provides a simpler and more intuitive formulation as these controls can be seen as Markov decision rules in the original MDP. We employ the Radon-Nikodym derivative to demonstrate that these two formulations are, in fact, equivalent. 

Second, we recast BS-deterministic MDPs within the {\it semicontinuous-semicompact} framework of the book from Hernández-Lerma and Laserre \cite{lerma}. In other words, we focus on the existence of Borel measurable, rather than semianalytic or universally measurable, value functions and optimal strategies. While this may appear less general, it offers three key advantages: (i) it enables us to derive {\it computable} results, at least in principle; (ii) our assumptions are reasonably mild, easy to verify, and sufficiently broad to encompass most practical control problems; and (iii) it simplifies the proofs. To elaborate on points (i) and (ii), \cite{BSdeterministic} assumes that the reward function \( r \) is upper semianalytic and that the space of feasible state-action pairs, \( \mathbb{K} \), is analytic. These conditions result in semianalytic optimal value functions and Markov \(\epsilon\)-optimal policies. While Borel-measurable value functions and policies can be derived within the {\it universally measurable} framework of \cite{BSdeterministic,bertsekas}, the assumptions required are often difficult to verify in practice. For example, \cite[Theorem 13]{BSdeterministic} requires \( r \) to be upper semicontinuous (u.s.c.), the function $(s,u) \to Q(ds'|s,u)$ to be continuous, \( \mathcal{U} \) to be compact, and \( \mathbb{K} \) to be expressible as an infinite union of nested closed sets \( \mathbb{K}^j \), with the additional condition that \( \lim_{j \to \infty} \inf_{(s,u)\in \mathbb{K}^j \setminus \mathbb{K}^{j-1}} r(s,u)=\infty \). In contrast, we adopt the standard assumptions of the {\it semicontinuous-semicompact} framework from \cite{lerma}. Specifically, while we also assume \( r \) is u.s.c. and $Q$ is strongly continuous, we replace the other assumptions in \citet[Theorem 13]{BSdeterministic} with the conditions that \( r \) is sup-compact in \( \mathbb{K} \) and that the transition kernel \( Q \) is strongly continuous. In addition, the semicontinuous-semicompact framework underpins point (iii)  as it {\it ``does not require a mathematical background beyond the graduate level"} \citep{lerma}. More specifically, proofs rely on the measurable selection theorem, which is more manageable than the exact selection theorem employed in \cite{BSdeterministic}. Furthermore, in \cite{BSdeterministic, bertsekas} BS-deterministic MDPs are used as a tool to prove theoretical results for the original MDP, such as the existence of stationary optimal policies and valid optimality equations. In contrast, we show that a measurized MDP inherits key assumptions from its baseline MDP when formulated in the semicontinuous-semicompact framework. This enables us to directly apply the established theory in \cite{lerma}, reversing the analysis: we derive properties of the lifted MDP by leveraging the simpler and already-proven results for the original MDP, further supporting point (iii).

Third, we extend the theory presented in \cite{BSdeterministic,bertsekas} to showcase the untapped potential of lifted MDPs. 

For example, we demonstrate how measurized MDPs connect to the Linear Programming (LP) formulation of the original MDP. In addition, we show that the sample average of the stochastic value function evaluated at states sampled from a particular state distribution converges to the measurized value of that distribution. This demonstrates that solving the measurized MDP accounts for simultaneously solving over (possibly infinite) paths of realizations of the original MDP. 
Moreover, we explore how lifting MDPs to higher-dimensional spaces introduces properties and constraints that are otherwise inaccessible in their original form. In particular, we illustrate how measurized MDPs enable the approximation of the value function $\overline{V}^*(\cdot)$ by a linear combination of basis functions $\overline{\phi}_k(\nu)$, which can easily accommodate approximations based on moments or divergence between measures.  For example, $\overline{\phi}_k(\nu)$ can measure the variance induced by the probability measure $\nu$ over the sample space $\mS$ or its distance to a distribution of reference $\nu_0$. 

To the best of our knowledge, such approximations have never been considered in the literature before. 

This feature, inherent to measurized MDPs, may potentially connect the prolific literature in Approximate Dynamic Programming \citep{powell,powell2,farias_ADP} with moments optimization \citep{moments1,moments2}. 

Furthermore, optimality equations \eqref{M-DCOE-simple} naturally allow for probabilistic constraints on space $\mathcal{S}$ by making $\Phi$ contingent on the current state distribution $\nu$. For instance, one may add risk constraints or bound the variance of future distribution $\nu'$, since this is known and can be computed using solely $\nu$, $\varphi$ and $Q$. This is in contrast to the standard MDP framework, where only deterministic constraints can be imposed at the present period \citep{altman}. Traditionally, one may also incorporate constraints that are fulfilled in expectation over the entire sample path \citep{altman,adelman_weakly}. More complex constraints have barely been studied in the literature. An exception is \cite{borkar}, which considered a finite-horizon MDP and bounded the CVaR of the accumulated costs at the terminal stage. 

The authors are not aware of other examples including such constraints in a DIH or AR MDP. 

Example \ref{ex:CVaR} in this paper shows how these constraints can be easily modelled in the measurized framework. In addition, \eqref{M-DCOE-simple} may enhance intuitiveness; for example, it allows us to understand what actions are taken most frequently according to the state distribution the MDP is in. This enables the interpretation of CVaR constraints as bounding the probability of taking {\it ``risky"} actions.

Fourth, we introduce a novel algebraic procedure for {\it ``measurizing"} a stochastic process. This structured approach enables the lifting of any MDP, ensuring consistency and enhancing intuitiveness. This lifting procedure allows us to generalize the observation made in \cite{bauerle} for MFMDPs by showing that any standard MDP with a transition kernel \( Q \) influenced by external random shocks leads to measure-valued MDPs with stochastic transitions. In other words, in this case the lifting produces a non-deterministic form of the optimality equations presented in \eqref{M-DCOE-simple}.

\paragraph{\textcolor{black}{ The AR criterion.}}  We further extend the theory in \cite{BSdeterministic} to consider the AR case. 
To our knowledge, lifted {\it deterministic} MDPs in the space of probability measures under the AR criterion have only been explored for MFMDPs without common noise, as in \cite{bauerle}. In contrast, our approach applies to any standard MDP.  A key component of the extension to the AR case is the equilibrium problem, which aims to find the best steady-state $(\nu^*,\varphi^*)$ maximizing $(1 - \alpha) \overline{V}^*(\nu)$. While the steady-state technique has been traditionally used for solving AR deterministic MDPs \citep{lerma_AR}, we extend its application to stochastic MDPs through the measurized framework. Unlike \cite{lerma_AR}, which requires the AR problem to be {\it ``dissipative"} and paths to converge to the optimal steady-state over time (Assumption 4.1 and Theorem 4.3 in \cite{lerma_AR}), our approach relies on more standard assumptions. Specifically, we assume the usual sufficient conditions for stochastic MDPs to satisfy the AR-optimality equations, along with the existence of an invariant state distribution for the AR-optimal decision rule.

Building on the equilibrium problem, we derive the measurized AR optimality equations as the Lagrangian dual of the equilibrium problem. In this derivation, the optimal bias function emerges as the dual variable associated with the equilibrium constraint, quantifying the gain in expected reward when deviating from an optimal equilibrium \((\nu^*, \varphi^*)\). This dual interpretation extends the perspective introduced by \cite{flynn} for deterministic MDPs in vector spaces of \(\mathbb{R}\) to Borel spaces and measurized MDPs, thus enabling its application to stochastic settings. While \cite{lerma_equilibrium} recognized the equilibrium problem as the dual of the stochastic AR optimality equations in Borel spaces, they neither considered the Lagrangian approach nor derived the measurized optimality equations from the former. Furthermore, our work highlights the equilibrium problem as a viable solution strategy for the stochastic AR problem, bridging a gap in the existing literature.  

Our approach also connects to the concept of minimum and maximum pairs\footnote{A maximum pair \((\tilde{\nu}, \tilde{\pi})\) consists of an initial state distribution \(\tilde{\nu} \in \mMpS\) and a policy \(\tilde{\pi} \in \Pi\) that jointly maximize the expected AR value function \eqref{AR}, evaluated with respect to \(\tilde{\nu}\) instead of a specific state \(s\). This makes \((\tilde{\nu}, \tilde{\pi})\) the best pair of state distribution and policy under the AR criterion.}, which have been studied in stochastic MDPs \citep{kurano,lerma,lerma_further,yu_minimumpair}. Maximum pairs are inherently linked to the steady-state technique because, under certain conditions, if an optimal equilibrium \((\nu^*, \varphi^*)\) exists, then \((\nu^*, \pi^*)\), where \(\pi^* = \{\varphi^*\}_{t \geq 0}\) is a stationary policy, forms a maximum pair \citep{kurano,lerma_further,yu_minimumpair}. However, our focus differs: rather than using equilibria to analyze the existence of AR solutions, we assume the necessary conditions for their existence and treat the problem of finding the best equilibrium as a direct solution method for the stochastic MDP.

We also examine how other classical approaches for analyzing AR problems, namely, the optimality equations, the LP formulation, and the vanishing-discount argument (see, e.g., Chapters 5 and 6 of \cite{lerma}), translate to the measurized framework. In particular, we show that whenever these techniques apply to the original stochastic MDP, they extend to its lifted counterpart. Interestingly, this extension can be established through the steady-state approach described above.

By unifying these contributions within the measurized framework, our work brings together previously disconnected ideas such as the steady-state technique, minimum pairs, and the Lagrangian approach for stochastic MDPs. As a result, the measurized framework not only opens new avenues for analyzing AR stochastic MDPs but also broadens the applicability of earlier contributions. 

The theory developed in this paper for the AR problem builds on the unichain case, in which, under any stationary policy, the induced Markov chain consists of a single recurrent class. 

Analogously to the DIH case, the theoretical contributions in this paper show that one can work in the (deterministic) AR measurized framework instead of in the classic (stochastic) one without loss of optimality,
except in the multi-chain setting and the case of purely transient MDPs.



\subsection{Structure of the paper}
The remainder of this paper is organized as follows. We begin with a brief introduction to the notation used throughout the paper. Section \ref{sec:MDP} provides an overview of the notation and theory related to classical stochastic MDPs. To ensure the paper is self-contained, we draw on definitions and results from \cite{lerma}, to whom we are greatly indebted. Specifically, Sections \ref{sec:discounted} and \ref{sec:LP} introduce DIH MDPs and their LP formulation, respectively. For more detailed information on MDPs in general Borel state and action spaces, readers are referred to \cite{lerma}, and for a comprehensive review of MDPs in countable state and action spaces, \cite{puterman} is recommended.

Section \ref{sec:measurized_def} formally defines both BS-deterministic and measurized MDPs. In Section \ref{sec:change_variables} we use the Radon-Nikodym derivative to show the equivalence between the two lifted models,

and in Section \ref{sec:assumptions} we demonstrate that the lifted MDP inherits the essence of the assumptions adopted for its original MDP.

Section \ref{sec:advantages} discusses the advantages of adopting the measurized framework. In Section \ref{sec:OE}, we derive the deterministic optimality equations for measurized DIH MDPs. Additionally, we explore the latent potential of the lifted framework by providing examples of augmented constraints (Section \ref{sec:constraints}) and value function approximations (Section \ref{sec:approximations}) that are not accessible within the original framework.

Section \ref{sec:measurized} establishes the mathematical connection between the lifted and original DIH MDPs. More specifically, Sections \ref{sec:states}, \ref{sec:valuefunctions}, and \ref{sec:policies} relate stochastic and measurized states, value functions, and policies, respectively.  Section \ref{sec:dualizing} demonstrates the connection between the LP formulation of the original MDP and the lifted MDPs. Section \ref{sec:measurizing} presents a novel algebraic procedure for measurizing a stochastic process, offering a lifting method for any MDP.

Finally, we address the AR case in Sections \ref{sec:Measurized_AR} and \ref{sec:T1-T3}. Section \ref{sec:Measurized_AR} begins with a brief review of the stochastic AR value function and the associated optimality equations; we relegate the standard assumptions and supporting theoretical results (which may also be found in \cite{lerma}) to Appendix A.1. Section \ref{sec:connection_discounted} introduces the equilibrium problem, from which we derive the measurized AR optimality equations in Section \ref{sec:M-AROE}. Section \ref{sec:measurizedAR} then establishes the connection between the stochastic and measurized optimal solutions. Section \ref{sec:T1-T3} shows how other classical techniques for analyzing the AR problem (namely, the optimality inequalities, the LP formulation, and the vanishing-discount approach) extend to the measurized framework. The paper concludes with a discussion of future research directions.

\subsection{Notation}

To make it easier to follow, we briefly introduce the notation used throughout this paper. We will use bold typeface to denote vectors and calligraphic font for sets and spaces. Table \ref{table:notation} contains the notation of sets, spaces, operators and abbreviations.

\begin{table}[h!]
	\centering	\begin{small}
		\begin{tabular}{rl}
		\hline
Notation & Definition\\
\hline
\hline
$\mA^c$ & complement of $\mA$\\
$\mBS$ & Borel $\sigma$-algebra of $\mS$\\
$\mMS$ & space of measures defined over a sample space $\mS$\\
$\mplusMS$ & space of {\it positive} measures defined over a sample space $\mS$\\
$\mMpS$ & space of {\it probability} measures defined over a sample space $\mS$\\
$\mathcal{K}(\mU|\mS)$ & space of stochastic kernels\tablefootnote{A stochastic kernel on $\mU$ given $\mS$ is a function $Q(\cdot|\cdot)$ such that $Q(\cdot|s)$ is a probability measure on $\mU$ for each fixed $s \in \mS$, and $Q(\mA|\cdot)$ is a measurable function on $\mS$ for each fixed $\mA \in \mB(\mU)$ (Definition C.1 in \cite{lerma}).} on $\mU$ given $\mS$\\
$\mathbb{P}_\nu(\cdot)$ & probability law induced by $\nu$\\
$\mathbb{E}_\nu[\cdot]$ & expectation operator taken with respect to measure $\nu$\\
a.s. & almost surely\\
\hline
\end{tabular}
\end{small}
\caption{Notation of sets, spaces, operators and abbreviations. \label{table:notation}}
\end{table}

Table \ref{table:notationMDPs} introduces the notation used for the different MDPs. Although these MDPs are formally introduced later (standard MDPs are presented in Section \ref{sec:MDP}, BS-deterministic and measurized MDPs are introduced in Section \ref{sec:measurized_def}), we provide this comprehensive table here for reference. We will use an overline to represent objects on the lifted space of measures, although measures themselves will be represented by Greek letters.

\begin{table}[h!]
	\centering	\begin{small}
		\begin{tabular}{rccc}
		\hline
& Standard MDP & \multicolumn{2}{c}{Lifted MDPs}   \\
&&\cite{BSdeterministic} & Measurized MDP\\
\hline
\hline
Abbreviation & MDP  & BS-deterministic MDPs & $\mM$-MDP\\
State space & $\mS$ &  $\mMpS$ &  $\mMpS$\\
States & $s$ &  $\nu$ &  $\nu$\\
Action space & $\mU$ & $\mM_{\mathbb{P}}(\mS\times\mU)$ &  $\Phi$ \\
Feasible action set & $\mU(s)$ & $ \Xi(\nu)$ &  $\Phi(\nu)$\\
Actions & $u$ & $\gamma$ & $\varphi$\\
Feasible set of state-action pairs & $\mathbb{K}$ & &  $\overline{\mathbb{K}}$\\
\textcolor{black}{Selectors} & $f(s)$ &  & $\overline{f}(\nu)$ \\
Markov decision rules & $\varphi(u|s)$ & &  $\psi(\varphi|\nu)$\\
Markov policies & $\pi=\{\varphi_t\}_{t\geq 0}$ & & $\overline{\pi}=\{\psi_t\}_{t\geq 0}$\\
Reward function & $r(s,u)$ &$\ovR(\gamma)$ &$\ovr(\nu,\varphi)$\\
Transition kernel & $Q(\mA|s,u), \ \mA\in \mB(\mS)$ &  $\oQ(\mathcal{P}|\gamma), \ \mathcal{P}\in \mB(\mMpS)$ &  $\oq(\mathcal{P}|\nu,\varphi), \ \mathcal{P}\in \mB(\mMpS)$\\
\textcolor{black}{Deterministic transition function} & & $F(\gamma)$ & $F(\nu,\varphi)$\\
\textcolor{black}{DIH optimal value function} & $V^*(s)$ & $\oV^*(\nu)$  &  $\oV^*(\nu)$\\
\textcolor{black}{AR optimal value function} & $J^*(s)$ & & $\oJ^*(\nu)$\\
\hline
\end{tabular}
\end{small}
\caption{Notation for the standard and lifted MDPs. \label{table:notationMDPs}}
\end{table}

\section{Standard stochastic MDPs \label{sec:MDP}}

This section often quotes and summarizes definitions and theoretical results in \cite{lerma}. Additional definitions and results coming from this resource can be found in Appendix A. Following the notation in Table \ref{table:notationMDPs}, we formally define Markov Decision Models as follows.\\

\begin{definition}[\cite{lerma}, Definition 2.2.1]\label{def:MDM}
A Markov Decision Model is a five-tuple $(\mathcal{S},\mathcal{U},\{\mathcal{U}(s)| \ s \in \mathcal{S}\},Q,r)$ consisting of

\begin{itemize}
\item[(i)] a Borel space $\mathcal{S}$, called the state space and whose elements are referred to as {\it states}
\item[(ii)] a Borel space $\mathcal{U}$, called the {\it control or action set}
\item[(iii)] a family $\{\mathcal{U}(s)| \ s \in \mathcal{S}\}$ of nonempty measurable subsets $\mathcal{U}(s)$ of $\mathcal{U}$, where $\mathcal{U}(s)$ denotes the set of feasible actions when the system is in state $s$, and with the property that the set 
$$\mathbb{K}=\{ (s,u)| \ s\in \mathcal{S},\ u \in \mathcal{U}(s)\}$$
 of feasible state-action pairs is a measurable subset of $\mathcal{S} \times \mathcal{U}$
\item[(iv)] a stochastic kernel $Q$ on $\mathcal{S}$ given $\mathbb{K}$ called the transition law
\item[(v)] a measurable function $r:\ \mathbb{K} \rightarrow \mathbb{R}$ called the reward-per-stage function\\
\end{itemize}
\end{definition}

Now assume that, at each period $t$, we observe the history $h_t=(s_0,u_0,...,s_{t-1},u_{t-1},s_t)$ of the Markov Decision Model, with $(s_k,u_k)\in \mathbb{K}$ for all $k=0,...,t-1$ and $s_t \in \mathcal{S}$. We denote the set of all admissible histories as $\mathcal{H}_t=\mathbb{K}^t \times \mathcal{S}$. Then we define a randomized policy as follows.\\

\begin{definition}[\cite{lerma}, Definition 2.2.3]
A randomized control policy is a sequence $\pi=\{\pi_t\}_{t=0,1,...}$ of stochastic kernels on the action set $\mathcal{U}$ given the set of all admissible histories $\mathcal{H}_t$, satisfying
$$\pi_t(\mathcal{U}(s_t)|h_t)=1 \qquad\forall h_t \in \mathcal{H}_t, \ t=0,1,... $$
The set of all randomized policies is denoted by $\Pi$.\\
\end{definition}

Under certain conditions, the optimal policy may depend solely on the current state of the Markov process. In addition, we may sometimes find deterministic policies yielding the same expected reward as randomized policies. The following definition introduces these concepts.\\

\begin{definition}[\cite{lerma}, Definitions 2.3.1 and 2.3.2]\label{def:Phi}
Let $\Phi$ be the set of all stochastic kernels $\varphi$ on $\mathcal{U}$ given $\mathcal{S}$ such that $\varphi(\mU(s)|s)=1$ for all $s \in \mS$; i.e.

\begin{equation} \label{eq:Phi}
\Phi:=\{\varphi \in \mathcal{K}(\mU|\mS): \ \varphi(\mU(s)|s)=1 \quad \forall s \in \mS\}.
\end{equation}
We call $\varphi \in \Phi$ an implementable {\it Markov decision rule} as it depends solely on the current state $s$ rather than the entire history $h$. Then a {\it randomized Markov policy} is a sequence $\pi=\{\varphi_t\}_{t\geq 0}$ where $\varphi_t \in\Phi$ for all $t\geq 0$. The set of all Markov randomized policies is denoted by $\Pi^{MR}$. 

Let $\mathbb{F}$ be the set of all measurable functions $f: \mS \rightarrow \mU$ verifying $f(s) \in \mU(s)$ for all $s \in \mS$. Any function $f$ in $\mathbb{F}$ is called a {\it selector}. A policy $\pi=\{\varphi_t\}_{t\geq 0}$ such that $\varphi_t \in\Phi$ and there exists a $f_t \in \mathbb{F}$ satisfying $\varphi_t(f_t(s)|s)=1$ for all $t\geq 0$ is called a {\it deterministic Markov policy}. The set of all Markov deterministic policies is denoted by $\Pi^{MD}$. Moreover, if $f_t(\cdot)=f(\cdot)$ for all $t\geq 0$, then $\pi$ is called a {\it deterministic stationary policy} and we use the notation $\pi \in \Pi^{D}$.\\
\end{definition}

By definition, we have that $\Pi^D\subset \Pi^{MD} \subset \Pi^{MR} \subset \Pi$. Now we use the definition of a control policy to lay out the unique probability distribution in the space of admissible histories according to Ionescu-Tulcea Theorem (see Proposition C.10 of \cite{lerma}). This allows us to properly introduce discrete-time MDPs. \\
\begin{definition}[\cite{lerma}, Definition 2.2.4]\label{def:MDP}
Let $(\Omega,\mathcal{F})$ be the measurable space consisting of the sample space $\Omega=(\mathcal{S}\times \mathcal{U})^\infty$ and $\mathcal{F}$ the corresponding product $\sigma$-algebra. Let $\pi=\{\pi_t\}_{t\geq 0}$ be an arbitrary control policy and $\nu_0$ the initial distribution on $\mathcal{S}$. Denote by $P_{\nu_0}^\pi$ the unique probability measure supported on $\mathcal{H}_{\infty}$ (see Ionescu-Tulcea Theorem), verifying
\begin{enumerate}
\item[(i)] $P_{\nu_0}^\pi(s_0 \in \mA)=\nu_0(\mA) \qquad \forall \mA \in \mB(\mS)$
\item[(ii)] $P_{\nu_0}^\pi(u_t \in \mC| h_t)=\pi_t(\mC|h_t) \qquad \forall \mC \in \mB(\mU), t >0$
\item[(iii)] $P_{\nu_0}^\pi(s_t \in \mA| h_t,u_t)=Q(\mA| s_t,u_t) \qquad \forall \mA \in \mB(\mS), \ t >0$.
\end{enumerate}
Then the stochastic process $(\Omega,\mathcal{F},P_{\nu_0}^\pi,\{s_t\}_{t\geq 0})$ is called a {\it discrete-time Markov Decision Process}.
\end{definition}


\subsection{Infinite-Horizon Discounted-Reward Problem \label{sec:discounted}}

The following definition introduces $\alpha$-discount optimal policies, which maximize the discounted expected revenue along an infinite horizon.\\

\begin{definition}\label{def:infinite_horizon}
We define the {\it value function} under policy $\pi \in \Pi$ as the infinite-horizon discounted reward
\begin{equation*}
V(\pi,s_{0}):= \ \mathbb{E}_s^\pi \left[ \sum_{t=0}^\infty \alpha^t r(s_t,u_t) \right]\qquad \forall s_{0} \in \mathcal{S}, 
\end{equation*}
where $\mathbb{E}_s^\pi$ is the expectation taken with respect to the probability $P^\pi_{\nu_0}$, being $\nu_0$ concentrated at $s$. We say $\pi^*\in \Pi$ is an {\it $\alpha$-discount optimal policy}, with $\alpha \in [0,1)$, if it verifies
\begin{equation}\label{infinite_horizon}
V(\pi^*,s_{0}):=\sup_{\Pi} \ V(\pi,s_{0})=\sup_{\Pi}  \ \mathbb{E}_s^\pi \left[ \sum_{t=0}^\infty \alpha^t r(s_t,u_t) \right]\qquad \forall s_{0} \in \mathcal{S}. 
\end{equation}
We denote the {\it optimal value function} as $V^*(s_{0}):=V(\pi^*,s_{0})$.\\
\end{definition}

If the expectation in \eqref{infinite_horizon} is well defined, then $V(\cdot)$ is bounded for each $s \in \mS$. Therefore, from now on we restrict ourselves to the space $\mathcal{V}(\mS)$  of bounded measurable functions\footnote{\textcolor{black}{Most of the results in this paper hold in the larger space of functions that are bounded with respect to a weighted norm. However, we chose to work under the slightly stronger assumption that the functions are measurable and bounded in order to keep the exposition streamlined.}} in $\mS$. In general, we assume that the state $s_t \in \mathcal{S}$ is observed and then an action $u_t\in \mathcal{U}(s_t)$ is chosen. After this control is implemented, the probability that the MDP transitions to a state $s_{t+1}$ in the set $ \mathcal{A}\in \mathcal{B}(\mathcal{S})$ is given by $Q(\mathcal{A}|s_t,u_t)$.

\cite{lerma} also provides a more generalized formulation for \eqref{infinite_horizon} when the probability distribution $\nu_0$ of the initial state $s_0$ is given, i.e. 

\begin{equation}\label{infinite_horizon_nu}
\mathbb{V}^*(\nu_0): =\sup_{\Pi} \ \int_{\mS} V(\pi,s) d\nu_0=\sup_{\Pi}  \ \mathbb{E}_{\nu_0}^\pi \left[ \sum_{t=0}^\infty \alpha^t r(s_t,u_t) \right]\qquad \forall \nu_0 \in \mMpS.
\end{equation}
%

\noindent The definition of $\mathbb{V}^*(\cdot)$ is related to our measurized value function, as it evaluates state distributions rather than stochastic states themselves. Section \ref{sec:valuefunctions} outlines the relationship between $\mathbb{V}^*(\cdot)$ and $
\oV^*(\cdot)$. 

We now introduce some sufficient assumptions for the existence of an optimal deterministic stationary policy \cite[Chapter 4]{lerma}.

\begin{assumption}[\cite{lerma}, Assumption 4.2.1]\label{assumption421} \phantom{aaa}
\begin{enumerate}
\item[(a)] The one-stage reward $r$ is u.s.c., upper bounded, and sup-compact in $\mathbb{K}$
\item[(b)] Q is either
\begin{itemize}
\item[(b1)] weakly continuous
\item[(b2)] strongly continuous
\end{itemize}
\end{enumerate}
\end{assumption}

\begin{assumption}[\cite{lerma}, Assumption 4.2.2]\label{assumption422} 
There exists a policy $\pi$ such that $V(\pi,s)< \infty$ for each $s \in \mathcal{S}$.
\end{assumption}

For definitions of u.s.c. functions, and strongly and weakly continuous kernels, the reader is referred to Appendix A. The consequences of working in the weakly case instead of in strongly continuous case is that the value function belongs to the set of u.s.c. rather than just measurable functions (see Theorem 3.3.5 in \cite{lerma}).  In this paper, we adopt Assumptions \ref{assumption421}(a),(b2) and \ref{assumption422} for the original MDP. Proposition \ref{prop:assumptions} will show that then the lifted MDP inherits Assumptions \ref{assumption421}(a),(b1) and \ref{assumption422} \textcolor{black}{under additional compactness assumptions on the original state-action space}. The transition kernel of the measurized MDP is not strongly continuous because the transition is deterministic in the lifted space of probability measures (i.e., the measurized kernel is a Dirac measure).  More details can be seen in the proof of Proposition \ref{prop:assumptions}, in Section \ref{sec:assumptions}. 

\textcolor{black}{We denote by $\Pi^0$ the family of policies for which Assumption \ref{assumption422} holds. Note that $\Pi^0=\Pi$ when the reward function $r$ is bounded.} The following theorem shows that if Assumptions \ref{assumption421} and \ref{assumption422} hold, we can restrict ourselves to the set $\Pi^{D}$ of deterministic stationary policies without loss of optimality, and reformulate the infinite horizon problem \eqref{infinite_horizon} using the optimality equations \eqref{DCOE}.

\begin{theorem}[\cite{lerma}, Theorem 4.2.3]\label{theorem423} 
Suppose Assumptions \ref{assumption421} and \ref{assumption422} hold. Then 
\begin{itemize}
\item[(a)] The $\alpha$-discount value function $V^*(\cdot)$ is the pointwise maximal solution to \eqref{DCOE}, and if $V(\cdot)$ is another solution to \eqref{DCOE}, then $V(\cdot)\leq V^*(\cdot)$.
\item[(b)] There exists a selector $f^* \in \mathbb{F}$ such that $f^*(s) \in \mU(s)$ attains the maximum in \eqref{DCOE}; i.e., the deterministic stationary policy $f_\infty^*=\{f^*\}_{t\geq 0}$ is $\alpha$-discount optimal; conversely, if $f_\infty^* \in \Pi^{DS}$ is $\alpha$-discount optimal, then it solves \eqref{DCOE}.
\item[(c)] If $\lim_{t \to \infty} \alpha^t \mathbb{E}_s^\pi \left[V(\pi^*,s_t)\right]=0$ for all $\pi \in \Pi^0$ and $s \in \mS$, then $\pi^*$ is $\alpha$-discount optimal if and only if $V(\pi^*,\cdot)$ satisfies \eqref{DCOE}.
\item[(d)] If an $\alpha$-discount optimal policy exists, then there exists one that is deterministic stationary.
\end{itemize}

\end{theorem}

%

That is to say, the previous theorem shows that under certain assumptions the deterministic stationary policy $\pi^*=\{\varphi_t^*\}_{t \geq 0}$, with $\varphi^*_t(f^*(s)|s)=1$  for all $t=0,1,...$, is an $\alpha$-discount optimal solution to the traditional MDP \eqref{infinite_horizon}. 
For what follows we frame our distributional approach to MDPs from Markov stochastic kernels; i.e., we replace $\Pi$ by $\Pi^{MR}$ in \eqref{infinite_horizon_nu}. Moreover, Theorem \ref{theorem423} and Lemma 4.2.7 in \cite{lerma} entail that $V^* \in \mathcal{V}$ is the unique bounded solution to \eqref{DCOE} for a bounded reward function $r$ (see Note 4.2.1 in \cite{lerma} for more details on this). 



Finally, we might sometimes make the following assumption so that the Monotone Convergence Theorem (MCT) can be applied.

\begin{assumption} \label{assumption:MCT}
The reward function $r$ is bounded below.
\end{assumption}
Given Assumption \ref{assumption421} and the newly introduced requirement, it follows, without loss of generality, that the reward function can be taken to be non-negative\footnote{Assumption \ref{assumption421} ensures that there exists $c>0$ such that $r(s,u)\leq c$ for all $s \in \mS$, $u \in \mU$, and Assumption \ref{assumption:MCT} guarantees the existence of $k>0$ such that $r(s,u)\geq -k$ for all $s \in \mS$, $u \in \mU$. The reward function $r(s,u)+k$ is therefore upper bounded and non-negative for all $s\in \mS$, $u \in \mU$.}.


%

\subsection{The linear programming formulation \label{sec:LP}}
In this section, we introduce the LP formulation of the optimality equations \eqref{DCOE}.  The reader is referred to Chapter 6 in \cite{lerma} for details on how to derive this formulation and the theoretical results showing its equivalency with the \eqref{DCOE}. Let $\nu_0\in \mMpS$ be any distribution over the states; then the following linear program retrieves the optimal value function $\nu_0$-almost surely

\begin{align}
\inf_{V(\cdot)} & \ \int_{\mathcal{S}} V(s) d\nu_0(s)  \label{LP} \tag{LP}\\
\mbox{s.t.} & \ V(s)  \geq  r(s,u) + \alpha  \int_{\mS}  V(s') Q(ds'|s,u)   \qquad \qquad \forall  s \in \mathcal{S}, \forall  u \in \mathcal{U}(s).\nonumber
\end{align}
Denote $\mU(\mA)=\cup_{s \in \mA} \mU(s)$ as the set of feasible actions that can be taken from set $\mA \in\mBS$; then the dual of \eqref{LP} is

\begin{align}
\sup_{\zeta \in \mM_+(\mathbb{K})} & \ \int_{\mS} \int_{\mU(\mS)} r(s,u) \ d\zeta(s,u)  \label{D_zeta} \tag{D$_{\zeta}$}\\
\mbox{s.t.} &  \ \zeta(\mA\times\mU(\mA))=\nu_0(\mA)+\alpha \int_{\mS} \int_{\mU(\mA)} Q(\mA|s,u) d\zeta(s,u) & \forall \mA \in \mB(\mS) \nonumber
\end{align}
where measure $\zeta\in \mM_+(\mathbb{K})$ is the dual variable associated with the constraints in \eqref{LP}. This measure is defined on $\mathbb{K}$ because the constraints in the primal apply exclusively to feasible state-action pairs. Note that $\zeta$ is not a probability measure because $\zeta(\mS\times\mU(\mS))=1/(1-\alpha)$ for the constraint with $\mA=\mS$.  Indeed, $\zeta$ can be interpreted as an expected discounted state-action frequency. The next theorem summarizes some of the theoretical results in \cite{lerma}

\begin{theorem} \label{th:HLdual}
Under Assumptions \ref{assumption421} and \ref{assumption422}, the optimal solution $V^*(\cdot)$ to \eqref{LP} coincides with the optimal value function retrieved from the optimality equation \eqref{DCOE} $\nu_0$-almost surely. Moreover, \eqref{D_zeta} is solvable and strong duality holds: i.e., 

$$\int_\mS V^*(s) d\nu_0(s)=  \int_{\mS} \int_{\mU(\mS)} r(s,u) \ d\zeta^*(s,u),$$
where $\zeta^*$ is the optimal solutions to \eqref{D_zeta}.
\end{theorem}


\vspace{-0.2cm}
\section{ Simplification of the BS-deterministic framework \label{sec:measurized_def}}

In this section, we simplify the BS-deterministic framework of \cite{BSdeterministic, bertsekas} by showing two key results: (a) that without loss of generality, we can work with controls \(\varphi \in \mathcal{K}(\mS|\mU)\), the Markov decision rules from the original MDP, instead of controls \(\gamma \in \mM_{\mathbb{P}}(\mS\times\mU)\), which require the additional constraint that the marginal coincide with the current state distribution \(\nu\); and (b) that in the semicontinuous-semicompact framework, the lifted MDP inherits key assumptions from the original MDP, making it possible to apply the theory from \cite{lerma}. We refer to (a) as the {\it ``decomposition property"}, demonstrated in Section \ref{sec:change_variables}, and (b) as the {\it ``inheritance property"}, shown in Section \ref{sec:assumptions}.

We begin by introducing the BS-deterministic MDP framework as described in \cite{bertsekas}. Following this, we present our {\it measurized} MDPs, referred to as $\mM$-MDPs. Both frameworks involve dealing with MDPs defined on the space of probability measures that generalize standard MDPs. Hence, these MDPs represent specific cases within the broader class of measure-valued MDPs, where states are treated as measures. For a more comprehensive discussion on measure-valued MDPs, including formal definitions and examples such as Partially-Observable MDPs (POMDPs) and MFMDPs, we direct the reader to Appendix \ref{sec:measure-valued}.

\begin{definition}[BS-deterministic MDPs, Definition 9.4 in \cite{bertsekas}]\label{def:deterministicMDP}
Let $(\mathcal{S},\mathcal{U},\{\mathcal{U}(s)| \ s \in \mathcal{S}\},Q,r)$ be a standard MDP. The {\it corresponding deterministic} MDP \\
$(\mMpS,\mM_{\mathbb{P}}(\mS\times\mU),\{\Xi(\nu)| \ \nu \in \mMpS\},\overline{Q},\ovr)$, henceforth called {\it BS-deterministic} MDP, is a measure-valued process such that:
\begin{itemize}
\item[(i)] a state $\nu \in \mMpS$ is a probability measure over the states $s\in \mS$ of the standard MDP
\item[(ii)] an action $\gamma \in \mM_{\mathbb{P}}(\mS\times\mU)$ is a probability measure over the state-action pairs of the standard MDP
\item[(iii)] $\Xi(\nu)$ is the set of admissible actions from state $\nu$; more specifically 

\begin{equation}\label{Xi}
\Xi(\nu)=\left\{ \gamma \in \mathcal{M}_{\mathbb{P}}(\mS\times\mU): \gamma(\mathbb{K})=1 \mbox{ and }  \ \int_{s \in \mA} \int_{\mU(s)} \gamma(ds,du)=\nu(\mA) \ \ \forall \mA \in \mB(\mS) \right\},
\end{equation}

\item[(iv)] the transition kernel $\overline{Q}$ is deterministic. More specifically, the next state $\nu'$ of the measurized MDP is a probability measure computed according to a function $F:  \mM_{\mathbb{P}}(\mS\times\mU) \rightarrow \mMpS$, defined as

\begin{equation*}
\nu'(\cdot)=F(\gamma)(\cdot):=  \int_{\mathcal{S}}  \int_{\mathcal{U}}Q(\cdot|s,u) \gamma(du,ds)  
\end{equation*}
Therefore, $\oQ$ can be expressed using F as

\begin{equation*}
\oQ(\mathcal{P}|\gamma)=\left\{ 
\begin{array}{ll}
1 & \quad \mbox{ if } F(\gamma)\in \mathcal{P}\\
0 & \quad \mbox{ otherwise}\\
\end{array}
\right.
\end{equation*}
\item[(v)] the reward function $\overline{R}$ is the expected revenue of the standard MDP computed with respect to the joint distribution $\gamma \in \mathcal{M}_{\mathbb{P}}(\mS\times\mU)$; i.e.,

\begin{equation*}
\ovR(\gamma):=\mathbb{E}_{\gamma} [r(s,u)]= \int_{\mathcal{S}} \int_{\mathcal{U}}   r(s,u) \ \gamma(du,ds) .\\
\end{equation*}
\end{itemize}
\phantom{aa}\\
\end{definition}

In BS-deterministic MDPs, both the transition kernel $\oQ$ and the reward function $\ovR$ implicitly depend on the current state $\nu$, as the state information is embedded within the action $\gamma$. In a measurized MDP, the transition kernel $\oq$ and the reward function $\ovr$ depend directly on the state-action pair, as in standard MDPs.\\

\begin{definition}\label{def:measurizedMDP} 
Let $(\mathcal{S},\mathcal{U},\{\mathcal{U}(s)| \ s \in \mathcal{S}\},Q,r)$ be a standard MDP. A {\it measurized} MDP\\ $(\mMpS,\Phi,\{\Phi(\nu)| \ \nu \in \mMpS\},\overline{q},\ovr)$ is a measure-valued MDP, where (i) is as in Definition \ref{def:deterministicMDP} and
\begin{itemize}
\item[(ii)] \textcolor{black}{an action $\varphi \in \Phi$, where} $\Phi:=\{ \varphi \in \mathcal{K}(\mU|\mS): \ \varphi(\mU(s)|s)=1 \ \forall s \in \mS\}$, is a {\it feasible Markov decision rule in the standard MDP}
\item[(iii)] $\Phi(\nu)$ is the set of admissible actions from state $\nu$ such that $\Phi(\nu) \subseteq \Phi$, which {\it may or may not coincide with the set $\Phi$}
\item[(iv)] the transition kernel $\overline{q}$ is also deterministic; specifically, the next state $\nu'$ can be computed as
\begin{equation}\label{F}
\nu'(\cdot)=F(\nu,\varphi)(\cdot):=  \int_{\mathcal{S}}  \int_{\mathcal{U}}Q(\cdot|s,u) \varphi(du|s) d\nu(s). 
\end{equation}
\item[(v)] the reward function $\ovr$ is the expected revenue of the standard MDP computed 

both with respect to the state distribution $\nu \in \mMpS$ and the randomized action $\varphi \in \Phi$; i.e.,

\begin{equation*}
\ovr(\nu,\varphi):=\mathbb{E}_{\nu}\mathbb{E}_{\varphi} [r(s,u)]= \int_{\mathcal{S}} \int_{\mathcal{U}}   r(s,u) \ \varphi(du|s) d\nu(s).\\
\end{equation*}

\end{itemize}

Additionally, when $\Phi(\nu)=\Phi$, we say that the MDP is the {\it measurized counterpart} of the original MDP. If, on the other hand, $\Phi(\nu)\subset\Phi$, we say that the lifted MDP is a {\it tightened} measurized MDP.\\
\end{definition}

In Definition \ref{def:measurizedMDP}(iv), we slightly abuse notation by using the same symbol $F$ as in Definition \ref{def:deterministicMDP}(iv) to denote the deterministic transition function.

In the next section, we introduce a change of variables that transforms a BS-deterministic MDP into a measurized MDP.

\subsection{The decomposition property \label{sec:change_variables}} 

The aim of this section is to show that, for a given measurized state-action pair $(\nu_t,\varphi_t) \in \mMpS \times \Phi$ at any period of time $t$, there exists a unique deterministic action $\gamma_t \in \Xi(\nu_t)$ such that

\begin{equation}\label{change_variables}
\gamma_{t}(\mA\times\mathcal{D})=\int_{\mA} \int_{\mathcal{D}} \varphi_{t}(du|s) d\nu_{t}(s)
 \qquad \forall \mA \in\mB(\mS), \mathcal{D}\in \mB(\mU)
\end{equation}
\noindent holds almost surely with respect to $(\nu_t,\varphi_t)$.  \textcolor{black}{While similar results are presented in \cite{bertsekas} (see Section 7.4.3), their proofs do not utilize the Radon-Nikodym derivative and appear more complex. In addition, those results do not state the uniqueness of $\varphi_t$. Therefore, we offer our own derivation here. }


%

The following lemma ensures that, under certain conditions, there is a one-on-one relationship between the probability measure $\gamma \in \Xi(\nu)$ and the pair $(\nu,\varphi) \in \mMpS \times \Phi$. To ease notation, denote the complement of a set $\mA$ as $\mA^c$.

\begin{lemma}\label{lemma:change_variables}
Assume that $(\mS,\mBS,\nu)$ and $(\mU(s),\mathcal{B}(\mU(s)),\varphi(\cdot|s))$ for every $s \in \mS$ are complete measure spaces. For all $t$, let $\nu_t \in \mMpS$ and $\gamma_t \in \Xi(\nu_t)$ be probability measures. Then there exists a $\nu_t$-almost unique stochastic kernel $\varphi_t \in \Phi$ such that \eqref{change_variables} holds. 
\end{lemma}

\proof{

Pick any $\phi \in \Phi$. First, we prove that $\gamma_t(\mA\times\cdot)$ is absolutely continuous with respect to measure

$$\eta_t(\mA \times \mathcal{D})=\int_{\mA}  \phi(\mathcal{D}|s) d\nu_t(s) \qquad \mathcal{D}\in \mathcal{B}(\mU)$$

\noindent for all $ \mA \in\mB(\mS)$. We use the fact that $\eta_t(\mA\times\mathcal{D})=0$ implies that either

\begin{enumerate}
\item[(i)] $\mathcal{D}\bigcap \left( \cup_{s \in \mA} \mU(s)\right)=\emptyset$, or
\item[(ii)] $\nu_t(\mA)=0$
\end{enumerate}
In case (i), $\gamma_t(\mA\times\mathcal{D})=0$ by definition of $\Xi(\nu_t)$. In case (ii), completeness yields

\begin{equation}\label{gammaneq}
\gamma_t(\mA\times\mathcal{D}) \leq \gamma_t(\mA\times\mU)=\nu_t(\mA)
\end{equation}

\noindent and hence $\gamma_t(\mA\times\mathcal{D}) =0$ if $\nu_t(\mA)=0$. Then the Radon-Nikodym theorem states that for all $ \mA \in\mB(\mS)$ there exists a measurable function $g: \mS \times \mU \rightarrow [0,\infty)$ such that 
 
$$\gamma_t(\mA\times\mathcal{D})=\int_{\mA} \int_{\mathcal{D}}  g(s,u) \eta_t(ds,du)=\int_{\mA} \int_{\mathcal{D}}   g(s,u) \phi(du|s) d\nu_t(s).$$

\noindent and $g$ is unique $\eta_t$-almost everywhere. 

Second, we prove that the stochastic kernel defined as
\begin{equation}\label{kernel_g}
\varphi_t(\mathcal{D}|s)=\int_{\mathcal{D}}   g(s,u) \phi(du|s)  \qquad \forall s \in \mS, \forall \mathcal{D}\in \mB(\mU)
\end{equation}
belongs to $\Phi$. Note that, by definition of $\Xi(\nu_t)$, the marginal of $\gamma_t$ should coincide with $\nu_t$; i.e.,

\begin{equation}\label{marginal_Si}
\gamma_t(\mA\times\mU)=\nu_t(\mA), \quad \forall \mA \in \mB(\mS).
\end{equation}
This implies that

$$\gamma_t(\mA\times\mU)=\int_{\mU}  \int_{\mA} g(s,u) \phi(du|s) d\nu_t(s)=\nu_t(\mA).$$
\textcolor{black}{Since $g$ is a nonnegative measurable function, Tonelli’s theorem allows us to exchange the order of integration, yielding } 

$$\int_{\mU}   g(s,u) \phi(du|s) =1 \qquad \forall s \in \mS,$$
and because $\gamma_t \in \Xi(\nu_t)$ and the measure spaces are complete we also have that

$$\int_{\mathcal{D}}  \int_{\mA} g(s,u) \phi(du|s) d\nu_t(s)=0 \qquad \forall \mathcal{D}\subseteq \left(\cup_{s \in \mA} \mU(s)\right)^c$$
which proves that $\varphi_t \in \Phi$. Since $g$ is unique $\eta_t$-a.e., this means that there exists a stochastic kernel $\varphi_t \in \Phi$ that is unique $\eta_t$-a.e and is defined as in \eqref{kernel_g} such that \eqref{change_variables} holds.  \hfill $\blacksquare$\\

%

}

Denote $\Gamma=\{\gamma_t\}_{t \geq 0}$ and define the possible path of joint state-action distributions starting at $\nu_0$ as 
$$\Omega(\nu_0):=\left\{\Gamma: \right.  \ \gamma_{t} \in \mM_{\mathbb{P}}(\mS\times\mU) \ \forall t \geq 0,  \ \ \gamma_0\in \Xi(\nu_0), \left. \gamma_{t}(\cdot\times\mU)=\int_{\mS} \int_{\mU} Q(\cdot|s,u)d\gamma_{t-1}(s,u) \ \forall t \geq 1 \right\}.$$
The following proposition proves that given an initial state distribution $\nu_0$ and a sequence $\Gamma \in \Omega(\nu_0)$, there exist stochastic kernels $\varphi_t\in \Phi$ for all $t\geq 0$ such that there is a one-on-one relationship between $\gamma_{t}$ and $(\nu_{t},\varphi_{t})$ for all $t \geq0$, where $\nu_t$ obeys the transition law F defined in \eqref{F} for all $t \geq 1$.

\begin{proposition}\label{prop:change_variables_t}
\textcolor{black}{Assume that $(\mS,\mBS,\nu)$ and $(\mU(s),\mathcal{B}(\mU(s)),\varphi(\cdot|s))$ are complete measure spaces for every $s \in \mS$.} Given a $\nu_0$ and a sequence $\Gamma\in \Omega(\nu_0)$, there exists a unique sequence of probability measures $\{ \nu_t\}_{t \geq 0}$ and a $\nu_t$-unique sequence of stochastic kernels $\{ \varphi_t\}_{t \geq 0}$ verifying $\varphi_t \in \Phi$ and $\nu_{t+1}=F(\nu_t,\varphi_t)$ for all $t\geq 0$ such that \eqref{change_variables} holds. 
\end{proposition}

\proof{

We prove this by induction.

{\it t=0: } Since $\gamma_0 \in \Xi(\nu_0)$ and \eqref{marginal_Si} holds, we can apply Lemma \ref{lemma:change_variables} to prove that there exists a unique kernel $\varphi_0 \in \Phi$ such that \eqref{change_variables} holds. \\
 
{\it  t=T: } Assume true; i.e., assume \eqref{change_variables}  holds for  $\nu_T$.\\

{\it t=T+1: } Since $\Gamma\in \Omega(\nu_0)$, then $\gamma_{T} \in \Xi(\nu_T)$. From Lemma \ref{lemma:change_variables}, it suffices to prove that $\gamma_{T+1}(\mA,\mU)=\nu_{T+1}(\mA)$, where $\nu_{T+1}=F(\nu_{T},\varphi_{T})$, to have \eqref{change_variables}. Because $\Gamma\in \Omega(\nu_0)$, we have that

\begin{align*}
\gamma_{T+1}(\mA\times\mU)&=\int_{\mS} \int_{\mU} Q(\mA|s,u) d\gamma_T(s,u)\\
&=\int_{\mS} \int_{\mU} Q(\mA|s,u) \varphi_T(du|s) d\nu_T(s)=F(\nu_T,\varphi_T)(\mA)\\
&=\nu_{T+1}(\mA)
\end{align*}
where the second equality comes from the fact that $\gamma_T(\mA\times\mathcal{D})=\int_{\mA} \varphi_T(\mathcal{D}|s) d\nu_T(s)$ for all $\mA \in \mBS,\ \mathcal{D} \in \mBU$.\hfill $\blacksquare$\\

}
%


The previous result ensures that one can work with the BS-deterministic or the measurized counterpart of the original MDPs interchangeably. 

In addition, the decomposition property facilitates establishing a connection between these lifted MDPs and POMDPs, a topic that was not explored in \cite{bertsekas} nor \cite{BSdeterministic}. In a standard MDP, the agent has complete information about the current state of the environment, allowing it to make optimal decisions based on that information. In contrast, in a POMDP, the agent does not have direct access to the true state of the environment; instead, it observes partial, noisy, or incomplete information through observations \citep{lovejoy}. This lack of complete information introduces uncertainty and makes decision-making more challenging. It is well known that a POMDP can be modeled as a Belief-state MDP (BMDP). In a BMDP the agent keeps track of a belief state $\nu\in \mMpS$, which is the current prior distribution over states $s\in\mS$, and thus these can be framed within the more general measure-valued framework (see Appendix C). The state transitions to $\nu'\in \mMpS$ according to Bayes rules. From this perspective, the measurized MDP is structurally isomorphic to a fully observable MDP formulated as a BMDP. This connection situates our measurized MDPs within a well-established framework and emphasizes the importance of the measure-valued state representation.

\vspace{-0.1cm}
\subsection{The inheritance property \label{sec:assumptions}}

\textcolor{black}{In contrast to \cite{BSdeterministic, bertsekas}, who use lifted MDPs primarily to derive results for infinite-horizon stochastic MDPs, we show that if the existing theory in \cite{lerma} applies to the original MDP, it can also be extended to the lifted MDP.} More generally, the so-called inheritance property that we introduce in this section ensures that the theory developed in \cite{lerma} applies not only to the measurized counterpart of an MDP, where $\varphi \in \Phi$, but also for any tightened measurized counterpart if  $\Phi(\nu)$ is closed. 
 However, the following additional compactness assumption on subsets of the original feasible state-action space is needed to ensure that the measurized reward function $\overline{r}$ preserves the sup-compactness of $r$. In particular, we require that for every compact set of states, the corresponding set of feasible actions is compact.

\begin{assumption} \label{ass:compactness}
The set $\mathbb{K}_\mathcal{C}:=\{(s,u) \in \mathbb{K}: \ s \in \mathcal{C}\}$ is compact for all $\mathcal{C} \subseteq \mS$ compact.
\end{assumption}

The following result demonstrates that the measurized counterpart of an MDP preserves the \textcolor{black}{essence of the} assumptions adopted for its original MDP, albeit the transition kernel becomes weakly continuous.

\begin{proposition} \label{prop:assumptions}
Suppose that an MDP $(\mathcal{S},\mathcal{U},\{\mathcal{U}(s)| \ s \in \mathcal{S}\},Q,r)$ follows Assumptions \ref{assumption421}(a), (b2) and \ref{assumption422}, and that Assumptions \ref{assumption:MCT} and \ref{ass:compactness} also hold. Then  
its tightened measurized MDP $(\mMpS,\Phi,\{\Phi(\nu)| \ \nu \in \mMpS\},\overline{q},\ovr)$ follows Assumptions \ref{assumption421}(a), (b1) and \ref{assumption422} if the set of feasible actions $\Phi(\nu)$ is closed for every $\nu \in \mMpS$.
 \end{proposition}

\proof{\phantom{aa}

{\bf A1.(a).i}: We first show that $\ovr$ is upper bounded. This is straightforward since $r$ is upper bounded, so $\mathbb{E}_\nu \mathbb{E}_\varphi [r(s,u)]<\infty$ as well. 

{\bf A1.(a).ii}: Second, we prove that $\ovr$ is u.s.c. We will use Proposition E.2 in \cite{lerma} (corresponding to Proposition \ref{prop:E2} in Appendix A below). Take the joint probability measure $\gamma\in\mM_\mathbb{P}(\mS\times\mU)$ from \eqref{change_variables}. In Section \ref{sec:change_variables} we show that this change of variables can be performed without loss of generality; i.e., for every $\nu \in \mMpS$ and $\varphi \in \Phi(\nu)$ there is a unique $\gamma$ verifying \eqref{change_variables}. One can therefore construct a sequence $\{\gamma_n\}_{n\in\mathbb{N}}$ converging weakly to $\gamma$. Since $r$ is u.s.c. by assumption, by Proposition E.2 in \cite{lerma} we get

$$ \lim\sup_{n \to \infty} \int_\mS  \int_\mU r(s,u) d\gamma_n(s,u) \leq  \int_\mS  \int_\mU r(s,u) d\gamma(s,u).$$
Because any pair $(\nu,\varphi)$ is uniquely characterized by a $\gamma$, this shows that $\ovr$ is upper semicontinuous.

{\bf A1.(a).iii}: To finish with Assumption \ref{assumption421}(a), we need to prove that $\ovr$ is sup-compact in  $\overline{\mathbb{K}}$. That is to say, we need to show that the set

$$\omU_c(\nu):=\{ \varphi \in \Phi(\nu): \ \ovr(\nu,\varphi) \geq c\}$$ 
is compact for every $\nu \in \mMpS$ and $c \in \mathbb{R}$. 

To simplify notation, we redefine this set in terms of joint state-action distributions:
$\overline{\mathcal U}_c(\nu):= \left\{ \gamma \in \tilde{\Xi}(\nu) : \overline{r}(\nu,\varphi) \ge c \right\}$,
where $\gamma \in \tilde{\Xi}(\nu)$ if $\gamma \in \Xi(\nu)$ (see \eqref{Xi}) and there exists $\varphi \in \Phi(\nu)$ such that \eqref{change_variables} holds. Thus, $\tilde{\Xi}(\nu)$ is the subset of $\Xi(\nu)$ consisting of joint state-action distributions $\gamma$ that admit a representation in terms of feasible measurized actions $\varphi$ in the tightened MDP (i.e., under $\Phi(\nu)$ rather than $\Phi$).

In infinite-dimensional spaces, compactness does not necessarily coincide with closedness and boundedness. We therefore work with relative compactness. Recall that a subset of a topological space is relatively compact if its closure is compact (see \cite{hitchhiker}, Section 2.8). Hence, it suffices to show that $\omU_c(\nu)$ is relatively compact and closed. For this proof, we endow $\mMp(\mS \times \mU)$ with the weak topology.


Since $\mS$ and $\mU$ are Borel spaces, they are separable and complete metric spaces (see Remark 2.1.1 in \cite{lerma}). As a consequenc, each $\nu \in \mMpS$ is {\it tight} \cite[Theorem 1.3]{billingsley}; i.e., for every $\epsilon>0$ there exists a compact set $\mathcal{C}_\epsilon \subset \mS$ such that $\nu(\mathcal{C}_\epsilon)>1-\epsilon$. By assumption, any set of the form $\mathbb{K}_{\mathcal{C}_\epsilon}$ is compact. In addition, we have that

$$\gamma(\mathbb{K}_{\mathcal{C}_\epsilon})=\int_{s \in \mathcal{C}_\epsilon} \int_{u \in \mU(s)} d\gamma(s,u)=\nu(\mathcal{C}_\epsilon) >1-\epsilon \quad \forall \gamma\in\Xi(\nu). $$

\noindent Therefore, all $\gamma \in \Xi(\nu)$ are tight and, as a consequence, $\Xi(\nu)$ is tight. Because $\mU$ is also a separable and complete metric space, Prohorov's Theorem (see Theorem E.6 in \citet{lerma}) ensures that  $\Xi(\nu)$ is relatively compact in $\mMp(\mS\times\mU)$ with the weak topology. Since $\omU_c(\nu) \subseteq \tilde{\Xi}(\nu) \subseteq \Xi(\nu)$, then $\omU_c(\nu)$ is relatively compact.


%

To show that $\omU_c(\nu)$ is closed, we build on the fact that the function $\varphi \mapsto \int_\mS  r(s,u) \varphi(\mU|s) d\nu$ is upper semicontinuous on $\Phi(\nu)$ for every $\nu \in \mMpS$. To see this, follow the proof of A1.(a).ii fixing $\nu$. According to Proposition A.1 in \cite{lerma} (see Proposition \ref{prop:A1} in Appendix A), this implies that $\omU_c(\nu)$ is closed. We have this result for every $\nu \in \mMpS$.

{\bf A1.(b)}: We now show that $\oq$ inherits weak continuity from the strong continuity of $Q$. We need to prove that

$$(\nu,\varphi) \mapsto \int g(\mu) \oq(d\mu| \nu,\varphi) $$
is continuous and bounded for every function $g$ continuous and bounded in $\mMpS$ (see Definition C.3 in \cite{lerma}, which can be found in the Appendix A as Definition \ref{def:stronglycontinuous} for convenience). Because $\oq$ is concentrated at $F(\nu,\varphi)$, we get that 

$$\int g(\mu) \oq(d\mu| \nu,\varphi)= g(F(\nu,\varphi)).$$
It suffices to show that $(\nu,\varphi) \mapsto F(\nu,\varphi)$ is continuous because the composition of continuous functions preserves continuity. Consider a sequence $\{(\nu_n,\varphi_n)\}_{n\in \mathbb{N}}$ that converges weakly to $(\nu,\varphi)$. The function $Q(\mA|\cdot,\cdot)$ is continuous for all $\mA \in \mBS$ because $Q$ is strongly continuous (see Proposition C.4 in \cite{lerma}). Therefore, by definition of weak convergence of measures (see Definition E.1 in \cite{lerma}) we have that

$$F(\nu_n,\varphi_n)(\cdot)=\int_\mS \int_\mU Q(\cdot|s,u) \varphi_n(du|s) d\nu_n(s) \rightarrow \int_\mS \int_\mU Q(\cdot|s,u) \varphi(du|s) d\nu(s)=F(\nu,\varphi)(\cdot),$$
yielding continuity.

{\bf A2}: We prove that exists a policy $\overline{\pi}$ such that $\oV(\overline{\pi},\nu)<\infty$ for each $\nu \in \mMpS$. 

It suffices to plug into \eqref{infinite_horizon_nu} the policy $\overline{\pi}_{f^*}=\{\varphi_{f^*}\}_t$, where $\varphi_{f^*}(f^*(s)|s)=1$ for all $s \in \mS$ and $f^*$ is the optimal selector of the standard MDP (see Theorem \ref{theorem423}(b)). This yields $\oV(\overline{\pi}_{f^*},\nu)=\int_\mS V^*(s) d\nu(s)<\infty$.

Let $\pi^*=\{f^*\}_{t \geq 0}$ be an optimal deterministic policy for the standard MDP, where $f^*$ is the optimal selector of the standard MDP (see Theorem \ref{theorem423}(b)). We evaluate the measurized DIH value function under the measurized policy $\overline{\pi}_{f^*}=\{\varphi_{f^*}\}_{t\geq 0}$, where $\varphi_{f^*}(f^*(s)|s)=1$ for all $s \in \mS$. Then we have

\begin{align*}
\oV(\overline{\pi}_{f^*},\nu) & = \mathbb{E}^{\overline{\pi}_{f^*}}_\nu \left[\sum_{t=0}^\infty \alpha^t \ovr(\nu_t,\varphi_t) \right]\\
&=\sum_{t=0}^\infty \alpha^t \ovr(\nu_t,\varphi_{f^*}) \\
&=\sum_{t=0}^\infty \alpha^t \mathbb{E}_{\nu_t} \mathbb{E}_{\varphi_{f^*}} [r(s,u)] \\
&= \mathbb{E}_{\nu} ^{\pi^*} \left[ \sum_{t=0}^\infty \alpha^t r(s,u)\right] \\
&=\int_\mS V(\pi^*,s) d\nu_0 <\infty
\end{align*}

\noindent where the first equality comes from Definition \ref{def:infinite_horizon}, the second from the fact that the transitions are deterministic in the measurized MDP, the third from the definition of the measurized reward function, and the fourth from the MCT.  \hfill $\blacksquare$ 

}

\textcolor{black}{In the following sections, we demonstrate how the previous result allows us to directly state the validity of the measurized equations \eqref{M-DCOE-simple}, which was achieved through a much more laborious process in Section 9.4 of \cite{bertsekas}. Additionally, we also show that one can retrieve the optimal measurized value function by solving the LP formulation of the measurized MDP, something that was not done in \cite{bertsekas}.}


\section{Advantages of working in the lifted framework \label{sec:advantages}}

In this section, we illustrate the benefits of working within the lifted framework. In Section \ref{sec:OE}, we derive the optimality equations for the measurized MDP, resulting in deterministic equations that can be addressed using deterministic optimization techniques. Sections \ref{sec:constraints} and \ref{sec:approximations} delve into how the augmented space allows for the incorporation of constraints and value function approximations that are challenging to model within the stochastic framework. Although these constraints and approximations could also be applied to BS-deterministic MDPs, to the best of our knowledge, they have not been explored in this context before. \textcolor{black}{Throughout this section, we assume that the conditions of Proposition \ref{prop:assumptions} hold.}



\subsection{Deterministic optimality equations \label{sec:OE}} 

Since measurized MDPs are a special class of MDPs, the optimal measurized value function can be defined as usual

\begin{equation} \label{measurized_infinitehorizon}
\oV^*(\nu_0):=\sup_{\overline{\pi} \in \overline{\Pi}} \ \mathbb{E}_{\nu_0}^{\overline{\pi}} \left[ \sum_{t=0}^\infty \alpha^t \ovr(\nu_t,\varphi_t)\right],
\end{equation} 
where $\overline{\Pi}$ is the set of all history-dependent randomized policies for the measurized MDP.

A measurized Markov decision rule $\psi$ is a stochastic kernel defined on the lifted space of measures; i.e. $\psi(\cdot|\nu)$ is a probability measure over $\mathcal{K}(\mU|\mS)$ for all $\nu$, and  $\psi(\varphi|\cdot)$ is a measurable function on $\mMpS$ for every $\varphi$.
Proposition \ref{prop:assumptions} ensures that Theorem \ref{theorem423} applies and we can restrict ourselves to the set of deterministic stationary policies $\overline{\pi}=\{\psi\}_{t \geq 0}$, where $\psi$ is concentrated around a selector $\overline{f}$. In other words, the supremum in the infinite-horizon problem above is attainable and there exists a selector $\overline{f}^*: \mMpS \to \Phi$ such that $\of^*(\nu) \in \Phi(\nu)$ for all $\nu \in \mMpS$ and $\of^*$ is the solution to \eqref{measurized_infinitehorizon}; i.e. we can replace $\overline{\Pi}$ by $\overline{\Pi}^D:=\{\psi_t \in \overline{\Pi}: \ \psi_t=\psi \ \forall t\geq 0 \mbox{ and } \forall \nu \in \mMS, \ \exists \varphi \in \Phi(\nu) \mbox{ s.t. } \psi(\varphi|\nu)=1\}$ in \eqref{measurized_infinitehorizon}. Since the transition to the next state is deterministic and the MCT allows us to interchange the expectation and the infinite sum, we get that 

\begin{equation} \label{oV*} 
\oV^*(\nu_0)=\sup_{\substack{\varphi_t \in \Phi(\nu_t) \\ \forall t \geq 0}} \left\{ \sum_{t=0}^\infty \alpha^t \ovr(\nu_t,\varphi_t) \right\},
\end{equation} 
where $\nu_{t}=F(\nu_{t-1},\varphi_{t-1})$ for all $t \geq 1$. Finally, Theorem \ref{theorem423} ensures that one can retrieve the optimal measurized value function $\oV^*(\cdot)$ and the optimal measurized actions $\of^*(\nu)$ from 

\begin{align}
\overline{V}^*(\nu)&=\sup_{\varphi \in \Phi(\nu)} \ \left\{ \ovr(\nu,\varphi) + \alpha   \overline{V}^*(F(\nu,\varphi))\right\} \qquad \forall \nu \in \mMpS. \label{measurized_DCOE} \tag{$\mM$-$\alpha$-DCOE}
\end{align}
%
\cite{bauerle} obtained similar deterministic equations for MFMDPs without common noise. More details on \cite{bauerle} and how it relates to the measurized MDPs can be found in Appendices \ref{sec:measure-valued} and \ref{appendix:bauerle}.
Since the measurized MDP inherits Assumption \ref{assumption421} and \ref{assumption422} from its original MDP, Theorem \ref{th:HLdual} also applies, yielding the following LP formulation that is essentially equivalent to \eqref{measurized_DCOE}

\begin{align}
\inf_{\oV(\cdot)} & \ \oV(\nu_0)  \label{measurized_LP} \tag{$\mM$-LP}\\
\mbox{s.t.} & \ \oV(\nu)  \geq  \ovr(\nu,\varphi) + \alpha \oV(F(\nu,\varphi))    \qquad \qquad \forall  \nu \in \mMpS, \forall  \varphi \in \Phi(\nu).\nonumber
\end{align}

One can think of the measurized MDP as aiming to optimize the distribution of the states through controls that are stochastic kernels over the original actions $u$. For any initial state distribution $\nu_0$ and any measurized policy $\overline{\pi}=\{\varphi_t\}_{t \geq 0}$ the trajectory of state distributions $\{F(\nu_t,\varphi_t)\}_{t\geq 1}$ can be computed deterministically using \eqref{F}, thus giving rise to a deterministic process. In particular, an optimal policy $\overline{\pi}^*:=\{\varphi_t^*\}_{t \geq 0}$ to the measurized MDP gives rise to an optimal trajectory of state distributions

\begin{equation}\label{trajectory}
\Upsilon^*:=\{\nu_t^*\}_{t\geq 0},
\end{equation}
where $\nu_t^*:=F(\nu_{t-1}^*,\varphi_{t-1}^*)$ for $t\geq 1$. In addition, the measurized framework facilitates the understanding of which actions we expect to take most frequently under a certain distribution of the states. More specifically, if at period $t$ the states are distributed according to $\nu_t$, we can define the distribution of actions under decision rule $\varphi_t$ as

\begin{equation}\label{rho}
\rho_{t}(\mathcal{D}):=\int_{\mS} \varphi_{t}(\mathcal{D}|s) d\nu_{t}(s).
 \qquad \forall \mathcal{D}\in \mB(\mU).
\end{equation}
Proposition \ref{prop:change_variables_t} ensures that $\rho_t$ is well defined, since for every pair $\varphi_t, \ \nu_t$ there exists a unique $\rho_t$ and vice-versa.


%
\subsection{Probabilistic constraints \label{sec:constraints}} 
Proposition \ref{prop:change_variables_t} implies that the BS-deterministic MDP associated with a standard MDP is equivalent to the measurized counterpart of that MDP; that is, the measurized MDP with $\Phi(\nu)=\Phi$. Interestingly, the set of admissible actions $\Phi(\nu)$ is allowed to be contingent on the current state distribution $\nu$\footnote{Although $\Xi(\nu)$ depends on $\nu$, this dependency arises solely from the requirement that the marginal of the action $\gamma$ matches $\nu$. However, to the best of our knowledge, no additional constraints in the augmented space have been explored until now.}.
This enables the modelling of probabilistic constraints on the states of the original MDP. For instance, one can set constraints on various risk measures related to the agent's future perceived costs or limit moments of future distributions over original states. Furthermore, one could also impose restrictions on the distribution of actions taken in the original MDP. These constraints do not arise naturally outside the lifted framework.

The following example illustrates how one could easily introduce Conditional Value at Risk (CVaR) constraints in the measurized framework. Handling such risk constraints is notably challenging in the conventional framework, leading much of the literature to focus on finite horizon MDPs and bounding the CVaR of discounted accumulated costs at a terminal stage, as seen in works like \cite{borkar}. A more recent work \cite{glynn} considers a long-run average MDP with finite state and action spaces, and optimizes the CVaR at the steady state. As the authors mention {\it ``dynamically optimizing CVaR is difficult since it is not a standard MDP and the principle of dynamic programming fails"}. In contrast, we offer here a more straightforward alternative within our framework. \\

\begin{example}[MDPs with CVaR constraints \label{ex:CVaR}] 
Consider a standard MDP with cost function $c: \mS \times \mU \rightarrow \mathbb{R}$. The optimality equations of its measurized counterpart with no additional constraints are

\begin{align*} 
\overline{V}^*(\nu)&=\inf_{\varphi \in \Phi} \ \left\{ \int_\mS \int_\mU c(s,u) \varphi(du|s) d\nu(s) + \alpha   \overline{V}^*(F(\nu,\varphi))\right\},  \qquad \forall \nu \in \mMpS.
\end{align*}
Here $\nu$ is the distribution of states of the original MDP, to which we want to add a CVaR constraint, and $\nu'(\cdot)=F(\nu,\varphi)(\cdot)$ is the subsequent state distribution. Define the Value at Risk (VaR$_\beta$) at the current period as

$$VaR_\beta(c;\nu,\varphi):=\arg\inf_{a \in \mathbb{R}} \left\{ \mathbb{P}_{\nu,\varphi}(c(s,u)\leq a)\geq \beta \right\}. $$
This measure quantifies the potential financial loss or risk within a specified confidence level $\beta>0$. 

We now explore how to write and interpret the VaR within an MDP. For simplicity, consider $\mS\subseteq \mathbb{R}$ and $\mU\subseteq \mathbb{R}$. For every $s \in \mS$, $a \in \mathbb{R}$, denote $\mathcal{U}_{s,a}:=\{ u\in \mU(s): c(s,u)\leq a\}$. Given $\nu$ and $\varphi$, we can express $VaR_\beta$ as 
\begin{align*}
VaR_\beta(c;\nu,\varphi)&= \arg\inf_{a \in \mathbb{R}} \left\{ \mathbb{P}_{\nu,\varphi}(c(s,u)\leq a)\geq \beta \right\} \\
&= \arg\inf_{a \in \mathbb{R}} \left\{ \int_{s \in\mS} \int_{u \in  \mathcal{U}_{s,a}}  \varphi(du|s) d\nu(s) \geq \beta \right\} \\
&= \arg\inf_{a \in \mathbb{R}} \left\{ \int_{s \in\mS}  \varphi(\mathcal{U}_{s,a}|s) d\nu(s) \geq \beta \right\} ,\\
&= \arg\inf_{a \in \mathbb{R}} \left\{   \rho(\mathcal{U}_{s,a})  \geq \beta \right\} ,
\end{align*}
where $\rho$ is the distribution of actions \eqref{rho}, which depends on $\nu$ and $\varphi$. If we think of the set $\mathcal{U}_{s,a}$ as the set of {\it ``safe actions"} to take from state $s$, we can see the $VaR_\beta$ as the minimum value $a^*$ such that the actions taken {\it ``are safe"} with a probability larger or equal to $\beta$. Increasing $\beta$ amounts to expanding that set so that the decisions taken are {\it``more often"} therein. In other words, the larger the $\beta$, the larger the set of actions we consider as {\it ``safe"}.

We now focus on the conditional VaR, which we define using the Acerbi's formula, to which we plug the VaR computed with respect to probability $\mathbb{P}_{\nu,\varphi}$

\begin{align}\label{CVaR}
CVaR_\upsilon(c)&:=\frac{1}{1-\upsilon} \int_\upsilon^1 VaR_\beta(c) d\beta\\
&=\frac{1}{1-\upsilon} \int_\upsilon^1 \arg\inf_{a \in \mathbb{R}} \left\{ \rho(\mathcal{U}_{s,a}) \geq \beta \right\} d\beta  ,\nonumber
\end{align}
In the equation above, threshold $\beta\in (0,1)$ is not given. Instead, we integrate with respect to $\beta$ for a ``reasonable" interval of values (e.g. if $\upsilon=0.95$ we consider $\beta\in[0.95,1)$). This provides an {\it ``average VaR"}, i.e., the expected threshold $a^*$ such that the actions taken {\it ``are safe"} within a reasonable range of probabilities.

Assume now that we want to bound the CVaR so it does not exceed a threshold $\theta>0$ at any decision period of the MDP. Integrating this requirement into the optimality equations yields the CVaR-constrained MDP

\begin{align*} 
\overline{V}^*(\nu)=\inf_{\varphi \in \Phi}  &\ \left\{ \int_\mS \int_\mU c(s,u) \varphi(du|s) d\nu(s) + \alpha   \overline{V}^*(F(\nu,\varphi))\right\},  \hspace{2.5cm} \forall \nu \in \mMpS\\
 \mbox{s.t.} & \ \frac{1}{1-\upsilon} \int_\upsilon^1 \arg\inf_{a \in \mathbb{R}} \left\{\rho(\mathcal{U}_{s,a})  \geq \beta \right\} d\beta \leq \theta  \hspace{5.6cm} \blacksquare
\end{align*}


%



%
\end{example}

\subsection{Probabilistic value function approximations \label{sec:approximations}}
In the usual MDP framework, value function approximations are often modelled as the weighted sum of basis functions $\phi_k: \mS \to \mathbb{R}$, yielding 

\begin{equation} \label{V_approx}
V(s)\approx \sum_{k=1}^K w_k \phi_k(s).
\end{equation}
Some choices for $\phi_k(\cdot)$ that have been explored in the literature are linear \cite{adelman2007}, separable piecewise linear \cite{vossen} or ridge exponential \cite{IJOC}. The theoretical results developed in the following section (and, more specifically, Theorem \ref{th:V}.(b)) ensures that imposing an approximation \eqref{V_approx} in the original MDP translates into an approximation in expected value; i.e., 

\begin{equation} \label{oV_Eapprox}
\oV(\nu)\approx \sum_{k=1}^K \overline{w}_k \mathbb{E}_\nu[\phi_k(s)].
\end{equation}
Interestingly, considering a measurized MDP whose states are probability measures in state-space $\mS$ allows one to model a broader class of basis functions that are contingent on $\nu$, yielding

\begin{equation} \label{oV_approx}
\oV(\nu)\approx \sum_{k=1}^K \overline{w}_k \overline{\phi}_k(\nu).
\end{equation}
Certain basis functions $\overline{\phi}_k(\cdot)$ may not be accessible within the original framework. To illustrate this, the subsequent example introduces some compelling measure-valued approximations.\\

\begin{example}[Some augmented basis functions] \label{ex:basis}
There are multiple measurized basis functions $\overline{\phi}$ that can be written as the expected value of some original basis function $\phi$. For instance, if $a>0$ is a scalar and the state space is one-dimensional, the following measurized basis functions can be considered in the augmented space:
\begin{itemize}
\item[(a)] Basis functions like $ \overline{\phi}_k(\nu)=\int_\mS s^a d\nu(s)=\mathbb{E}_\nu[s^a]$ provide a moment approximation, generated by $\phi_k(s)=s^a$.
\item[(b)] One could consider Laplace transforms as basis functions; i.e., $ \overline{\phi}_k(\nu)=\int_\mS e^{as} d\nu(s)=\mathbb{E}_\nu[e^{as}]$, generated by $\phi_k(s)=e^{as}$.
\item[(c)] In addition, basis functions could be chosen as generating functions of the state $s$, having $ \overline{\phi}_k(\nu)=\int_\mS a^{s} d\nu(s)=\mathbb{E}_\nu[a^{s}]$, when $\phi_k(s)=a^s$.
\item[(d)] Furthermore, consider a state $s$ taking values in $\mathbb{R}$. For instance, $s$ might measure the risk of cardiovascular disease in a patient, the losses incurred by a portfolio... etc. One could include probabilistic basis functions like $ \overline{\phi}_k(\nu)=\mathbb{P}_\nu(s>a)$, with $a\in[0,100]$, to approximate the value of the current risk distribution $\oV^*(\nu)$. Here $\phi_k(s)=\mathbf{1}_{\{s>a\}}$, where $\mathbf{1}_{\{\cdot \}}$ is the indicator function. 
\end{itemize}
However, many other basis functions in the lifted space are different in nature, in the sense that they cannot be written merely as expected values of some function $\phi_k(s)$ of the random state $s$. For example:
\begin{itemize}
\item[(e)] One could use the CVaR \eqref{CVaR} as a basis function. Diverse basis functions arise for different values of the threshold $\upsilon$.
\item[(f)] Basis functions could measure the distance to a benchmark distribution $\mu$; for instance:
\begin{itemize}
\item[(f.1)] $\overline{\phi}(\nu)$ can be the Wasserstein $p$-distance between $\nu$ and $\mu$
$$\overline{\phi}(\nu)=\left( \inf_{\eta \in \Delta(\nu,\mu)} \int_{\mS\times \mS} \|s-x \|_p d\eta(s,x) \right)^{1/p}$$
 where $p\geq 1$ and $\Delta(\nu,\mu)$ is the set of all joint distributions $\eta$ in $\mS \times \mS$ with marginals $\nu$ and $\mu$. Note that different values of $p$ yield different basis functions. 
 \item[(f.2)] Instead, one could consider the Kullback–Leibler divergence as a basis function
 $$\overline{\phi}(\nu)= \int_\mS log\left( \frac{d\nu}{d\mu}\right) d\nu(s),$$
 where $d\nu/d\mu$ is the Radon-Nikodym derivative of $\nu$ with respect to $\mu$.
\end{itemize}
\end{itemize}

Note that (a)-(d) are available in the standard Approximate Dynamic Programming context, because they derive from standard basis functions. They become moment approximations in the measurized framework but do not provide any additional approximation power. In contrast, our framework allows (e)-(f).
\phantom{aa} \hfill $\blacksquare$\\
\end{example}

%

\section{Connection between the stochastic MDP and its measurized counterpart \label{sec:measurized}} 

In this section we analyze the connection between the original MDP and its measurized counterpart; i.e., the lifted MDP with $\Phi(\nu)=\Phi$ for all $\nu \in \mMpS$.




\subsection{Relationship between states. \label{sec:states}}

We first examine the relationship between the measurized and original states, summarized in Figure \ref{fig:measurized_states}.
\begin{center}
\vspace{-0.5cm}
\begin{figure}[H]
\centering
\begin{tabular}{c}
\includegraphics[width=250pt]{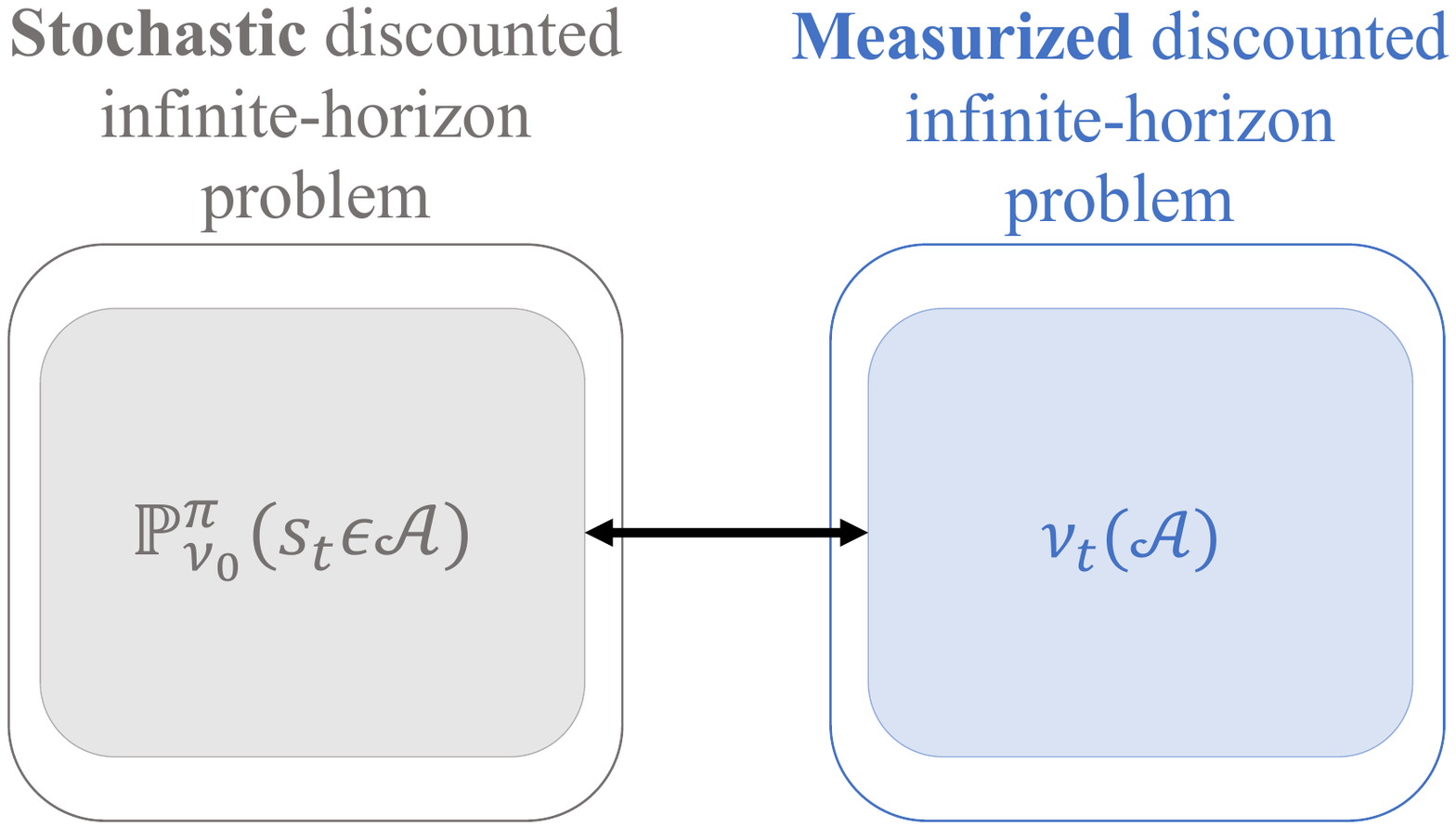}\\
\end{tabular}
\vspace{-0.5cm}
\caption{ Equivalence between the states $\nu$ of the measurized MDP and the states $s$ of the original MDP from which it was lifted. \label{fig:measurized_states}}
\end{figure}
\vspace{-0.5cm}
\end{center}

Let $\pi=\{\varphi_t\}_{t\geq 0}$; Theorem \ref{theorem423} enables one to consider exclusively either deterministic or Markov policies without loss of optimality. Then using the definition of the revenue function $\ovr$, we can rewrite the infinite horizon problem \eqref{infinite_horizon_nu} as the discounted measurized MDP problem as follows

\begin{align}
\mathbb{V}^*(\nu_0) &=\sup_{\pi \in \Pi^{MR}}\mathbb{E}_{\nu_0}^\pi \left[ \sum_{t=0}^\infty \alpha^t r(s_t,u_t)\right] &\nu_0 \in \mMS\nonumber\\
 &=\sup_{\substack{\varphi_t \in \Phi, \\ t\geq 0}}  \   \sum_{t=0}^\infty  \alpha^t  \mathbb{E}_{\nu_t} \mathbb{E}_{\varphi_t} [r(s,  u)] & \nu_0 \in \mMS \label{M_infinite}\\
  &=\sup_{\substack{\varphi_t \in \Phi, \\ t\geq 0}}  \   \sum_{t=0}^\infty  \alpha^t  \ovr(\nu_t,  \varphi_t) &\nu_0 \in \mMS \nonumber \\
  &=\oV^*(\nu_0) \nonumber
\end{align}
where the second equality arises from the MCT (which needs Assumptions \ref{assumption421} and \ref{assumption:MCT}) and the last equality comes from \eqref{oV*}. This proves that the optimal value function $\mathbb{V}^*(\cdot)$ defined in \eqref{infinite_horizon_nu} coincides with the optimal measurized value function $\oV^*(\cdot)$ at every state distribution $\nu \in \mMS$. Through \eqref{M_infinite}, one can intuit a connection between the distributions $\nu_t$ controlled through $\varphi_t$  and the probability distribution $P_{\nu_0}^{\pi}$. We explicitly state this relationship in the following proposition.

\begin{proposition}\label{prop:distribution}
\textcolor{black}{Suppose Assumptions \ref{assumption421} and \ref{assumption422} hold and} define the policy $\pi=\{\varphi_t\}_{t \geq 0}$. Then the sample path $\Upsilon=\{\nu_t\}_{t\geq 0}$, with $\nu_{t+1}=F(\nu_t,\varphi_{t})$  for all $t\geq 0$, is related to the probability distribution $P_{\nu_0}^{\pi}$ as follows
\begin{equation}\label{Pnu}
P_{\nu_0}^{\pi}(s_{t+1} \in \mA)=\nu_{t+1}(\mA) \qquad \forall \mA \in \mB(\mS)
\end{equation}
\end{proposition}

\proof{

According to Remark C.11 in \cite{lerma}, the probability measure $P_{\nu_0}^{\pi}$ has the following expression for all $t \geq 0$

\begin{align}\label{Pnupi}
P_{\nu_0}^{\pi}(s_{t+1} \in\mA_{t+1})&=\medmath{\int_{\mS} \nu_0(ds_0) \int_{\mU} \int_{\mS}Q(ds_1|s_0,u_0) \varphi_0(du_0|s_0) \int_{\mU} \int_{\mS}Q(ds_2|s_1,u_1) \varphi_1(du_1|s_1)  \hdots \int_{\mU} Q(\mA_{t+1} |s_t,u_t) \varphi_t(du_t|s_t) }
\end{align}

We prove this proposition by induction

{\it Prove for t=0:}  \quad $P_{\nu_0}^{\pi}(s_0 \in \mA_0)=\int_{\mA_0} \nu_0(ds_0)=\nu_0(\mA_0).$

{\it Prove for t=1:} \quad \begin{align*}
P_{\nu_0}^{\pi}(s_1 \in \mA_1)&=\int_{\mS} \int_{\mU} \int_{\mA_1} Q(ds_1|s_0,u_0) \varphi_0(du_0|s_0) \nu_0(ds_0)\\
&=\int_{\mS} \int_{\mU}  Q(\mA_1|s_0,u_0) \varphi_0(du_0|s_0) \nu_0(ds_0)\\
&=F(\nu_0,\varphi_0)(\mA_1)\\
&=\nu_1(\mA_1)
\end{align*}
where the last equality comes from \eqref{F}.\\

{\it Assume true for t:}  \quad $P_{\nu_0}^{\pi}(s_t \in \mA_t)=\nu_t(\mA_t).$

{\it Prove for t+1:} \begin{align*}
\nu_{t+1}(\mA_{t+1})&=  \int_{\mathcal{S}}  \int_{\mathcal{U}}Q(\mA_{t+1}|s_t,u_t) \varphi_t(du_t|s_t) d\nu_t(s_t)\\
&\medmath{=\int_{\mS} \nu_0(ds_0) \int_{\mU} \int_{\mS}Q(ds_1|s_0,u_0) \varphi_0(du_0|s_0) \int_{\mU} \int_{\mS}Q(ds_2|s_1,u_1) \varphi_1(du_1|s_1)  \hdots \int_{\mU} Q(\mA_{t+1} |s_t,u_t) \varphi_t(du_t|s_t) }\\
&=P_{\nu_0}^{\pi}(s_{t+1} \in \mA_{t+1})
\end{align*}
where the second equality comes from the fact that $\nu_t(\mA_t)=P_{\nu_0}^\pi(s_t \in \mA_t)$ and hence we can plug the expression of $P_{\nu_0}^\pi$ as defined in \eqref{Pnupi}. \hfill $\blacksquare$\\

}

\subsection{Relationship between optimal value functions. \label{sec:valuefunctions}}
In this section, we provide various theoretical results that link the measurized value function $\oV^*(\cdot)$ with the stochastic value function $V^*(\cdot)$ of the original MDP from which it was lifted. Figure \ref{fig:measurized_valuefunction} summarizes the results.

\begin{center}
\vspace{-0.5cm}
\begin{figure}[H]
\centering
\begin{tabular}{c}
\includegraphics[width=300pt]{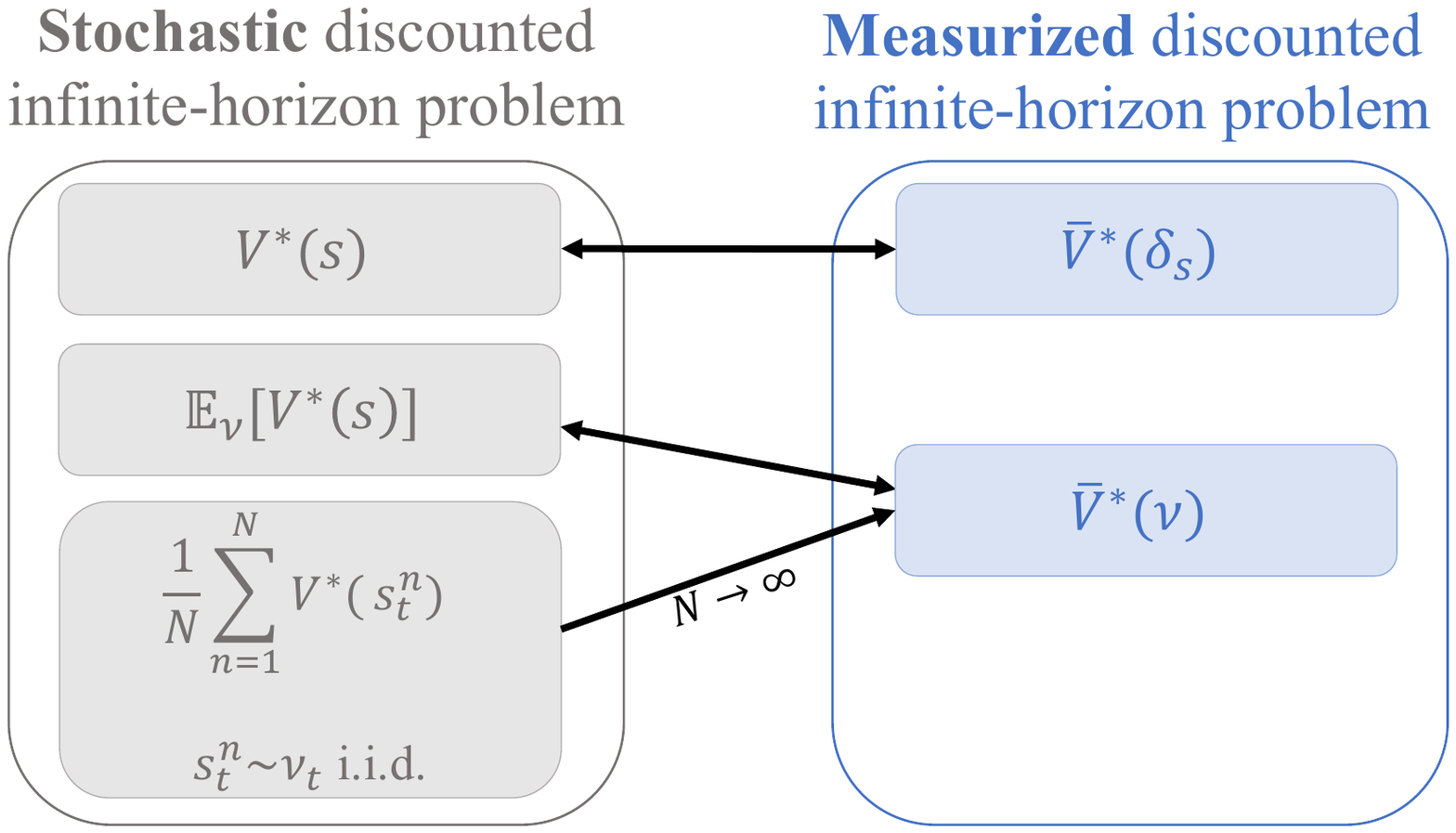}\\
\end{tabular}
\vspace{-0.5cm}
\caption{Connection between the original value function $V^*(\cdot)$ and the measurized value function $\oV^*(\cdot)$. Bidirectional arrows show equivalence; the unidirectional arrow indicates an asymptotic result. \label{fig:measurized_valuefunction}}
\end{figure}
\vspace{-0.5cm}
\end{center}

The following theorem shows that the stochastic and measurized value functions coincide if the latter is evaluated at a Dirac measure $\delta_s$ concentrated at the stochastic state $s$. In addition, it also shows that the supremum in \eqref{infinite_horizon_nu} can be interchanged with the integral, hence yielding that the measurized value function is the expected stochastic value function. \textcolor{black}{Analogous results for BS-deterministic MDPs are Corollary 9.3.2 and Proposition 9.5 in \cite{bertsekas}. Both results rely on the universally measurable version of the Ionescu-Tulcea Theorem \cite[Proposition 7.45]{bertsekas} along with a change of variables like \eqref{change_variables}, which, as previously noted, is presented in a more complicated form in \cite{bertsekas}. Additionally, Proposition 9.5 draws on results from the finite-horizon MDP and employs the exact selection theorem, leading to \(\epsilon\)-optimal policies. In contrast, we take advantage of the fact that the measurized MDP inherits assumptions from the original MDP (Proposition \ref{prop:assumptions}), allowing us to apply directly a result in the semicontinuous-semicompact framework (Theorem \ref{theorem423} in this paper).}


\begin{theorem} \label{th:V}
\textcolor{black}{Suppose all the assumptions in Proposition \ref{prop:assumptions} hold, and} let $V^*(\cdot)$  and $\V^*(\cdot)$ be the solutions to \eqref{DCOE} and \eqref{measurized_DCOE}, respectively. Then

\begin{flalign*}
(a) & \qquad \V^*(\delta_s)=V^*(s) & \forall  s \in \mS  
\\
(b) & \qquad  \oV^*(\nu)=\int_{\mS} V^*(s) d\nu(s)=\mathbb{E}_\nu [ V^*(s) ]  & \nu \in \mMpS 
\end{flalign*}

\end{theorem}

\proof{

We start by proving (a). Let $\nu_0=\delta_{s}$ and denote $\nu_t$ as the distribution of the states at time $t$, i.e. $s_t \sim \nu_t$. Proposition \ref{prop:distribution} showed that these distributions can be computed through the recursive equation $\nu_t=F(\nu_{t-1},\varphi_{t-1})$, for all $t=1,2,...$. According to Theorem \ref{theorem423}, we can restrict ourselves to deterministic or randomized Markov policies without loss of optimality. Then we can use Equation \eqref{infinite_horizon} to find an $\alpha$-discount optimal policy as

\begin{align*}
V^*(s)&=\sup_{\Pi^{MR}}  \ \mathbb{E}_s^\pi \left[ \sum_{t=0}^\infty \alpha^t r(s_t,u_t) \right] & s \in \mathcal{S}\\
&=\sup_{\Pi^{MR}}  \   \mathbb{E}_{P_s^\pi}  \left[ \sum_{t=0}^\infty \alpha^t r(s_t,u_t) \right] & s \in \mathcal{S}\\
&=\sup_{\Pi^{MR}}  \   \sum_{t=0}^\infty \alpha^t \ \mathbb{E}_{P_s^\pi}  \left[ r(s_t,u_t) \right] & s \in \mathcal{S}\\
&=\sup_{\Pi^{MR}}  \   \sum_{t=0}^\infty  \alpha^t  \int_{\mU} \int_{\mS} r(s_t,u_t)  \varphi_t(du|s) d\nu_t(s)&    s \in \mathcal{S},\ \nu_0=\delta_s,\\
&=\sup_{\Pi^{MR}}  \   \sum_{t=0}^\infty  \alpha^t  \ovr(\nu_t,  \varphi_t) &  s \in \mathcal{S}, \ \nu_0=\delta_s, \ \nu_t=F(\nu_{t-1},\varphi_{t-1}) \ t\geq 1\\
&= \V^*(\delta_s) &  s \in \mathcal{S},
\end{align*} 

\noindent where $ \mathbb{E}_{P_s^\pi}$ is the expected value taken with respect to the Ionescu-Tulcea probability measure $P_s^\pi$ introduced in Definition \ref{def:MDP}, the third equality is a consequence of the MCT and we use the definition of the revenue function $\ovr$ in the fifth equality. 

We now show (b). Putting together the definition of $\mathbb{V}^*(\cdot)$ in Equation \eqref{infinite_horizon_nu}, and its equivalence with $\overline{V}^*(\cdot)$ demonstrated in \eqref{M_infinite}, we have that 

\begin{equation}\label{sup_int}
\overline{V}^*(\nu)=\sup_{\pi} \int_\mS V(\pi,s) d\nu(s).
\end{equation}
Hence it suffices to prove that 

\begin{equation}\label{sup_int}
\sup_{\pi} \int_\mS V(\pi,s) d\nu(s)=\int_\mS V^*(s) d\nu(s).
\end{equation}
and that the supremum is attainable. The latter is true due to Theorem \ref{theorem423} applying to the measurized MDP, which inherits the assumptions of the original MDP as is shown in Proposition \ref{prop:assumptions}. To demonstrate \eqref{sup_int}, we first show that $\sup_{\pi} \int_\mS V(\pi,s) d\nu(s) \geq \int_\mS V^*(s) d\nu(s)$. This is obvious, since

\begin{align*}
\sup_{\pi} \int_\mS V(\pi,s) d\nu(s)& \geq \int_\mS V(\pi',s) d\nu(s) &\forall \pi' \in \Pi 
\end{align*}
In particular, this holds true for $\pi^*=\arg\max_\pi V(\pi,s)$, which can be assumed to be deterministic and stationary according to Theorem \ref{theorem423}. Hence we have that
\begin{align*}
\sup_{\pi} \int_\mS V(\pi,s) d\nu(s)& \geq \int_\mS V(\pi^*,s) d\nu(s) = \int_\mS V^*(s) d\nu(s) & 
\end{align*}

Finally, we show the opposite inequality $\sup_{\pi} \int_\mS V(\pi,s) d\nu(s) \leq \int_\mS V^*(s) d\nu(s)$ holds. Because $V(\pi,s) \leq V^*(s) $ for all $\pi \in \Pi$ and $s \in \mS$, we have that

\begin{align*}
\hspace{4cm} \int_\mS V(\pi,s) d\nu(s) &\leq \int_\mS V^*(s) d\nu(s) & \forall \pi \in \Pi. & \hspace{3.5cm} \blacksquare \\
\end{align*}
%

}


Under some mild conditions, the strong law of large numbers ensures that an empirical distribution converges to the true distribution almost surely as the sample size increases. Hence one naturally wonders if the measurized value function is also endowed with this property; i.e., if the average of the stochastic value functions evaluated at a sample of states converges to the measurized function of the sampling distribution. \textcolor{black}{This result has no parallel in \cite{bertsekas} and allows us to interpret the measurized value function as solving over infinite realizations of the MDP.} The following corollary indeed proves this.

\begin{corollary}\label{prop:asymptotics1}
\textcolor{black}{Suppose all the assumptions in Proposition \ref{prop:assumptions} hold.} Let $\nu_t \in \mMpS$ be the state distribution at time $t$, and let $s^1_t,...,s^N_t$ be $N$ states sampled from $\nu_t$. Then
\begin{equation}\label{eq:asymptotics1}
\oV^*(\nu_t)=\lim_{N\rightarrow \infty} \ \frac{1}{N} \sum_{n=1}^N V^*(s^n_t)
\end{equation}
\end{corollary}

\proof{
Denote the empirical initial state distribution as
\begin{align*}
\hat{\nu}_t^N(\mA)=\frac{1}{N} \sum_{n=1}^N \mathbf{1}_{\left\{s^n_t \in \mA\right\}} & &\forall \mA \in \mBS,
\end{align*}
where $\mathbf{1}_{\left\{s \in \mA\right\}}$ is the indicator function of set $\mA$. Then we have that

\begin{align*}
\frac{1}{N} \sum_{n=1}^N V^*(s^n_t) &=\int_\mS V^*(s) d\hat{\nu}^N_t(s)\\
&=\oV^*\left(\hat{\nu}_t^N\right),
\end{align*}
where the first equality comes by construction of $\hat{\nu}_t^N$ and the second comes from Theorem \ref{th:V}. Since $\hat{\nu}_t^N \to \nu_t$ as $N \rightarrow \infty$, we have that

\begin{align*}
\lim_{N \to \infty} \frac{1}{N} \sum_{n=1}^N V^*(s^n_t)&= \lim_{N \to \infty} \int_\mS V^*(s) d\hat{\nu}^N_t(s)\\
&=\int_\mS V^*(s) d(\lim_{N \to \infty} \hat{\nu}^N_t(s))\\
&= \oV^*(\nu_t),
\end{align*}
where the second equality follows from \textcolor{black}{Proposition C.12 in \cite{lerma}. This result applies because the sequence of measures $\{\hat{\nu}_t^N\}_{N \geq 1}$ converges setwise to $\nu_t$. Moreover, since $V^* \in \mathcal{V}(\mS)$, there exists a constant $M$ such that $V^*(s) \leq M$ for all $s \in \mS$. Thus, the integrand $V^*(\cdot)$ can be viewed as a degenerate sequence bounded above by $M$, with $\int_{\mS} M \, d\hat{\nu}_t^N = \int_{\mS} M \, d\nu_t$. Therefore, the conditions of Proposition C.12 in \cite{lerma} are satisfied.}
 \hfill $\blacksquare$\\

}

\subsection{Relationship between optimal policies. \label{sec:policies}}
In this section, we show the equivalence between the measurized and original actions, summarized in Figure \ref{fig:measurized_actions}.

\begin{center}
\vspace{-1cm}
\begin{figure}[H]
\centering
\begin{tabular}{c}
\includegraphics[width=250pt]{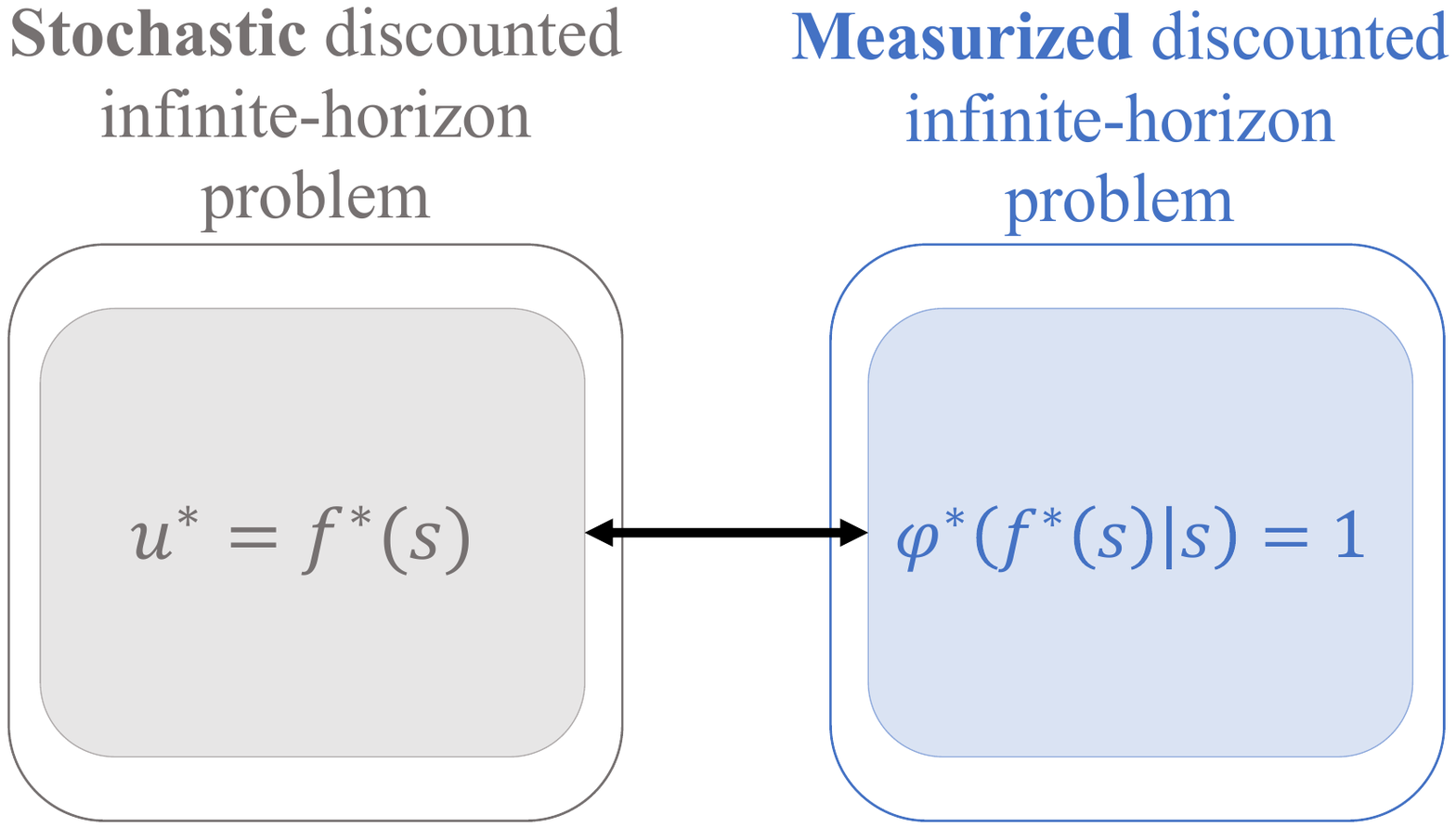}\\
\end{tabular}
\caption{Equivalence between the optimal actions $\varphi^*$ of measurized MDP and the optimal actions $u^*$ of the original MDP it was lifted from. \label{fig:measurized_actions}}
\end{figure}
\vspace{-0.5cm}
\end{center}

The following result shows that the optimal policy of the measurized optimality equations \eqref{measurized_DCOE} coincides with the optimal policy to \eqref{DCOE}; hence it belongs to $\Pi^D$ and is attainable. \textcolor{black}{Although not explicitly stated in \cite{bertsekas}, this result can be inferred from Corollary 9.5.1 and Definition 9.9 therein.}


\begin{theorem}\label{prop:unique_phi}
\textcolor{black}{Suppose all the assumptions in Proposition \ref{prop:assumptions} hold.}  Then the optimal decision rule $\varphi^*$ obtained by solving \eqref{measurized_DCOE} does not depend on the current state distribution $\nu$ and is concentrated around the optimal selector $f^* $ attaining the maximum in \eqref{DCOE};  that is to say $\varphi^*(f^*(s)|s)=1$ for all  $s \in \mS$.
\end{theorem}

\proof{
The optimality equations corresponding to the infinite horizon {\it measurized} problem \eqref{M_infinite} are \eqref{measurized_DCOE}, where

\begin{align*}
\overline{V}^*(\nu)&=\sup_{\varphi \in \Phi(\nu)} \ \ovr(\nu,\varphi) + \alpha   \overline{V}^*(F(\nu,\varphi)) & \forall \nu \in \mMpS\\
&=\sup_{\varphi \in \Phi(\nu)} \  \int_{\mS} \int_\mU \left\{ r(s,u)+\alpha \int_\mS V^*(s') Q(ds'|s,u) \right\} \varphi(du|s) d\nu(s) & \forall \nu \in \mMpS\\
&=\sup_{\varphi \in \Phi(\nu)} \  \int_{\mS} \int_\mU \left\{ r(s,u)+\alpha \mathbb{E}_Q [V^*(s')|s,u] \right\} \varphi(du|s) d\nu(s) & \forall \nu \in \mMpS
\end{align*}
where the second equality comes from Theorem \ref{th:V} (b) and the definition of the deterministic transition F. Theorem \ref{theorem423} claims that, if Assumptions \ref{assumption421} and \ref{assumption422} hold, then there exists a selector $f^* \in \mathbb{F}$ such that
\begin{align*}
V^*(s)&=\max_{u\in \mU(s)} \left\{r(s,u)+\alpha \mathbb{E}_Q[V^*(s')|s,u] \right\}=r(s,f^*(s))+\alpha \mathbb{E}_Q[V^*(s')|s,f^*(s)] & \forall s \in \mS,
\end{align*}
which means that for all $u \in \mU(s)$
\begin{align*}
r(s,u)+\alpha \mathbb{E}_Q[V^*(s')|s,u] &\leq r(s,f^*(s))+\alpha \mathbb{E}_Q[V^*(s')|s,f^*(s)]=V^*(s) & \forall s \in \mS.
\end{align*}
Hence, for any fixed $\nu\in \mMpS$ and $\varphi \in \Phi$ we have that
\begin{align*}
\int_{\mS} \int_\mU \left\{ r(s,u)+\alpha \mathbb{E}_Q[V^*(s')|s,u]\right\} \varphi(du|s) d\nu(s) &\medmath{\leq \int_{\mS} \int_\mU \left\{ r(s,f^*(s))+\alpha \mathbb{E}_Q[V^*(s')|s,f^*(s)] \right\} \varphi(du|s) d\nu(s)} & \\
&= \int_{\mS}  \left\{ r(s,f^*(s))+\alpha \mathbb{E}_Q[V^*(s')|s,f^*(s)] \right\}  d\nu(s) & \\
&= \int_{\mS} \int_\mU \left\{ r(s,u)+\alpha \mathbb{E}_Q[V^*(s')|s,u)] \right\} \varphi^*(du|s) d\nu(s), & \\
\end{align*}
where $\varphi^* \in \Phi$ is the stochastic kernel concentrated around the selector as $\varphi^*(f^*(s)|s)=1$ for all $s \in \mS$. Therefore,

\begin{align*}
\sup_{\varphi \in \Phi(\nu)} \int_{\mS} & \int_\mU \left\{ r(s,u)+\alpha \mathbb{E}_Q[V^*(s')|s,u]\right\} \varphi(du|s) d\nu(s) &\forall \nu \in\mMpS\hspace{2cm}\\
 &\leq  \sup_{\varphi \in \Phi} \int_{\mS} \int_\mU \left\{ r(s,f^*(s))+\alpha \mathbb{E}_Q[V^*(s')|s,f^*(s)] \right\} \varphi^*(du|s) d\nu(s)  &\forall \nu \in\mMpS\hspace{2cm}\\
&  =\ovr(\nu,\varphi^*) + \alpha   \overline{V}^*(F(\nu,\varphi^*))=\oV^*(\nu) &\forall \nu \in\mMpS \hspace{2cm} \blacksquare
\end{align*}
}

As a consequence, we can replace the supremum by a maximum in \eqref{measurized_DCOE}, over $\varphi \in\Phi$ independently of $\nu$.


\subsection{Connection to the Linear Programming formulation. \label{sec:dualizing}}
In this section we demonstrate that the dual variables $\zeta$ of the LP formulation \eqref{LP} can be interpreted as a discounted sum of state distributions $\{\nu_t\}_{t\geq 0}$ and measurized controls $\{\varphi_t\}_{t\geq 0}$. More specifically, we implement the change of variables \eqref{change_variables} introduced in Section \ref{sec:change_variables} to connect $\zeta$ with the states and actions of the measurized MDP. Figure \ref{fig:change_variables} visually demonstrates the connection.

 \begin{center}
 \vspace{-0.3cm}
\begin{figure}[H]
\centering
\begin{tabular}{c}
\includegraphics[width=270pt]{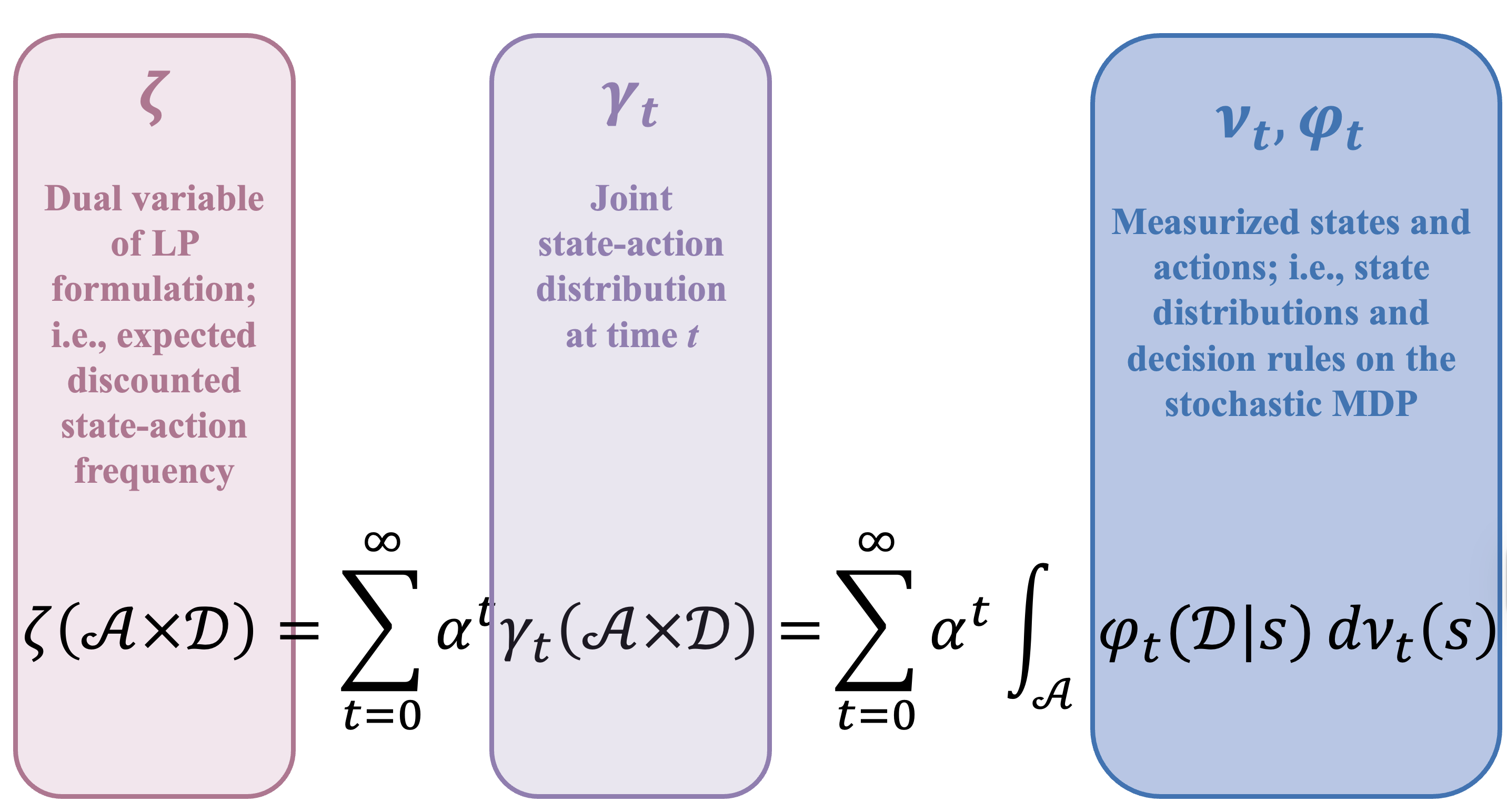}\\
\end{tabular}
\vspace{-0.5cm}
 \caption{Visual representation of the relationship between the dual variables $\zeta$ of the LP formulation of the stochastic MDP, the joint state-action distribution at time $t$, $\gamma_t$, and the measurized states and actions. \label{fig:change_variables}}
\end{figure}
\vspace{-0.5cm}
\end{center}

%
This allows us to think of measurizing as dualizing the standard LP formulation of \eqref{DCOE} and performing some changes of variables, as illustrated in Figure \ref{fig:LPmeasurizing}. \textcolor{black}{Although $\gamma$ coincides with the actions in BS-deterministic MDPs, which incorporate information about the state distribution, the connection between the lifting process and the dual of the LP formulation was not explicitly addressed in \cite{bertsekas}.}

\begin{center}
\vspace{-0.1cm}
\begin{figure}[h!]
\centering
\begin{tabular}{c}
\includegraphics[width=350pt]{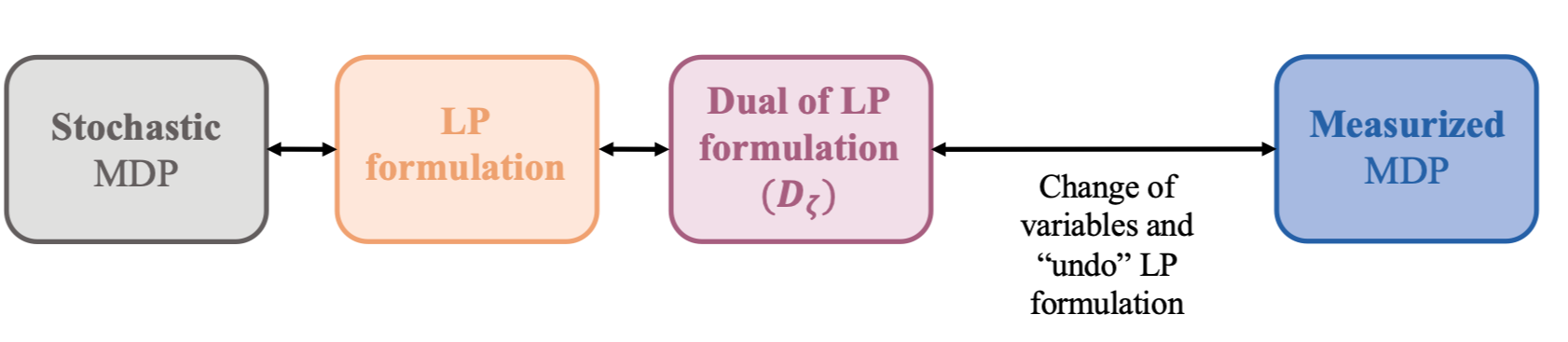}\\
\end{tabular}
 \caption{Illustration of the relationship between the LP formulation of the stochastic MDP and the measurized MDP. This offers an alternative recipe for {\it measurizing} any MDP. \label{fig:LPmeasurizing}}
\end{figure}
\vspace{-0.5cm}
\end{center}

%



We now establish the connection between the measurized value function and the dual problem \eqref{D_zeta}.
Using Theorem \ref{theorem423}, Theorem \ref{th:V} and Proposition \ref{prop:change_variables_t}, one can perform the change of variables \eqref{change_variables} in the discounted infinite-horizon problem \eqref{oV*}, yielding

\begin{equation}\label{optimal_gamma}
\oV^*(\nu_0)=\sup_{\Gamma \in \Omega(\nu_0)}   \ \sum_{t=0}^\infty \alpha^t  \mathbb{E}_{\gamma_t} \left[  r(s,u) \right]\qquad  \nu_0 \in \mMpS.
\end{equation}
Because the reward function $r$ does not depend on time $t$, intuitively one could perform in \eqref{optimal_gamma} the following change of variables 
\begin{equation}\label{change_zeta} \tag{$\zeta$}
\zeta(\mA\times\mathcal{D})=\sum_{t=0}^\infty \alpha^t \gamma_t(\mA\times\mathcal{D}) \qquad \forall \mA \in\mB(\mS), \mathcal{D}\in \mB(\mU).
\end{equation}
Recall that $\zeta(\mS\times\mU)=1/(1-\alpha)$, so $\zeta$ is not a probability measure. Performing the change of variables \eqref{change_zeta} would lead to an optimization problem 

\begin{equation}\label{optimal_zeta}
\oV^*(\nu_0)=\sup_{\zeta \in Z(\nu_0)}   \  \int_\mS \int_\mU r(s,u) d\zeta(s,u)  \qquad  \nu_0 \in \mMpS, 
\end{equation}
where the set $Z(\nu_0)$ gathers the transition law of $\gamma_t$ and thus $\nu_t$, having
\begin{align*}
Z(\nu_0):=&\left\{\zeta \in \mM_+(\mS\times\mU):  \  \zeta (\mA\times\mU(\mA)^c)=0 \ \ \forall \mA \in \mB(\mS), \phantom{\int_{\mS}}  \right.\\
& \left. \quad \zeta(\mA\times\mU)=\nu_0(\mA) +\alpha  \int_{\mS} \int_{\mU} Q(\mA|s,u) d\zeta(s,u)  \ \ \forall \mA \in \mB(\mS)\right\}. 
\end{align*}
With this specification, Problem \eqref{optimal_zeta} is actually the dual \eqref{D_zeta} of the LP formulation of the stochastic MDP. This means that, if we can perform the change of variables \eqref{change_zeta} without loss of optimality, one can view the measurized MDP as equivalent to the {\it dual} of the stochastic MDP. The following proposition shows this.

\begin{proposition}\label{prop:zeta}
\textcolor{black}{Suppose all the assumptions in Proposition \ref{prop:assumptions} hold.} Let $\zeta^*$ be the optimal solution to \eqref{D_zeta}. Let $\oV^*(\cdot)$ and $\varphi^*$ be the optimal measurized value function and decision rule solving \eqref{measurized_DCOE}. Denote the optimal state-distribution path that $\varphi^*$ gives rise to as $\Upsilon^*=\{\nu_t^*\}_{t \geq 0}$, where $\nu_0^*$ coincides with the initial state distribution $\nu_0$. Then:

\begin{itemize}
\item[(a)] $\oV^*(\nu_0)$ coincides with the optimal objetive of the dual problem \eqref{D_zeta}.
\item[(b)] Without loss of optimality we can assume that $\zeta^*(\mA\times\mD)=\sum_{t=0}^\infty \alpha^t \int_\mA \varphi^*(\mD|s) d\nu_t^*(s)$ for all $\mA \in \mBS$, $\mD \in \mBU$.
\end{itemize}
\end{proposition}

\proof{

The proof of (a) is straightforward, since we have that 

$$\oV^*(\nu_0)=\int_S V^*(s) d\nu_0(s)=\int_\mS \int_\mU r(s,u) d\zeta^*(s,u),$$
where the first equality comes from Theorem \ref{th:V}.(b) and the second equality is given by Theorem \ref{th:HLdual}.

To show (b), we construct the following measure for any initial state distribution $\nu_0 \in \mMpS$ 

$$\tilde{\zeta}(\mA\times\mD):=\sum_{t=0}^\infty  \alpha^t \int_\mA \varphi^*(\mD|s) d\nu_t^*(s), \qquad \forall \ \mA \in \mBS, \ \mD \in \mBU.$$
To prove (b), it suffices to show that

\begin{itemize}
\item[(i)] $\tilde{\zeta}$ is feasible to \eqref{D_zeta}.
\item[(ii)] $\tilde{\zeta}$ gives the same objective value as $\zeta^*$ in \eqref{D_zeta}; i.e. $\int_\mS \int_\mU r(s,u) d\tzeta(s,u)=\int_\mS \int_\mU r(s,u) d\zeta^*(s,u)$.
\end{itemize}

We start by proving (i):

\begin{align*}
\tzeta(\mA,\mU)& = \sum_{t=0}^\infty \alpha^t \int_\mA \varphi^*(\mU|s) d\nu_t^*(s)=\sum_{t=0}^\infty \alpha^t \nu_t^*(\mA)\\
&=\nu_0^*(\mA) + \sum_{t=1}^\infty \alpha^t  F(\varphi^*,\nu^*_{t-1})(\mA)\\
&=\nu_0^*(\mA) + \sum_{t=1}^\infty \alpha^t  \int_\mS \int_\mU Q(\mA|s,u) \varphi^*(du|s) d\nu^*_{t-1}(s)\\
&=\nu_0^*(\mA) +  \alpha \int_\mS \int_\mU Q(\mA|s,u) \sum_{t=0}^\infty \alpha^t  \varphi^*(du|s) d\nu^*_{t}(s)\\
&=\nu_0^*(\mA) +  \alpha \int_\mS \int_\mU Q(\mA|s,u) d\tzeta(s,u),
\end{align*}
where the fifth equality comes from the MCT. In addition,

$$\tzeta(\mA,\mU(\mA)^c)=\sum_{t=0}^\infty \alpha^t \int_\mA \varphi^*(\mU(\mA)^c|s) d\nu_t^*(s)=0$$
because $\varphi^*(\mU(s)|s)=1$ for all $s \in \mS$. Therefore, $\tzeta \in \mM_+(\mathbb{K})$ and it fulfills the constraint in \eqref{D_zeta}.

We now prove (ii):
\begin{align*}
\int_\mS \int_\mU r(s,u) d\zeta^*(s,u)&=\oV^*(\nu_0)=\sup_{\substack{\varphi_t \in \Phi, \\ t\geq 0}}  \mathbb{E}_{\nu_0}^\pi \left[ \sum_{t=0}^\infty \alpha^t r(s_t,u_t)\right] \\
 &= \int_\mS \int_\mU  \sum_{t=0}^\infty \alpha^t r(s,u) \varphi^*(du|s) d\nu_t^*(s) \\
  &= \int_\mS \int_\mU   r(s,u) d\tzeta(s,u), 
\end{align*}
where the first equality comes from (a) and the rest from plugging the optimal decision rule $\varphi^*$ and the optimal path of state distributions $\Upsilon^*$ in the measurized discounted infinite horizon expression \eqref{M_infinite}. \hfill $\blacksquare$\\
}






\section{An algebraic procedure for lifting the stochastic optimality equations. \label{sec:measurizing}}
In this section, we introduce a novel method for lifting the {\it stochastic} MDP to the measure-valued framework, offering a more intuitive approach compared to the one described in Section \ref{sec:dualizing}.
Figure \ref{fig:measurizing} outlines the process.

\begin{center}
\vspace{-0.1cm}
\begin{figure}[H]
\centering
\begin{tabular}{c}
\includegraphics[width=350pt]{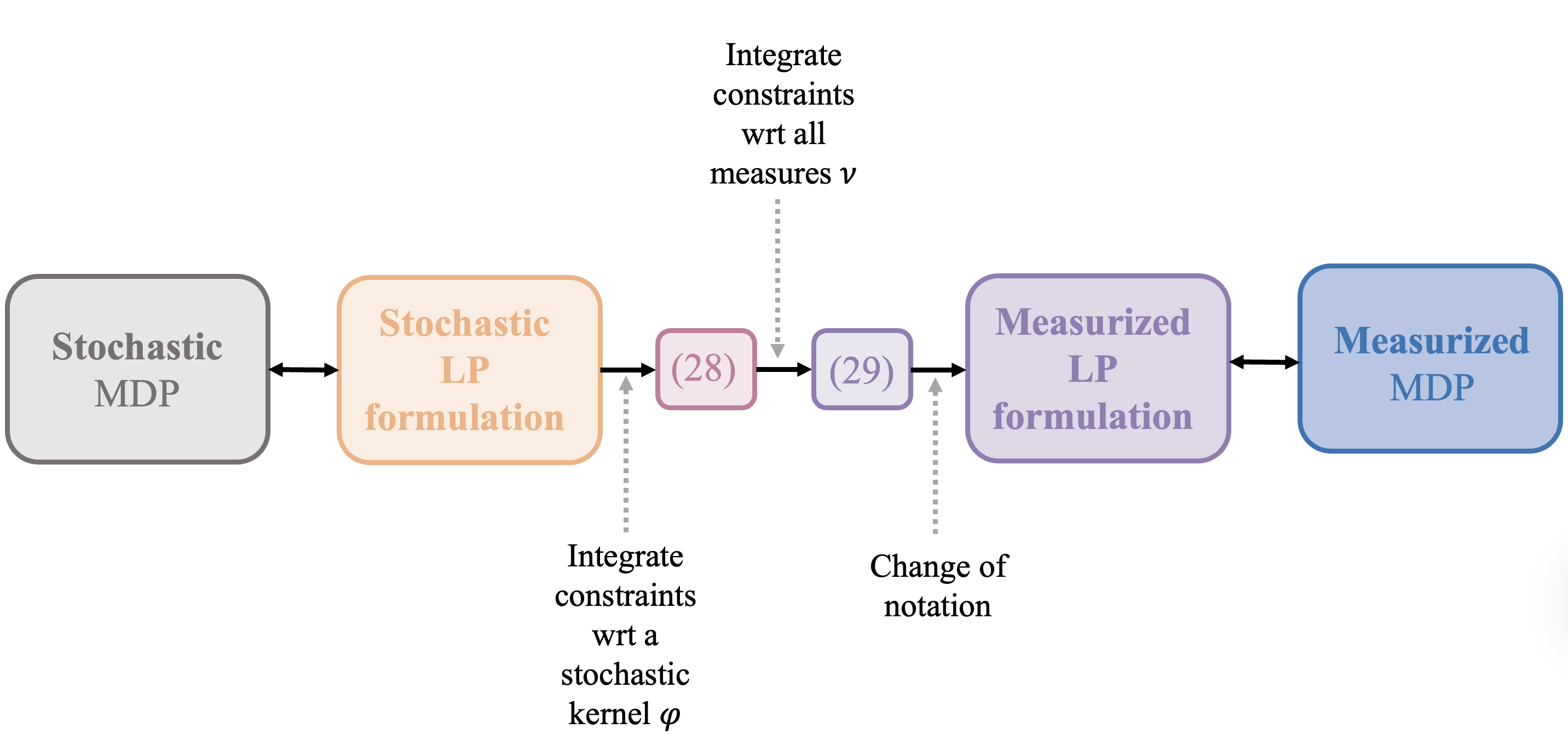}\\
\end{tabular}
 \caption{Illustration of the {\it measurizing} process; i.e. how one can intuitively lift any stochastic MDP to the measure-valued framework. \label{fig:measurizing}}
\end{figure}
\vspace{-0.5cm}
\end{center}

We start with formulation \eqref{LP}. The first step is to integrate the constraints over the action space $\mU$ using a stochastic kernel $\varphi \in \Phi$, where $\Phi$ is defined as in \eqref{eq:Phi}:

\begin{align}
\inf_{\substack{ V \in \mathcal{V}}} & \ \int_{\mathcal{S}} V(s) d\nu_0(s) \label{measurized_DCOE0} \\
\mbox{s.t.} & \ \ V(s)  \geq \int_{\mathcal{U}}  r(s,u) \varphi(du|s)+ \alpha \int_{\mathcal{U}}  \mathbb{E}_Q \left[ V(s')| s,u\right]  \varphi(du|s) \qquad \qquad \forall  s \in \mathcal{S}, \varphi \in \Phi. \nonumber
\end{align}

\noindent The second step involves aggregating the constraints in Problem \eqref{measurized_DCOE0} over all probability measure $\nu \in \mathcal{M}(\mathcal{S})$:

\begin{align}
\inf_{\substack{V \in \mathcal{V}}} & \ \int_{\mathcal{S}} V(s) d\nu_0(s) \label{measurized_DCOE1} \\
\mbox{s.t.} & \int_{\mathcal{S}} V(s) d\nu(s) \geq \int_{\mathcal{S}} \int_{\mathcal{U}}  \left\{ r(s,u) + \alpha \mathbb{E}_Q \left[ V(s')| s,u\right] \right\} \varphi(du|s) d\nu(s) \qquad \qquad \forall \nu \in \mMpS, \varphi \in \Phi. \nonumber
\end{align}

\noindent This constrained LP can be interpreted as encompassing all convex combinations of the constraints defined for each feasible state-action pair \((s, u)\). As the final step, we express this problem using the measurized notation. Since we optimize \(V(\cdot)\) within the space of bounded measurable functions on \(\mS\), denoted as \(\mathcal{V}(\mS)\), the objective in \eqref{measurized_DCOE1} remains well-defined and bounded. Denote

\begin{equation}\label{eq:measurizedV}
\overline{V}(\nu):= \int_{\mathcal{S}} V(s) d\nu(s)=\mathbb{E}_\nu [V(s)] < \infty;
\end{equation}
then $\overline{V} \in \mathcal{V}(\mMpS)$, where $\mathcal{V}(\mMpS)$ denotes the space of bounded measurable functions in the space of distributions $\mMpS$. Then Problem \eqref{measurized_DCOE1} can be rewritten as

\begin{align}
\inf_{\substack{\overline{V} \in \mathcal{V}(\mMpS) }} & \ \overline{V}(\nu_0) \label{measurized_DCOE2} \\
\mbox{s.t.} &\ \overline{V}(\nu) \geq \ovr(\nu,\varphi) + \alpha \int_{\mathcal{S}} \int_{\mathcal{U}}  \mathbb{E}_Q \left[ V(s')| s,u\right]  \varphi(du|s) d\nu(s) \qquad \qquad \forall \nu \in \mMpS, \varphi \in \Phi \nonumber
\end{align}

\noindent Now decompose the expectation in \eqref{measurized_DCOE2} as

\begin{align*} 
\int_{\mathcal{S}} \int_{\mathcal{U}}  \mathbb{E}_Q \left[ V(s')| s,u\right]  \varphi(du|s) d\nu(s) &=\int_{\mathcal{S}} \int_{\mathcal{U}}  \int_{\mathcal{S}} V(s') Q(ds'|s,u) \varphi(du|s) d\nu(s)\\ 
&=\int_{\mathcal{S}}  V(s') \int_{\mathcal{U}}  \int_{\mathcal{S}} Q(ds'|s,u) \varphi(du|s) d\nu(s) \nonumber\\ 
&= \int_{\mathcal{S}} V(s') dF(\nu,\varphi)(s') \nonumber \\
&= \overline{V}(F(\nu,\varphi))= \overline{V}(\nu'), \nonumber
\end{align*}
where $F$ is as defined in \eqref{F} and Fubini's Theorem ensures that the order of integration can be switched\footnote{Fubini's Theorem can be applied because $V\in \mathcal{V}$ is a measurable function and $Q(\cdot|s,u)$, $\varphi(\cdot|s)$ and $\nu(\cdot)$ are probability measures for all $s \in \mS$ and $u \in \mU(s)$, thus they are $\sigma$-finite measures.}. 
This new notation is consistent with the notion of the optimal measurized value function being the expected value of the optimal original value function, as demonstrated in Theorem \ref{th:V}.
Hence we can rewrite Problem \eqref{measurized_DCOE1} as

\begin{align}
\inf_{\substack{\overline{V} \in \mathcal{V}(\mMpS) }} & \ \overline{V}(\nu_0) \label{measurized_DCOE3} \\
\mbox{s.t.} &\ \overline{V}(\nu) \geq \ovr(\nu,\varphi) + \alpha \oV(F(\nu,\varphi)) \qquad \qquad \forall \nu \in \mMpS, \varphi \in \Phi ,\nonumber
\end{align}

\noindent which coincides with the measurized LP formulation \eqref{measurized_LP} when $\Phi(\nu)=\Phi$ for all $\nu$. As discussed in Section \ref{sec:OE}, the optimal measurized value function can be obtained by solving \eqref{measurized_DCOE3}. Moreover, the theory in Section \ref{sec:measurized} ensures there is no optimality gap between \eqref{LP} and \eqref{measurized_DCOE3}.

This lifting procedure provides consistency to derive the measurized counterpart of an MDP and enhances interpretability. For instance, it is clear that a tightened measurized MDP entails a relaxation of the measurized LP; i.e., we do not consider the constraints for all pairs $(\nu,\varphi) \in \mMpS \times \Phi$ but a subset of these as $\Phi(\nu) \subseteq \Phi$ for all $\nu$. However, if all the constraints $(s,u) \in \mathbb{K}$ are still present in \eqref{measurized_LP}, one still obtains $\oV^*(\delta_s)=V^*(s)$. To see this, it suffices to solve \eqref{measurized_LP} with $\nu_0$ concentrated at $s$.

\subsection{Connection to measure-valued MDPs with stochastic transitions. \label{sec:measure-valued}}

In this section, we explore how the lifting procedure helps us understand general measure-valued MDPs as stochastic MDPs with external sources of randomness. We illustrate this through Figure \ref{fig:measure-valued-MDP} and the following example.

\begin{center}
\vspace{-0.1cm}
\begin{figure}[H]
\centering
\begin{tabular}{c}
\includegraphics[width=350pt]{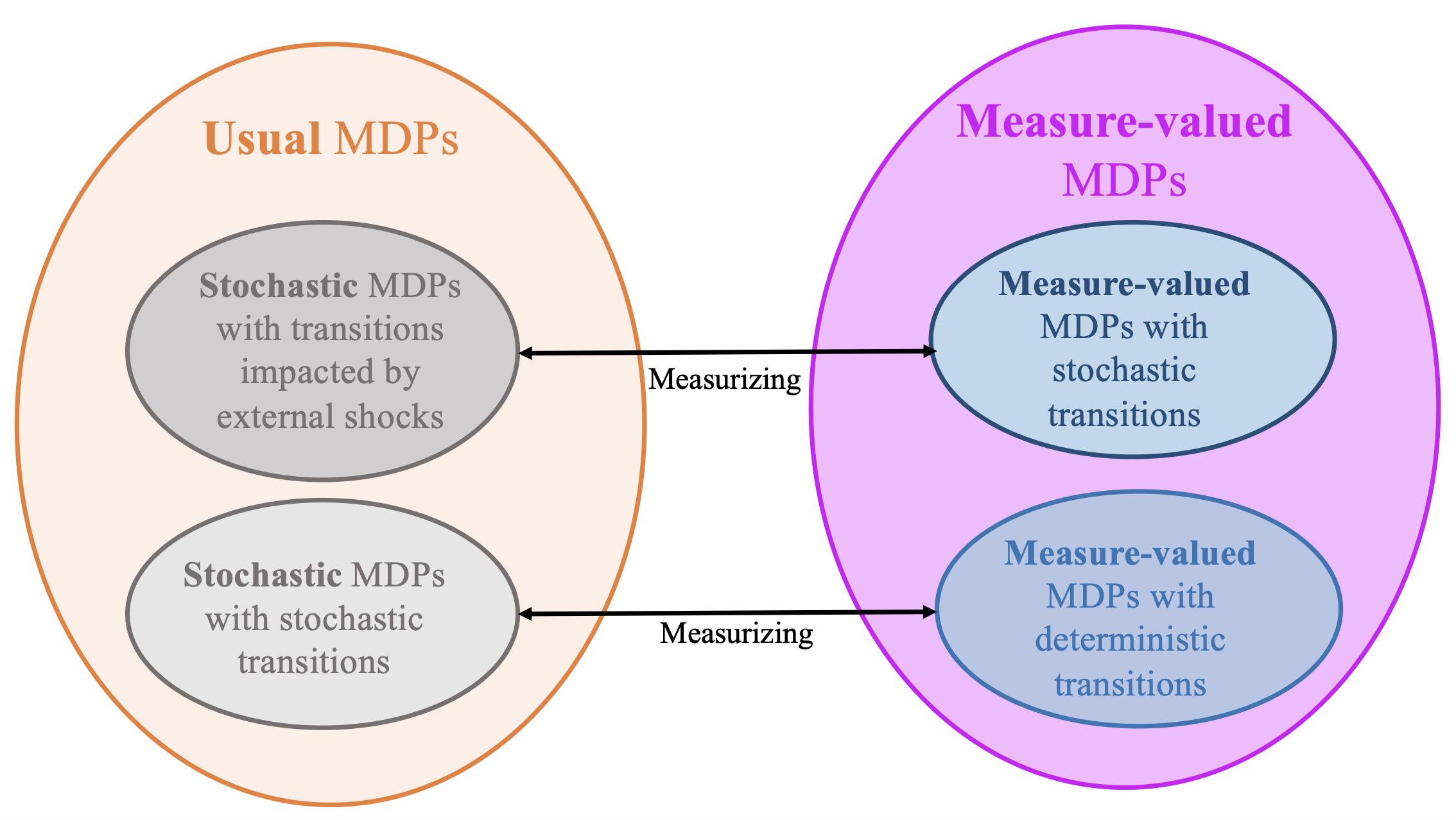}\\
\end{tabular}
 \caption{Illustration of how the {\it measurizing} process lifts different types of MDP into the measure-valued framework. \label{fig:measure-valued-MDP}}
\end{figure}
\vspace{-0.5cm}
\end{center}

\begin{example}[Measure-valued MDPs as lifted MDPs with random transitions] \label{ex:MVMDP}
Consider the usual stochastic MDP with transition kernel $Q$ contingent on random shocks $z \in \mathcal{Z}$, following a distribution $\mu$. More specifically, $Q(\mA|\cdot)$ is a measurable function defined on $\mathbb{K}\times \mathcal{Z} $ for all $\mA \in \mBS$, and $Q(\cdot|s,u,z) $ is a probability measure on $\mS$ for every fixed $((s,u),z) \in \mathbb{K}\times \mathcal{Z}$. The optimality equations are 

\begin{equation*}
V^*(s) =\max_{u \in \mathcal{U}(s)} \  \left\{ r(s,u) + \alpha \int_{\mathcal{Z}} \int_\mS \ V^*(s') Q(ds'|s,u,z)  d\mu(z) \right\}, \qquad \forall s \in \mS. 
\end{equation*}

\noindent The LP formulation is

\begin{align*}
\inf_{\substack{ V \in \mathcal{V}}} & \ \int_{\mathcal{S}} V(s) d\nu_0(s) 
\\
\mbox{s.t.} & \ \ V(s)  \geq r(s,u) + \alpha \int_{\mathcal{Z}} \int_\mS \ V^*(s') Q(ds'|s,u,z)  d\mu(z) \qquad \qquad \forall  s \in \mathcal{S}, u \in \mU(s). \nonumber
\end{align*}

\noindent Implementing the lifting procedure yields

\begin{align*}
\inf_{\substack{ V \in \mathcal{V}}} & \ \int_{\mathcal{S}} V(s) d\nu_0(s) 
\\
\mbox{s.t.} & \ \ \medmath{ \int_\mS V(s) d\nu(s) \geq \int_\mS \int_\mU r(s,u)\varphi(du|s) d\nu(s) + \alpha  \int_\mS \int_\mU \int_{\mathcal{Z}} \int_\mS \ V^*(s') Q(ds'|s,u,z)  d\mu(z) \varphi(du|s) d\nu(s)  \quad \forall \nu \in \mMpS, \varphi \in \Phi}. \nonumber
\end{align*}
Translating this into measurized notation, we get

\begin{align}
\inf_{\substack{\overline{V} \in \mathcal{V}(\mMpS) }} & \ \overline{V}(\nu_0) \label{measurized_MDPz} \\
\mbox{s.t.} &\ \overline{V}(\nu) \geq \ovr(\nu,\varphi) + \alpha \int_{\mathcal{Z}}  \oV(\nu_z)  d\mu(z), \qquad \qquad \forall \nu \in \mMpS, \varphi \in \Phi \nonumber
\end{align}

\noindent where $\nu_{z,\nu,\varphi}(\cdot):=\int_\mS \int_\mU    Q(\cdot|s,u,z)  \varphi(du|s) d\nu(s)$ is the {\it deterministic} distribution of states in the subsequent period of time given the shock $z \in \mathcal{Z}$, a state distribution $\nu$ and a decision rule $\varphi$. The {\it random} distribution $\nu'$ of next states has a probability law given by the measure-valued stochastic kernel $\oq$. However, as seen in \eqref{measurized_MDPz} $\nu'$ is also tied to the distribution $\mu$ of the random shocks. More specifically, given a set $\mathcal{P}$ of the Borel $\sigma$-algebra defined on $\mMpS$, we have that

$$\mathbb{P}_{\oq}(\nu' \in \mathcal{P})=\mathbb{P}_\mu(\nu_{z,\nu,\varphi} \in \mathcal{P}) \qquad \forall \mathcal{P} \in \mathcal{B}(\mMpS), \nu \in \mMpS, \varphi \in \Phi. $$
In particular, we can write

$$\oq(\mathcal{P}|\nu,\varphi):= \int_{\mathcal{Z}} \mathbf{1}_{\{\nu_{z,\nu,\varphi} \in \mathcal{P}\}} d\mu(z) \qquad \forall \mathcal{P} \in \mathcal{B}(\mMpS), \nu \in \mMpS, \varphi \in \Phi,$$
thus having 

\begin{align}
\inf_{\substack{\overline{V} \in \mathcal{V}(\mMpS) }} & \ \overline{V}(\nu_0) \label{measurized_MDPq} \\
\mbox{s.t.} &\ \overline{V}(\nu) \geq \ovr(\nu,\varphi) + \alpha \mathbb{E}_{\oq}[\oV^*(\nu')| \nu,\varphi ], \qquad \qquad \forall \nu \in \mMpS, \varphi \in \Phi. \nonumber
\end{align}
According to \cite{lerma}, under certain assumptions the optimal solution to \eqref{measurized_MDPq} coincides with the optimal measure-valued function solving the optimality equations

\begin{align}
\overline{V}^*(\nu)&=\sup_{\varphi \in \Phi} \ \left\{ \ovr(\nu,\varphi) + \alpha  \mathbb{E}_{\oq}[ \overline{V}^*(\nu') | \nu,\varphi ] \right\} \qquad \forall \nu \in \mMpS. \label{measure-valued_DCOE}  \
\end{align}
 \phantom{aa} \hfill   \raisebox{1cm}{$\blacksquare$}

 \end{example}

 \vspace{-1.75cm}


\section{The Average-Reward problem \label{sec:Measurized_AR}} 
Under the AR criterion one seeks a policy that maximizes the average reward per unit time as the time horizon approaches infinity. Formally, the optimal AR value function can be retrieved through solving

\begin{equation}\label{AR}
J^*(s_0)=\sup_{\pi \in \Pi}  \left\{ \liminf_{n\rightarrow \infty} \mathbb{E}_{s_0}^{\pi} \left[ \frac{1}{n}\sum_{t=0}^{n-1} r(s_t,u_t)\right] \right\}  \qquad \forall s_0 \in \mS.
\end{equation}
Similarly to the discounted case, under certain assumptions the optimal action to take from any state $s \in \mS$ is deterministic and can be derived by the AR Optimality Equations (AROE):

\begin{align}\label{AROE}\tag{AROE}
\rho^* +h^*(s) & = \sup_{\mU(s)} \left\{ r(s,u) + \mathbb{E}_Q [ h^*(s')| s, u] \right\} & \forall s \in \mS, \\
\rho^* +h^*(s) &= r(s,f^*(s)) + \mathbb{E}_Q [ h^*(s')| s, f^*(s)] & \forall s \in \mS.  \nonumber
\end{align}

\noindent where $\rho^*$ is a constant called the {\it gain}, $h^*(\cdot)$ a real-valued measurable function called the optimal {\it bias} function, and $f^*$ is a {\it selector}. 
The optimal AR value function $J^*(\cdot)$ doesn't appear directly in \eqref{AROE}, but it can be shown to be constant and equal to $\rho^*$; i.e., $J^*(s) = \rho^*$ for all $s \in \mS$. For more details on solving the standard MDP under the AR criterion in Borel state and action spaces, see Appendix \ref{appendix:MDP} and also \cite{lerma}.

We now analyze the AR criterion within the context of measurized MDPs;

this criterion was not explored in \cite{bertsekas,BSdeterministic}.

Using the definition of measurized MDP, we now define the measurized AR-optimal value function.

\begin{equation}\label{M-AR-optimal}
\oJ^*(\nu_0)=\sup_{\overline{\pi} \in \overline{\Pi}} \left\{ \liminf_{n\rightarrow \infty} \frac{1}{n} \sum_{t=0}^{n-1} \ovr(\nu_t,\varphi_t)\right\} \qquad \forall \nu_0 \in \mMpS.
\end{equation}
where every $\overline{\pi} \in \overline{\Pi}$ yields a corresponding sequence $\{\nu_t,\varphi_t\}_{t\geq 0}$. We say that  a policy $\overline{\pi}^*$ is $\mM$-AR-optimal if it is an optimal solution to \eqref{M-AR-optimal}. The measurized AR MDP gives rise to the optimality equations 

\begin{align}\label{measurized_AROE}\tag{$\mM$-AROE}
\overline{\rho}^* +\ovh^*(\nu) & = \sup_{\varphi \in \Phi(\nu)} \left\{ \ovr(\nu,\varphi) + \ovh^*(F(\nu,\varphi)) \right\} \qquad \forall \nu \in \mMpS \\
\orho^* +\ovh^*(\nu) &= \ovr(\nu,\overline{f}^*(\nu)) + \ovh^*(F(\nu,\of^*(\nu)))  \qquad \ \forall \nu \in \mMpS,  \nonumber
\end{align}
where $\orho^*$ is the measurized optimal gain, $\ovh^*(\cdot)$ is a real-valued measurable function on $\mMpS$ and $\of^* \in \overline{\mathbb{F}}$ a measurized selector; i.e., $\of^*(\nu)=\varphi \in \Phi(\nu)$ for all $\nu \in \mMpS$. Finally, the measurized LP formulation is given by

\begin{align}
\inf_{\orho,\ovh(\cdot)} & \ \ \orho \label{measurized_LPrho} \tag{$\mM$-LP$_\rho$} \\
\mbox{s.t.} & \ \ \orho \geq  \ovr(\nu,\varphi)+  \ovh(F(\nu,\varphi)) -\ovh(\nu)   \qquad \qquad \forall  \nu \in \mMpS, \varphi \in \Phi. \nonumber
\end{align}

In the following section, we derive a nexus formulation from the DIH case that allows us to simply derive \eqref{measurized_AROE} and \eqref{measurized_LPrho}.

\subsection{The equilibrium problem \textcolor{black}{(the steady-state technique)} \label{sec:connection_discounted}}

 In this section, we establish an intuitive connection between the DIH and AR discounted cases, arising from seeking to find the optimal steady state within the discounted MDP framework. This rationale is analogous to the {\it steady-state technique} used to solve the AR problem in deterministic MDPs \citep{lerma_AR}. 
%
More specifically, the DIH optimality equations \eqref{measurized_DCOE}, with a discount factor $\alpha \in (0,1)$, can be expressed as

\begin{align*}
\overline{V}_\alpha^*(\nu)=\max_{\varphi \in \Phi} &\ \left\{ \ovr(\nu,\varphi) + \alpha   \overline{V}_\alpha^*(\nu')\right\} \qquad \qquad \forall \nu \in \mMpS 
\\
\mbox{s.t. } &\ \nu'=F(\nu,\varphi). \nonumber
\end{align*}
\textcolor{black}{Here, we index the measurized DIH value function by $\alpha$ to make explicit its dependence on the discount factor, and retain this notation henceforth when working in the AR setting}. A natural next step is to consider this problem in steady state, where $\nu=\nu'$, i.e.
\begin{align*}
(1-\alpha)\overline{V}_\alpha^*(\nu)=\sup_{\varphi \in \Phi} &\ \ovr(\nu,\varphi) \qquad \qquad \qquad \forall \nu \in \mMpS  
\\
\mbox{s.t. } &\ \nu=F(\nu,\varphi). \nonumber
\end{align*}
The constraint $\nu=F(\nu,\varphi)$ indicates that the current state distribution $\nu$ is a solution to the fixed point equation given by $F(\cdot,\varphi)$, where $\varphi \in \Phi$ is the implemented decision rule, and therefore we can consider the Markov chain with transition given by $F(\cdot,\varphi)$ to be in equilibrium. We then aim to find the best steady state in the MDP by solving the equilibrium problem,

\begin{flalign}
\sup_{\substack{\nu \in \mMpS \\ \varphi \in \Phi}} & \ \ \ovr(\nu,\varphi)  \label{equilibrium} \tag{E}\\
\mbox{s.t.}& \ \ \nu=F(\nu,\varphi). \nonumber
\end{flalign}

We highlight the bilinear structure of Problem \eqref{equilibrium}, which is linear in the state distribution \(\nu\) when fixing the policy \(\varphi\), and vice versa. This bilinearity may facilitate an optimization approach that systematically finds equilibria by iteratively fixing one variable while optimizing the other. 


\subsection{Deriving the measurized optimality equations from the equilibrium problem \label{sec:M-AROE}}

\begin{center}
\vspace{-0.5cm}
\begin{figure}[h!]
\centering
\begin{tabular}{c}
\includegraphics[width=250pt]{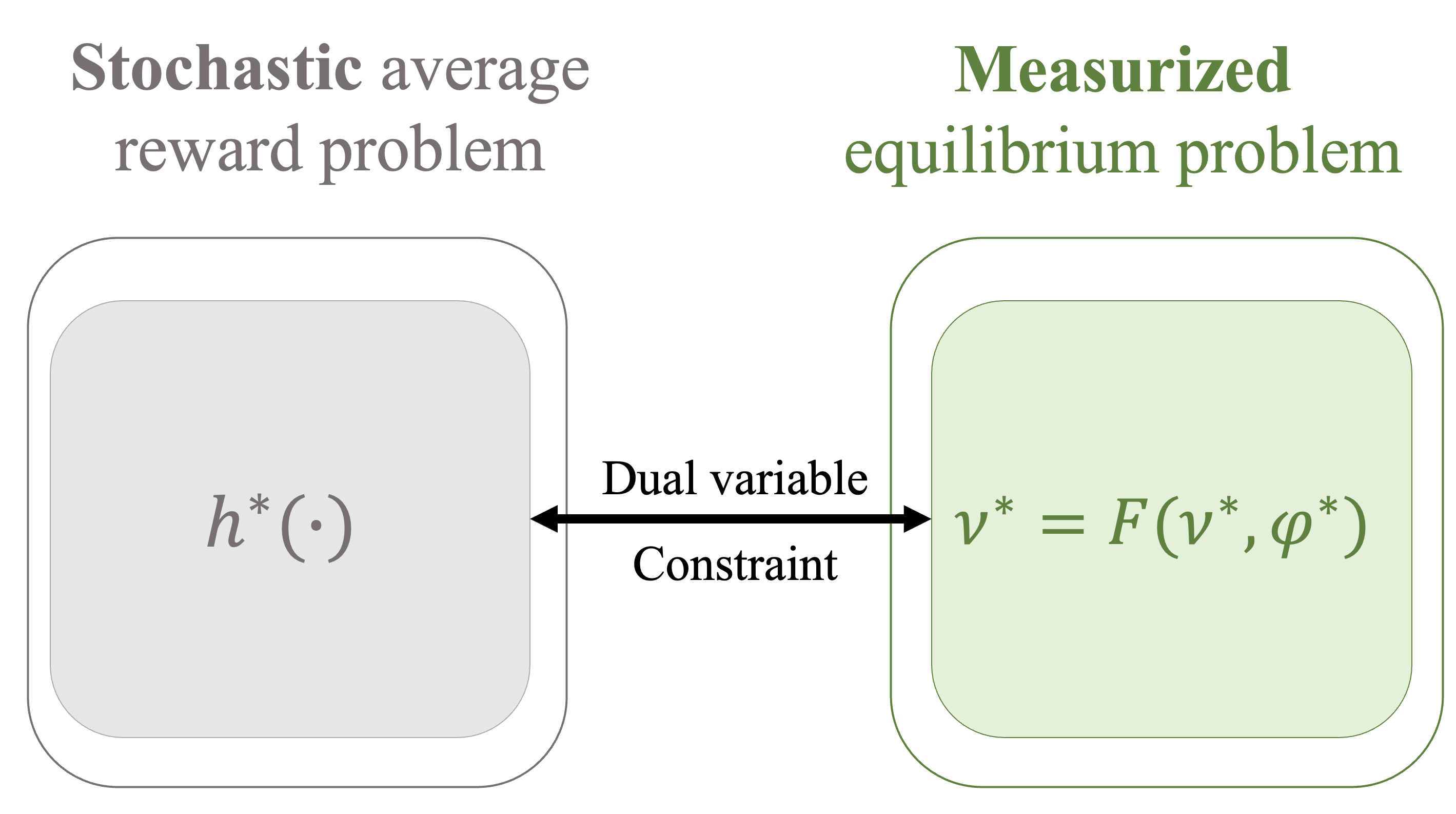}\\
\end{tabular}
\caption{\small Relationship between bias function $h(\cdot)$ in \eqref{AROE} and the equilbrium constraint in \eqref{equilibrium} \label{fig:h_dualvariable}}
\end{figure}
\vspace{-0.5cm}
\end{center}

In this section we derive the measurized optimality equations from the equilibrium problem and show their validity\footnote{In Appendix \ref{sec:measurizedAROE}, we provide an alternative derivation of the optimality equations using the lifting procedure outlined in Section \ref{sec:measurizing}.}. 
We start noting that the $\ovh(\cdot)$ terms in \eqref{measurized_AROE} cancel out for state distributions and decision rules that are in equilibrium, i.e., for those $\varphi \in \Phi$ such that there exists a state distribution $\nu_\varphi=F(\nu_\varphi,\varphi)$. This implies that $\orho^*  = \sup_{\varphi \in \Phi(\nu)} \left\{ \ovr(\nu,\varphi) + \ovh^*(F(\nu,\varphi)) -\ovh^*(\nu)\right\} \geq \ovr(\nu_\varphi,\varphi)$ for all equilibria $(\nu_\varphi,\varphi)$. As a consequence, $\orho^*$ gives an upper bound for \eqref{equilibrium}; i.e., $\orho^* \geq  \sup_{\varphi \in \Phi} \ovr(\nu_\varphi,\varphi)$, suggesting that the measurized optimality equations might be related to the dual of \eqref{equilibrium}. In what follows, we characterize the Lagrangian function and the Lagrangian dual formulation associated with the equilibrium problem. 

First, we rewrite the equilibrium constraint as $F(\nu,\varphi)(\mA)-\nu(\mA)$ for all $\mA \in \mBS$. \textcolor{black}{Assume that the Lagrange multiplier $h: \mS \to \mathbb{R}$ associated with this constraint exists, and write the Lagrangian function for the equilibrium problem as} 

\begin{align*} 
L(\nu,\varphi,h):=\ovr(\nu, \varphi)+\int_\mS h(s) dF(\nu, \varphi)(s) - \int_\mS h(s) d\nu(s).
\end{align*}
Now denote $\ovh(\nu):=\int_\mS h(s) d\nu(s)$ and rewrite the Lagrangian function above as

\begin{align} \label{L}
L(\nu,\varphi,\ovh):=\ovr(\nu, \varphi)+\ovh(F(\nu, \varphi)) - \ovh(\nu).
\end{align}
Then the Lagrangian dual of the equilibrium problem reads

\begin{align} \label{Lagrangian_equilibrium} \tag{L$_E$}
\inf_{\ovh(\cdot)} \sup_{\substack{\nu \in \mMpS\\ \varphi \in \Phi}} \ L(\nu,\varphi,\ovh),
\end{align}
and its optimal value gives an upper bound to the optimal reward in equilibrium. The Lagrangian problem \eqref{Lagrangian_equilibrium} is closely related to the measurized LP formulation associated with \eqref{measurized_AROE}. Indeed, the optimal solution $\orho^*$ to \eqref{measurized_LPrho} satisfies the constraints therein, thus $\orho^*\geq \sup_{\{ \nu \in \mMpS, \varphi \in \Phi\}} \ L(\nu,\varphi,\ovh)$. However, \eqref{measurized_LPrho} is infimizing over $\ovh$, thus intuitively yielding \eqref{Lagrangian_equilibrium}. Thereafter, we explore in more detail conditions under which the optimal measurized gain can be retrieved by solving these problems.

Now, we are interested in analyzing when there is no duality gap between \eqref{equilibrium} and \eqref{Lagrangian_equilibrium}, and understanding the role of the measurized optimality equations, especially since the Lagrangian function \eqref{L} resembles the right side of \eqref{measurized_AROE}. The following proposition links 
triplets \((\nu^*, \varphi^*, \ovh^*)\) solving \eqref{measurized_AROE} to saddle points that satisfy

\begin{equation} \label{L_sandwich}
L(\nu, \varphi, \oh^*) \leq L(\nu^*, \varphi^*, \ovh^*) \leq L(\nu^*, \varphi^*, \ovh),
\end{equation}
providing conditions for strong duality between \eqref{Lagrangian_equilibrium} and \eqref{equilibrium}.

\begin{proposition} \label{lemma:L_E}
Let $(\orho^*,\ovh^*,\of^*)$ be a 
triplet solving \eqref{measurized_AROE} such that $\of^*(\nu)=\varphi^*$ for all $\nu \in \mMpS$. Assume that there exists a $\nu_{\varphi^*} \in \mMpS$ that is invariant with respect to $\varphi^*$. Then
\begin{itemize}
\item[(a)] The function $L(\nu,\varphi^*,\ovh^*)$ is constant for all $\nu \in \mMpS$ and equal to $\orho^*=\ovr(\nu_{\varphi^*},\varphi^*)$
\item[(b)] $(\nu_{\varphi^*},\varphi^*,\ovh^*)$ is a saddle point satisfying \eqref{L_sandwich}
\item[(c)] There is no duality gap between \eqref{Lagrangian_equilibrium} and \eqref{equilibrium} and strong duality holds
\item[(d)] $(\nu_{\varphi^*},\varphi^*)$ is an optimal equilibrium.
\end{itemize}
\end{proposition}

\proof{

We first show (a). Since $(\orho^*,\ovh^*,\varphi^*)$ \textcolor{black}{solves \eqref{measurized_AROE}}, we have that

\begin{equation} \label{M-AROE_triplet}
\orho^* = \ovr(\nu,\varphi^*) + \ovh^*(F(\nu,\varphi^*))-\ovh^*(\nu)=L(\nu,\varphi^*,\ovh^*)  \qquad \forall \nu \in \mMpS.
\end{equation}
In particular, the equation above holds for $\nu_{\varphi^*}$, thus yielding that $\orho^*=\ovr(\nu_{\varphi^*},\varphi^*)$ because the bias-function terms cancel out. 

Now we demonstrate (b) and (c) by showing that \eqref{L_sandwich} holds.  First note that 

$$L(\nu_{\varphi^*},\varphi^*,\ovh^*) = L(\nu_{\varphi^*},\varphi^*,\ovh) \qquad \forall \ovh$$ 
since $\ovh(\nu_{\varphi^*})=\ovh(F(\nu_{\varphi^*},\varphi^*))$ for any $\ovh$ and thus these terms cancel out. Second, we show that $L(\nu,\varphi,\oh^*) \leq L(\nu_{\varphi^*},\varphi^*,\ovh^*)$ for all $\nu \in \mMpS$ and $\varphi \in \Phi(\nu)$. Indeed

\begin{align*}
L(\nu,\varphi,\oh^*) & = \ovr(\nu,\varphi) + \ovh^*(F(\nu,\varphi))-\ovh^*(\nu) \\
& \leq \max_{\varphi \in \Phi(\nu)} \left\{ \ovr(\nu,\varphi) + \ovh^*(F(\nu,\varphi))-\ovh^*(\nu) \right\}\\
&= \ovr(\nu,\varphi^*) + \ovh^*(F(\nu,\varphi^*))-\ovh^*(\nu) \\
&=\orho^*\\
&=\ovr(\nu_{\varphi^*},\varphi^*)\\
&=L(\nu_{\varphi^*},\varphi^*,\ovh^*),
\end{align*}
which shows (b). Finally, since $(\orho^*,\ovh^*,\varphi^*)$ satisfies \eqref{measurized_AROE}, we get that

\begin{align*} 
\ovr(\nu_{\varphi^*},\varphi^*)&=\ovr(\nu,\varphi^*)-\ovh^*(F(\nu,\varphi^*)) +\ovh^*(\nu) & \forall \nu \in \mMpS \\
&=L(\nu,\varphi^*,\ovh^*) & \forall \nu \in \mMpS \\
&= \sup_{ \varphi \in \Phi}  \ L(\nu,\varphi,\ovh^*)& \forall \nu \in \mMpS \\
&\geq \inf_{\ovh(\cdot)}  \sup_{ \varphi \in \Phi}  \ L(\nu,\varphi,\ovh)& \forall \nu \in \mMpS \\
& = \inf_{\ovh(\cdot)}  \sup_{\substack{\nu \in \mMpS\\ \varphi \in \Phi}}  \ L(\nu,\varphi,\ovh), & 
\end{align*}
but we had established before that $\ovr(\nu_{\varphi^*},\varphi^*)\leq \inf_{\ovh(\cdot)}  \sup_{\nu,\varphi}  \ L(\nu,\varphi,\ovh)$ because they were dual problems. As a consequence, we get (c) and (d). \hfill $\blacksquare$
}

This result allows us to interpret the optimality equations as the search for a saddle point of the Lagrangian function corresponding to the equilibrium problem. One can see $\ovh(\cdot)$ as quantifying the gain in terms of expected reward when deviating from the equilibrium $(\nu_{\varphi^*},\varphi^*)$; i.e.,

\begin{equation}\label{H*}
\ovh^*(F(\nu,\varphi^*))-\ovh^*(\nu)=\ovr(\nu_{\varphi^*},\varphi^*)-\ovr(\nu,\varphi^*).
\end{equation}

This aligns with the interpretation of $h^*(\cdot)$ as dual variable associated with the equilibrium constraint\footnote{This resembles the role of the {\it storage function} $\lambda$ in the ``dissipative" assumption adopted by \cite{lerma_AR} to show that the steady-state technique applies to a deterministic MDP. In particular, $\lambda$ is is required to satisfy $\lambda(F(\nu,\varphi))-\lambda(\nu)\leq \ovr(\nu^*,\varphi^*)-\ovr(\nu,\varphi)$.}. \textcolor{black}{The challenge we address in subsequent sections is that of showing that a 
triplet $(\orho^*,\ovh^*,\of^*)$ solving \eqref{measurized_AROE} exists for the lifted MDP and how it is connected with a triplet  $(\rho^*,h^*,f^*)$ solving \eqref{AROE}. }



\subsection{Connection between stochastic and measurized solutions \label{sec:measurizedAR}} 

\begin{center}
\vspace{-0.1cm}
\begin{figure}[H]
\centering
\begin{tabular}{c}
\includegraphics[width=250pt]{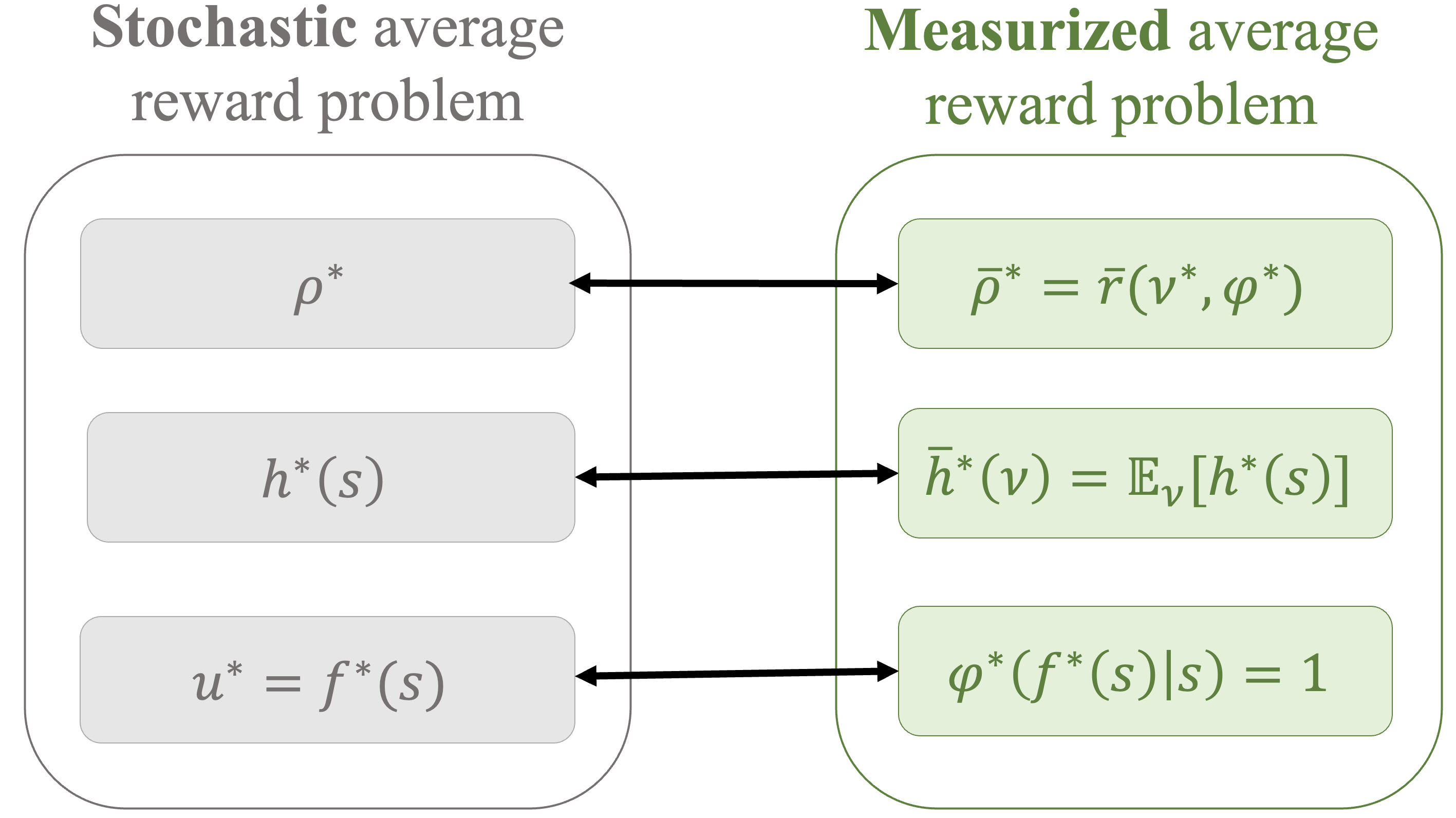}\\
\end{tabular}
\caption{Connection between stochastic and triplets solving \eqref{AROE} and \eqref{measurized_AROE}, respectively. \label{fig:canonical_triplets}}
\end{figure}
\vspace{-0.5cm}
\end{center}

In this section, we connect stochastic and measurized AR MDPs. First, we observe that Proposition \ref{prop:distribution} applies beyond the DIH case, meaning the states of the measurized AR MDP are also tied to the Ionescu-Tulcea probability. We then establish a connection between the optimal gain, bias function, and policy. The following proposition demonstrates that if the stochastic MDP can be solved via \eqref{AROE}, its measurized counterpart can similarly be solved via \eqref{measurized_AROE}. It also links the stochastic and measurized triplets \textcolor{black}{solving these optimality equations}, showing that the optimal gains are identical and equal to the optimal reward in equilibrium. Furthermore, it demonstrates that an optimal equilibrium is characterized by the AR-optimal decision rule and its invariant measure.

\begin{theorem}\label{th:measurized_AROE}
Let $(\rho^*,h^*,f^*)$ solve \eqref{AROE} \textcolor{black}{and assume that $h^*$ is integrable. Denote by $\varphi_{f^*}$ the measurized action concentrated around $f^*$ and let $\nu_{f^*}$ be an invariant distribution with respect to $\varphi_{f^*}$}. Define

\begin{equation}\label{measurized_h*}
\oh^*(\nu):=\int_{\mS} h^*(s) d\nu(s).
\end{equation}
Then:
\begin{itemize}
\item[(a)] $(\rho^*,\oh^*,\varphi_{f^*})$ solves the measurized optimality equation \eqref{measurized_AROE} with $\Phi(\nu)=\Phi$ for all $\nu\in \mMpS$.
\item[(b)] $(\nu_{f^*},\varphi_{f^*})$ is an optimal equilibrium
\item[(c)] The optimal stochastic and measurized AR value functions coincide and are equal to the optimal reward in equilibrium; i.e., 

\begin{equation}\label{Jopt}
\oJ^*(\nu)=\overline{\rho}^*=\rho^*=\ovr(\nu_{f^*},\varphi_{f^*}) \qquad \forall \nu \in \mMpS.
\end{equation}

\end{itemize}

\end{theorem}

\proof{  First we show (a). Since $(\rho^*,h^*,f^*)$ solves \eqref{AROE} we have that 

\begin{align*}
\rho^* +h^*(s) &= r(s,f^*(s)) + \mathbb{E}_Q [ h^*(s')| s, f^*(s)] & \forall s \in \mS.  \nonumber
\end{align*}
Integrating on the left and right-hand sides yield 

\begin{align*}
\rho^* + \oh^*(\nu) &=\rho^* + \int_{\mS} h^*(s) d\nu(s) & \forall \nu \in \mMpS \\
&= \int_{\mS} \left[ r(s,f^*(s)) + \int_{\mS}   h^*(s') Q(ds'| s, f^*(s)) \right] d\nu(s)&  \forall \nu \in \mMpS \\
&= \int_{\mS} \int_{\mU}  \left[ r(s,u) +  \int_{\mS} h^*(s') Q(ds'| s, u)\right] \varphi_{f^*}(du|s) d\nu(s)&  \forall \nu \in \mMpS \\
&=\ovr(\nu,\varphi_f)+\oh^*(F(\nu,\varphi_{f^*})) &  \forall \nu \in \mMpS,
\end{align*}
thus having that $(\rho^*,\oh^*,\varphi_{f^*})$ solves the second optimality equation in \eqref{measurized_AROE}. To show that it also solves the first equation, it suffices to prove that

\begin{align*}
\rho^* + \oh^*(\nu) & \geq \sup_{\varphi \in \Phi} \left\{ \ovr(\nu,\varphi) + \ovh^*(F(\nu,\varphi)) \right\} &  \forall \nu \in \mMpS.
\end{align*}
To show this, we integrate on the left and right-hand sides of the first \eqref{AROE} equation
\begin{align*}
\rho^* + \oh^*(\nu) &=\int_{\mS} \left\{ \rho^* + h^*(s)\right\} d\nu(s) & \forall \nu \in \mMpS \\
&\geq \int_{\mS} \sup_{\mU(s)}\left\{ r(s,u) + \int_{\mS}   h^*(s') Q(ds'| s, u) \right\} d\nu(s)&  \forall \nu \in \mMpS \\
&= \int_{\mS} \left\{ r(s,f^*(s)) + \int_{\mS}   h^*(s') Q(ds'| s, f^*(s)) \right\} d\nu(s)&  \forall \nu \in \mMpS \\
&= \int_{\mS} \int_{\mU}  \left[ r(s,u) +  \int_{\mS} h^*(s') Q(ds'| s, u)\right] \varphi_{f^*}(du|s) d\nu(s)&  \forall \nu \in \mMpS \\
&= \max_{\varphi \in \Phi} \int_{\mS} \int_{\mU}  \left[ r(s,u) +  \int_{\mS} h^*(s') Q(ds'| s, u)\right] \varphi(du|s) d\nu(s)&  \forall \nu \in \mMpS, \\
\end{align*}
where the inequality comes from $(\rho^*,h^*,f^*)$ solving \eqref{AROE}. The maximum in the last line comes from the fact that $f^*$ is attainable, and thus $\varphi_{f^*}$ is. The equality in this last line is due to $f^*$ maximizing the right-hand side of the optimality equations, thus having a strict inequality would contradict the optimality of $f^*$.

Now we show (b) and (c). It suffices to integrate with respect to $\nu_{f^*}$ and $\varphi_{f^*}$ on both sides of \eqref{AROE}

\begin{align*}
\orho^*&=\int_{\mS} \int_{\mU} \rho^* \  \varphi_{f^*}(du|s) d\nu_{f^*}(s) \\
& = \int_{\mS} \int_{\mU} r(s,f^*(s))  \varphi_{f^*}(du|s) d\nu_{f^*}(s)+ \int_{\mS} \int_{\mU}  \int_\mS h^*(s')Q(ds'|s,f^*(s))  \varphi_{f^*}(du|s) d\nu_{f^*}(s) \\
& \qquad - \int_{\mS} \int_{\mU}  h^*(s) \varphi_{f^*}(du|s) d\nu_{f^*}(s) \\
&= \ovr(\nu_{f^*},\varphi_{f^*}) + \ovh^*(F(\nu_{f^*},\varphi_{f^*}) ) -\ovh^*(\nu_{f^*} )\\
&= \ovr(\nu_{f^*},\varphi_{f^*}) 
\end{align*}

As we discussed at the beginning of the section $\orho^*\geq \ovr(\nu^*,\varphi^*)$; therefore $(\nu_{f^*},\varphi_{f^*}) $ is an optimal equilibrium. \hfill $\blacksquare$
}

As a consequence of the previous theorem and Proposition \ref{lemma:L_E}, we can claim that if there exists \textcolor{black}{a triplet $(\rho^*, h^*, f^*)$ solving \eqref{AROE}} and there exists an invariant measure $\nu_{f^*}$, then strong duality holds between the equilibrium problem and its Lagrangian dual. Moreover, we can establish a connection between \textcolor{black}{solutions of \eqref{AROE}}, optimal equilibria and saddle points of the Lagrangian function \eqref{L}. In addition, Theorem \ref{th:measurized_AROE} also allows us to intuitively analyze the implications of working in the unichain versus multichain case. For \(J^*(\cdot)\) to be unique, either (a) we are in the unichain case, where there is at most one equilibrium measure \(\nu_f\) for any deterministic rule \(\varphi_f\), and it is reachable regardless of the initial state distribution, or (b) we are in a special multichain case, where, for a fixed \(\varphi_f\), all equilibrium points \((\nu_f, \varphi_f)\) yield the same revenue. In both scenarios, we can find an optimal solution to \eqref{equilibrium} that is reachable from the initial distribution \(\nu_0\). However, if we fall outside cases (a) and (b) into the general multichain case, \eqref{equilibrium} may produce an optimal solution \((\nu^*, \varphi^*)\) that is not necessarily reachable from every starting distribution \(\nu_0\) and $J^*(\cdot)$ may not be constant.

\section{Translation of further classical techniques to the measurized AR problem \label{sec:T1-T3}}

While we have demonstrated how the steady-state technique can be extended to solve standard stochastic MDPs and shown that the lifted MDP can be solved through \eqref{measurized_AROE} if the original MDP can be solved via \eqref{AROE}, there are additional techniques for addressing AR problems:

\begin{enumerate}
\item[(T1)] {\it The AR optimality inequalities}: Deriving optimal decisions by solving equations such as \eqref{AROE}, where the equality is replaced with \( \leq \).
\item[(T2)] {\it The LP formulation}: Obtaining the AR-optimal value function by solving an LP formulation.
\item[(T3)] {\it The vanishing discount approach \citep[Chapter 5.3]{lerma}}: Solving the AR problem by taking the limit as the discount factor \( \alpha \to 1 \) in the DIH case.
\end{enumerate}

These approaches are well established for stochastic MDPs (see Chapters 5 and 6 in \cite{lerma}). \cite{lerma_AR} provides conditions under which T1 and T3 extend to deterministic MDPs. Although Section \ref{sec:MDP} presents the classical assumptions and theoretical results for the standard DIH case, we refer the reader to Appendices A.1 and A.2 for the corresponding assumptions and results for standard AR MDPs in order to keep the exposition self-contained without overburdening the main text. In this section, we demonstrate that if T1-T3 can be applied to the original MDP, they also extend to its lifted MDP. Notably, we achieve this through the steady-state technique, which, while widely recognized for deterministic MDPs \citep{lerma_AR}, has not been commonly applied to general stochastic MDPs. 

\subsection{Solvability of the measurized AR optimality inequalities (T1) }

In this section, we extend the inheritance property to the AR case. In particular we show that if the stochastic MDP satisfies the commonly known sufficient assumptions for solving the AR problem through the optimality inequalities (see Equations \eqref{AROI}, Assumption \ref{assumption541}, Lemma \ref{lemma:vanishing_discount} and Theorem \ref{theorem543} in Appendix A.1), then the measurized MDP inherits these properties. Consequently, in the lifted framework, the optimal measurized AR value function and policy can be obtained by solving the measurized optimality inequalities. Using the definition of measurized MDP, we can write the measurized AROI as

\begin{align}\label{measurized_AROI}\tag{$\mM$-AROI}
\overline{\rho}^* +\ovh^*(\nu) & \leq \sup_{\varphi \in \Phi(\nu)} \left\{ \ovr(\nu,\varphi) + \ovh^*(F(\nu,\varphi)) \right\} \qquad \forall \nu \in \mMpS \\
\orho^* +\ovh^*(\nu) &\leq \ovr(\nu,\overline{f}(\nu)) + \ovh^*(F(\nu,\of(\nu)))  \qquad \ \forall \nu \in \mMpS.  \nonumber\end{align}

The following lemma shows that Assumption \ref{assumption541} is inherited by the measurized MDP. 



\begin{lemma}\label{lemma:assumption541}
\textcolor{black}{Suppose all the assumptions in Proposition \ref{prop:assumptions} hold.} Then if the stochastic MDP satisfies Assumption \ref{assumption541}, its measurized counterpart also does.
\end{lemma}

\proof{
We define $\ovh_{\alpha} (\nu):=\oV_{\alpha}^* (\delta_z) - \oV_{\alpha}^* (\nu)$ and $\overline{b}(\nu):= \int_{\mS} b(s) d\nu(s)$.
\begin{itemize}
\item[(a)]Theorem \ref{th:V} implies that 
$$(1-\alpha)V_\alpha^*(z) =(1-\alpha)\oV_\alpha^*(\delta_z)\geq M \qquad \forall \alpha \in [\beta,1) $$
which proves part (a) of Assumption \ref{assumption541}.
\item[(b)] If $-b(s) \leq h_{\alpha}(s) \leq N$ for all $s \in \mS$ and $\alpha \in [\beta,1)$, then

\begin{align}
-\int_{\mS} b(s) d\nu(s)&\leq \int_{\mS} h_{\alpha} (s) d\nu(s) \leq \int_{\mS} N d\nu(s) &\qquad \forall \nu \in \mMpS. \label{expected_assumption}
\end{align}
Once again making use of Theorem \ref{th:V} we have for all $\nu \in \mMpS$
\begin{align*}
\overline{h}_{\alpha} (\nu)&=\int_{\mS}V_{\alpha}^* (z) d\nu(s)- \int_{\mS}V_{\alpha}^* (s) d\nu(s) \\
&=\int_{\mS} \left(V_{\alpha}^* (z) - V_{\alpha}^* (s) \right)d\nu(s) \\
&=\int_{\mS} h_{\alpha} (s) d\nu(s) .
\end{align*}
We can rewrite \eqref{expected_assumption} as
$$ -\overline{b}(\nu)  \leq \oV_{\alpha}^* (\delta_z) - \oV_{\alpha}^* (\nu)  \leq N \qquad \forall \nu \in \mMpS$$
hence having that part (b) of Assumption \ref{assumption541} is also satisfied. \hfill $\blacksquare$\\
\end{itemize}
}

The inheritance property ensures that there is a triplet $(\orho^*,\ovh^*,\of^*)$ solving \eqref{measurized_AROI} with $\of^*$ being an AR-optimal selector and $\orho^*=\oJ^*(\nu)$ for all $\nu \in \mMpS$.

\subsection{Solvability of the measurized AR LP formulation (T2)}


In Section \ref{sec:dualizing}, we showed that \eqref{measurized_DCOE} can be obtained by dualizing the LP formulation of \eqref{DCOE} and performing the change of variables \eqref{change_variables}. Following a similar approach, the next proposition establishes that the equilibrium problem is the dual of the AR LP formulation (see Equation \eqref{LP_rho} in Appendix A.2). This result deepens the intuition behind the LP theory in the AR context \cite[Chapter 6]{lerma} and shows that technique T2 applies in the measurized framework. Furthermore, it ensures solvability under different assumptions than those in \cite{lerma}.

\begin{proposition}\label{prop:Edual}
The equilibrium problem \eqref{equilibrium} is the dual problem of the standard \textcolor{black}{AR LP formulation} \eqref{LP_rho} after performing change of variables \eqref{change_variables}.
\end{proposition}

\proof{

We start by associating a measure $\gamma \in \mM_+(\mS\times\mU)$ with the constraints in \eqref{LP_rho}, and formulate its dual problem as

\begin{align*}
\sup_{\gamma \in \mM_+(\mS\times\mU)} & \ \int_\mS \int_{\mU(s)} r(s,u) d\gamma(s,u) \\
\mbox{s.t.} &\  \int_\mS \int_{\mU(s)} d\gamma(s,u)=1\\
&\  \int_{\mA} \int_{\mU(s)} d\gamma(s,u) - \int_\mS \int_{\mU(\mS)} Q(\mA|s,u) d\gamma(s,u)=0 \qquad \forall \mA \in \mB(\mS)
\end{align*}
where the first constraint is associated to variable $\rho$ and the second to $h(\cdot)$. Moreover, the first constraint states that $\gamma$ is actually a probability measure, and the second constraint says that $\gamma$ is invariant through the stochastic kernel $Q$, so we can rewrite the previous problem as 

\begin{align*}
\sup_{\gamma \in \mM_\mathbb{P}(\mS\times\mU)} & \ \int_\mS \int_{\mU(s)} r(s,u) d\gamma(s,u) \\
\mbox{s.t.} &\ \  \gamma(\mA \times \mU(\mA))= \int_\mS \int_{\mU(\mS)} Q(\mA|s,u) d\gamma(s,u) \qquad \forall \mA \in \mB(\mS)
\end{align*}
Proposition \ref{prop:change_variables_t} ensures that we can perform the change of variables \eqref{change_variables} without loss of optimality, having

\begin{align*}
\sup_{\substack{\nu \in \mM_\mathbb{P}(\mS)\\\varphi \in \Phi}} & \ \int_\mS \int_{\mU(s)} r(s,u) \varphi(du|s) d\nu(s) \\
\mbox{s.t.} &\ \  \nu(\mA)= \int_\mS \int_{\mU(\mS)} Q(\mA|s,u)  \varphi(du|s) d\nu(s) \qquad \forall \mA \in \mB(\mS)
\end{align*}
which using the definition of $\ovr$ and $F$ can be rewritten as \eqref{equilibrium}. \hfill $\blacksquare$\\
}

Although the linear program \eqref{LP_rho} is the dual of \eqref{equilibrium}, there might be a duality gap. The following result shows that this gap is indeed zero\footnote{A similar result is shown in Chapter 6 of \cite{lerma} under different assumptions.}.

\begin{corollary} \label{lemma:strong_duality}
Consider Assumptions \ref{assumption421} and \ref{assumption551}. If there exists and invariant measure $\nu_{f^*}$ for the decision rule $\varphi_{f^*}$ concentrated around an AR-optimal selector $f^*$ of the stochastic MDP, then the optimal objective values of \eqref{LP_rho} and \eqref{equilibrium} coincide and strong duality holds.
\end{corollary}

\proof{

Let $\rho'$ denote the optimal objective of \eqref{LP_rho}. Since this problem is the dual of \eqref{equilibrium}, then its optimal solution is larger or equal to the optimal reward in equilibrium; i.e., 

\begin{equation}\label{rho>E}
\rho'\geq \ovr(\nu^*,\varphi^*). 
\end{equation}

However, under Assumptions \ref{assumption421} and \ref{assumption551}, Theorem \ref{theorem554} guarantees the existence of a constant $\rho^*$, a function $h^*(\cdot)$ and a selector $f^* \in \mathbb{F}$ solving the \eqref{AROE} such that $\rho^*=J^*(s)$ for all $s \in \mS$. The triplet $(\rho^*,h^*,f^*)$ is then a feasible solution to \eqref{LP_rho}, having that $\rho^*\geq \rho'$. However, Theorem \ref{th:measurized_AROE}.(c) entails that $\rho^*=\ovr(\nu^*,\varphi^*)$, thus yielding

\begin{equation}\label{rho<E}
\rho'\leq \ovr(\nu^*,\varphi^*). 
\end{equation}

Therefore, \eqref{rho>E} and \eqref{rho<E} yield $\rho'=\rho^*=\ovr(\nu_{f^*},\varphi_{f^*})$. Since $(\nu_{f^*},\varphi_{f^*})$ exists and it is feasible to \eqref{equilibrium}, strong duality holds. \hfill $\blacksquare$\\
}

As a consequence of this result and Theorem \ref{th:measurized_AROE} we have that under these assumptions we can construct a feasible solution to the measurized LP formulation  \eqref{measurized_LPrho}. Moreover, the optimal value for \eqref{LP_rho} coincides with the optimal value of \eqref{measurized_LPrho}, thus ensuring the solvability of T2 in the measurized framework. The practical relevance of these results lies in the potential of the LP approach for Approximate Dynamic Programming, commonly used for high-dimensional MDPs. For instance, the bias function $\ovh(\cdot)$ can be approximated via a linear combination of basis functions $\overline{\phi}_k(\nu)$, incorporating moment-based or divergence-based measures; see Section \ref{sec:approximations}.

\subsection{Validating and extending the vanishing discount technique (T3)}

In this section, we partially resolve the open question in \cite{lerma_AR} regarding the conditions under which the vanishing discount leads to AR-optimality equations in deterministic MDPs, though our answer pertains specifically to measurized MDPs. Moreover, while the traditional vanishing discount technique focuses on the convergence of $V^*_\alpha(\cdot)$ to $J^*(\cdot)$ as \(\alpha \to 1\), the measurized framework naturally extends this to include optimal policies and their invariant measures. Specifically, if \(\varphi^*_{\alpha}\) is the optimal decision rule for the \(\alpha\)-DIH MDP and \(\nu_{\alpha}^*\) is its best invariant measure (i.e., \(\nu_{\alpha}^*\) solves the equilibrium problem with decision rule fixed to \(\varphi^*_{\alpha}\)), then  \((\nu_{\alpha}^*, \varphi^*_{\alpha})\) converges to an optimal equilibrium \((\nu^*, \varphi^*)\) as \(\alpha \to 1\). Consequently, the measurized AR problem can also be solved via its DIH optimality equations as \(\alpha\) approaches 1. Interestingly, we rely on the equilibrium problem \eqref{equilibrium} to establish these findings. Figure \ref{fig:vanishing_discount} provides a visual summary.

\begin{center}
\begin{figure}[h!]
\centering
\begin{tabular}{c}
\includegraphics[width=250pt]{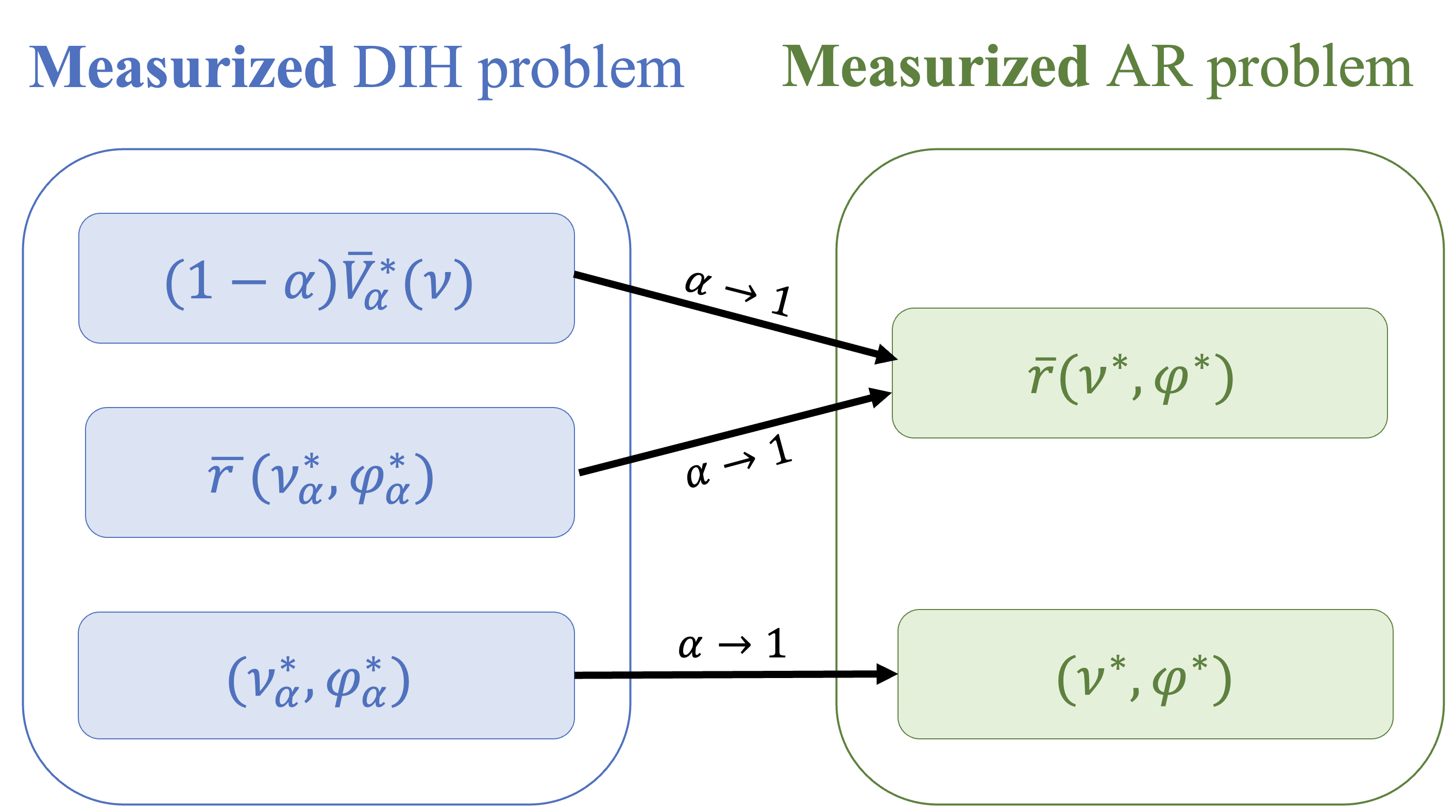}\\
\end{tabular}
\caption{Connection between measurized DIH, AR and equilibrium solutions as $\alpha \to 1$. \label{fig:vanishing_discount}}
\end{figure}
\vspace{-1cm}
\end{center}

\subsubsection{Convergence of discounted-optimal measurized value function as $\alpha \to 1$}

In this section, we replicate the traditional result that ensures the convergence of the discounted-optimal value functions to the optimal gain (and hence the optimal value of \eqref{equilibrium}) as $\alpha \to 1$. We extend this result and show that the rewards of the equilibria induced by the DIH-optimal policies also converge to the optimal gain. To show the latter we will make use of the following definition.

\begin{definition}\label{discounted_equilibria}
Let \( f_{\alpha}^* \) denote the optimal selector for the DIH MDP with discount factor \( \alpha \). Define \( \varphi_{f_{\alpha}^*} \) as the decision rule concentrated around \( f_{\alpha}^* \), and \( \nu_{f_{\alpha}^*} \) as the best invariant measure associated with \( \varphi_{f_{\alpha}^*} \); i.e., \( \nu_{f_{\alpha}^*} \) solves the equilibrium problem \eqref{equilibrium} with the decision rule fixed to \( \varphi_{f_{\alpha}^*} \). We refer to the pair \( (\nu_{f_{\alpha}^*}, \varphi_{f_{\alpha}^*}) \) as the \(\alpha\)-discounted-optimal equilibrium.
\end{definition}

The following theorem shows that $ (1-\alpha)\oV^*_{\alpha}(\nu)$ and the reward of the $\alpha$-discounted-optimal equilibria, $\ovr(\nu_{f_{\alpha}^*}, \varphi_{f_{\alpha}^*})$, converge to the optimal value of \eqref{equilibrium} (and thus $\orho^*$, by virtue of Theorem \ref{th:measurized_AROE})  as $\alpha \to 1$ . \textcolor{black}{To show this, we need to adopt the usual assumptions for the DIH and AR MDPs, but also we need to adopt Assumption \ref{assumption:welldefined}(a), which is usually made to guarantee that the Policy Iteration Algorithm for the AR problem is {\it well defined} \citep{PIA}.}

\begin{theorem}\label{prop:equilibrium_AROI}
Let $(\nu^*,\varphi^*)$ be an optimal solution to the equilibrium problem \eqref{equilibrium}. Suppose that the stochastic MDP satisfies \textcolor{black}{all the assumptions in Proposition \ref{prop:assumptions}}, and also Assumptions \ref{assumption551} and \ref{assumption:welldefined}(a). Then:

\begin{enumerate}

%

\item[(a)] There exists a  nondecreasing sequence of discount factors $\alpha(n) \uparrow 1$ such that 
\begin{equation}\label{measurized_vanishing_discount}
\lim_{n \rightarrow \infty} (1-\alpha(n)) \oV^*_{\alpha(n)}(\nu)=\ovr(\nu^*,\varphi^*) \qquad \forall \nu \in \mMpS.
\end{equation}

\item[(b)] Under Assumption \ref{assumption:MCT} there also exists a sequence of discount factors $\alpha(n) \uparrow 1$ such that
\begin{equation*} 
\lim_{n \rightarrow \infty} \ovr\left(\nu_{f_{\alpha(n)}^*},\varphi_{f_{\alpha(n)}^*}\right)=\ovr(\nu^*,\varphi^*),
\end{equation*}
\end{enumerate}
\end{theorem}

\proof{
Lemma \ref{lemma:assumption541} states that if an MDP satisfies Assumption \ref{assumption541}, then the measurized MDP also does. Since the stochastic MDP follows Assumption \ref{assumption551}, it also follows Assumption \ref{assumption541}. As a consequence, Lemma \ref{lemma:vanishing_discount} and Theorem \ref{theorem543}  hold for the measurized MDP, that is to say, there exists a constant $\orho^*$ such that

\begin{equation}\label{eq:measurized_vanishing}
\lim_{n \rightarrow \infty} (1-\alpha(n)) \oV^*_{\alpha(n)}(\nu)=\orho^* \qquad \forall \nu \in \mMpS.
\end{equation}
Moreover, Theorem \ref{th:measurized_AROE} ensures that $\orho^*$ is the measurized AR-value function, i.e., 

$$\oJ^*(\nu) =\orho^*=\ovr(\nu^*,\varphi^*) \qquad \nu \in \mMpS,$$

\noindent which shows (a).

Now we prove (b). In Theorem \ref{prop:unique_phi}, we show that $\pi^*_\alpha$ does not depend on the initial state distribution $\nu$ and it is concentrated around the optimal stochastic policy; i.e., $\pi^*_\alpha=\{ \varphi_{f^*_\alpha} \}_{t \geq 0}$. We now evaluate the discounted-optimal measurized value function at the best invariant measure with respect to $\varphi_{f^*_\alpha}$

\begin{align*}
(1-\alpha(n))\oV^*_{\alpha}(\nu_{f^*_\alpha(n)})&=(1-\alpha(n))\sum_{t=0}^\infty \alpha(n)^t \ovr(\nu_{f^*_\alpha(n)},\varphi_{f^*_\alpha(n)})  &\quad \forall \alpha\in[0,1) \nonumber\\
&=\ovr(\nu_{f^*_\alpha(n)},\varphi_{f^*_\alpha(n)}) & \quad \forall \alpha\in[0,1) \\
&\leq\ovr(\nu^*,\varphi^*) \nonumber
\end{align*}
where the second equality comes from $(\nu_{f_\alpha(n)^*},\varphi_{f_\alpha(n)^*}) $ being an equilibrium, and the last inequality from $(\nu^*,\varphi^*) $ being an optimal equilibrium. Using (a), we have that 

\begin{align*}
\ovr(\nu^*,\varphi^*)&=\lim_{n\rightarrow \infty} (1-\alpha(n)) \oV_{\alpha(n)}^* \left(\nu_{f_{\alpha(n)}^*}\right)\\
&=\lim_{n\rightarrow \infty} \ovr\left(\nu_{f_{\alpha(n)}^*},\varphi_{f_{\alpha(n)}^*}\right) \\
&\leq \ovr(\nu^*,\varphi^*).
\end{align*}
Thus yielding (b). \hfill $\blacksquare$\\
}

\subsubsection{Convergence of discounted-optimal measurized policies and their steady-state distributions as $\alpha \to 1$} 

In this section, we further extend the traditional vanishing discount technique to show that the $\alpha$-discounted-optimal equilibria $\left(\nu_{f_{\alpha}^*}, \varphi_{f_{\alpha}^*}\right)$ converge to an optimal equilibrium $(\nu^*,\varphi^*)$ as $\alpha \to 1$.

\begin{proposition}
Consider 
all the assumptions made in Theorem \ref{prop:equilibrium_AROI}. In addition, assume that there exists a subsequence $\{\alpha(n)\}_{n \in \mathbb{N}}$, $\alpha(n)\rightarrow 1$, such that there exists a limiting kernel $\varphi_{1}^*:=\lim_{n\rightarrow \infty} \varphi_{f_{\alpha(n)}^*} \in \Phi$ and a limiting state distribution $\nu_1^*:=\lim_{n\rightarrow \infty} \nu_{f_{\alpha(n)}^*}$. Then the supremum in the equilibrium problem \eqref{equilibrium} is attainable and we can assume that the optimal equilibrium is $(\nu_1^*,\varphi_1^*)$ without loss of optimality.
\end{proposition}

\proof{

The proof follows this structure:
\begin{itemize}
\item[(i)] First, we prove that the limiting measurized state-action pair $(\nu_1^*,\varphi_1^*)$ yields the same revenue as any optimal equilibrium $(\nu^*, \varphi^*)$; i.e., we demonstrate that $\ovr(\nu_1^*,\varphi_1^*)=\ovr(\nu^*, \varphi^*)$.
\item[(ii)] Second, we show that $(\nu_1^*,\varphi_1^*)$ is a feasible solution to the equilibrium problem; i.e., we show that $\nu_1^*=F(\nu_1^*,\varphi_1^*)$.
\end{itemize}

\begin{itemize}
\item[(i)] We use Theorem \ref{prop:equilibrium_AROI}. Since \eqref{measurized_vanishing_discount} holds for any state distribution $\nu \in \mMpS$, in particular it holds for the invariant measure $\nu_{f_{\alpha(n)}^*}$, which verifies $\nu_{f_{\alpha(n)}^*}=F(\nu_{f_{\alpha(n)}^*},\varphi_{f_{\alpha(n)}^*})$. Therefore we have

\begin{align*}
\ovr(\nu^*,\varphi^*) &=\lim_{n \rightarrow \infty} (1-\alpha(n)) \oV^*_{\alpha(n)} (\nu_{f_{\alpha(n)}^*}) &\\
 & = \lim_{n \rightarrow \infty} \ovr(\nu_{f_{\alpha(n)}^*},\varphi_{f_{\alpha(n)}^*}) & \\\
&=  \lim_{n \rightarrow \infty}  \int_{\mS} \int_{\mU} r(s,u) \varphi_{f_{\alpha(n)}^*}(du|s) d\nu_{f_{\alpha(n)}^*}(s) &  \\
&=   \int_{\mS} \int_{\mU} r(s,u)  \lim_{n \rightarrow \infty} \varphi_{f_{\alpha(n)}^*}(du|s) d\nu_{f_{\alpha(n)}^*}(s) &  \\
&=   \int_{\mS} \int_{\mU} r(s,u)  \varphi^*_{1}(du|s) d\nu^*_{1}(s) &  \\
&=\ovr(\nu_1^*,\varphi_1^*) & \\
\end{align*}
where the third and sixth equalities comes from the definition of $\ovr$, the third equality comes from the MCT, the fourth arises from applying the algebra of limits. 

\item[(ii)] Now we prove the feasibility of $(\nu_1^*,\varphi_1^*)$, i.e., we prove $\nu_1^*$ is invariant for $\varphi_1^*$ by showing that
\begin{equation*}
\nu_1^*(\cdot)=\int_{\mS} \int_{\mU} Q(\cdot|s,u) \varphi_1^*(du|s) d\nu_1^*(s)
\end{equation*}
\end{itemize}
For all $n \in \mathbb{N}$ we have

\begin{align*}
\nu_{f_{\alpha(n)}^*}(\mA)-\int_{\mS} \int_{\mU} Q(\mA|s,u) \varphi_{f_{\alpha(n)}^*}(du|s) d\nu_{f_{\alpha(n)}^*}(s)&=0\\
\lim_{n\rightarrow \infty} \left\{ \nu_{f_{\alpha(n)}^*}(\mA)-\int_{\mS} \int_{\mU} Q(\mA|s,u) \varphi_{f_{\alpha(n)}^*}(du|s) d\nu_{f_{\alpha(n)}^*}(s)\right\} &= \lim_{n\rightarrow \infty} 0=0\\
\end{align*}
Because the limit of product of sequences is the product of the limits of the sequences; i.e., $\lim_{n \rightarrow \infty} a_n b_n =  \left(\lim_{n \rightarrow \infty} a_n \right) \left(\lim_{n \rightarrow \infty} b_n \right)$, the previous equality reads
\begin{align*}
\nu_{1}^*(\mA)-\int_{\mS} \int_{\mU} Q(\mA|s,u) \varphi_{1}^*(du|s) d\nu_{1}^*(s)&=0.
\end{align*}
Therefore, there exists a feasible solution $(\nu^*_1,\varphi_1^*)$ to the equilibrium problem attaining the optimal objective. Therefore, without loss of optimality we can assume that $\nu_1^*=\nu^*$ and $\varphi_1^*=\varphi^*$ and the supremum in \eqref{equilibrium} can be replaced by a maximum. \hfill $\blacksquare$\\


}


\section{Concluding remarks and extensions.}

In this paper, we have revisited the concept of general lifted MDPs in the space of measures, which was described, studied, and used in \cite{BSdeterministic} and Chapter 9 of \cite{bertsekas}. While these works laid the foundational theory for what we refer to as BS-deterministic MDPs, it has seen limited practical application since their introduction in 1979. 

\textcolor{black}{We make BS-deterministic MDPs more accessible by placing them} within the semicontinuous-semicompact framework of \cite{lerma}, focusing on Borel-measurable optimal value functions and policies rather than semianalytic or universally measurable ones. This approach, while seemingly less general, offers several benefits: (i) it allows for computable results; (ii) it relies on mild, easily verifiable assumptions suitable for most practical control problems; and (iii) it simplifies proofs. Unlike \cite{BSdeterministic,bertsekas}, which uses BS-deterministic MDPs to derive results for the original MDP in the universally measurable framework, we show that a measurized MDP inherits key assumptions from its baseline MDP in the semicontinuous-semicompact framework. This enables us to apply existing theory from \cite{lerma} directly, relying on simpler mathematics.

Another contribution of our work is the simplification of BS-deterministic MDPs, whose controls are joint state-action probability measures, constrained so that their marginals match the current state distribution. In contrast, the measurized MDPs introduced in Section \ref{sec:measurized_def} use actions that are Markov decisions in the original MDP. In Section \ref{sec:change_variables} we use the Radon-Nikodym derivative to show that these two models are equivalent. We further enhance the interpretability of the measurized framework by providing a sampling interpretation, which allows us to view the measurized value function as solving over infinite realizations of the standard (stochastic) MDP.

Beyond making the lifted framework more interpretable, we showcase its practical advantages. As demonstrated in Section \ref{sec:advantages}, the lifted framework facilitates the incorporation of probabilistic constraints and value function approximations, such as moment-based constraints and Conditional Value at Risk (CVaR). 

These were not considered or exploited in the original works \cite{bertsekas,BSdeterministic}, highlighting the untapped potential of this approach. In addition, in Sections \ref{sec:dualizing} and \ref{sec:measurizing} we contribute to the practicality of lifted MDPs by introducing two alternative lifting methods that are both easy to implement and understand. 

We illustrate how the latter procedure reveals that non-deterministic measure-valued MDPs arise when lifting MDPs whose transitions are impacted by external random shocks.

In Section \ref{sec:Measurized_AR}, we demonstrate how the measurized framework unifies and extends several ideas previously explored in the AR literature. Specifically, we show that the steady-state technique, traditionally applied only to deterministic MDPs, can also solve the stochastic AR problem within the measurized framework. This perspective connects with the concept of minimum pairs, reframing it from an analysis technique to a direct solution method. Additionally, we establish that the measurized optimality equations are the Lagrangian dual of \eqref{equilibrium}. This generalizes the interpretation of bias functions as dual variables associated with equilibrium constraints from deterministic MDPs in real spaces \citep{flynn} to stochastic MDPs on Borel spaces. Consequently, the measurized framework not only provides a new approach to analyzing AR stochastic MDPs but also broadens the applicability of these foundational techniques.

Finally, Section \ref{sec:T1-T3} shows how additional tools for analyzing AR problems (namely, the AR optimality inequalities, the LP formulation, and the vanishing discount method) extend to the measurized framework. Some of these techniques also admit a broader interpretation through the steady-state approach. In particular, whereas the classical vanishing discount method studies the convergence of the optimal DIH value function to the optimal AR value function as the discount factor approaches one, the measurized framework naturally extends this perspective to include optimal policies and their associated invariant state distributions.

For future applications, the measurized framework shows significant potential for managing populations. The sampling interpretation presented in Corollary \ref{prop:asymptotics1} demonstrates that measurized MDPs are particularly suitable when the population size is sufficiently large. More theoretical works may focus on extending the measurized perspective to other frameworks within the MDP field. 
Finally, the equilibrium problem that we discuss in this paper is a bilinear program: it is linear on the state distribution $\nu$ when we fix the policy $\varphi$, and vice-versa. This structure is new in the MDP context, and a Primal-Dual Decomposition Algorithm could be designed to address the resolution of such an infinite-dimensional bilinear program. Since the equilibrium problem is deterministic, one could use (infinite-dimensional) optimization techniques to solve it.

\section*{Appendix}
\appendix
\section{\label{appendix:lerma} Additional definitions and theoretical results from \cite{lerma} }

\begin{definition}[sup-compact]
A function $g:\mathbb{K} \to  \mathbb{R}$ is said to be upper sup-compact on $\mathbb{K}$ if the set

$$\{ u \in \mU(s): \ g(s,u)\geq c\}$$ 
is compact for every $s \in \mS$ and $c\in \mathbb{R}$.
\end{definition}

\begin{definition}[upper semicontinuous] \label{def:usc}
Let $\mX$ be a metric space and $g:\mX \to \mathbb{R} \cup \{+\infty\}$ a function such that $g(x)<\infty$ for all $x \in \mX$. Then $g$ is said to be upper semicontinuous at $x \in \mX$ if
$$\lim\sup_{n\to \infty} g(x_n) \leq g(x)$$ 
for all sequences $\{x_n\}_{n\in\mathbb{N}}$ in $\mX$ that converge to $x$.
\end{definition}

\begin{definition}[Strongly and weakly continuous kernels] \label{def:stronglycontinuous}
The stochastic kernel $Q \in \mathcal{K}(\mS| \mathbb{K})$ is called strongly (resp. weakly) continuous if the function $(s,u) \mapsto \int_\mS g(s') Q(ds'|s,u)$ is continuous and bounded in $\mS$ for every measurable (resp. continuous) and bounded function $g$.
\end{definition}

\begin{proposition}[\cite{lerma}, Proposition A.1] \label{prop:A1}

Let $g$ be a function as in Definition \ref{def:usc}. Then the following statements are equivalent:
\begin{itemize}
\item[(a)] g is upper semicontinuous
\item[(b)] the set $\{(x,c)\in \mX\times\mathbb{R}: \ g(x)\geq c\}$ is closed
\item[(c)] the sets $S_c(g):=\{x\in \mX: \ g(x)\geq c\}$ are closed for all $c\in \mathbb{R}$. 
\end{itemize}
\end{proposition}

\begin{proposition}[\cite{lerma}, Proposition A.2] \label{prop: A2}
Let $L(\mX)$ be the family of all functions on $\mX$ that are lower semicontinuous and bounded below. Then $g \in L(\mX)$ if and only if there exists a sequence $\{g_n\}_n$ of continuous and bounded functions on $\mX$ such that $g_n \uparrow g$.
\end{proposition}

\begin{proposition}[\cite{lerma}, Proposition E.2] \label{prop:E2}
Let $\mu,\mu_1,\mu_2.... \in \mM_\mathbb{P}(\mX)$ such that the sequence $\{\mu_n\}_{n\in\mathbb{N}}$ converges weakly to $\mu$. If $g:\mX \rightarrow \mathbb{R}$ is upper semicontinuous and bounded above, then

$$\lim\sup_{n \to \infty} \int_\mX g(x) d\mu_n(x) \leq \int_\mX g(x) d\mu(x).$$
In other words, the function $\mu \mapsto \ \int_\mX g(x) d\mu(x)$ inherits the upper semicontinuity from $g$.
\end{proposition}

\subsection{Standard stochastic AR MDPs \label{appendix:MDP}} \phantom{aa}\\

This section introduces the assumptions adopted throughout the paper for the stochastic AR problem \cite[Chapter 5]{lerma}, and relevant theoretical results pertaining to techniques T1-T3 (see Section \ref{sec:T1-T3}) in the semicontinuous-semicompact framework of \cite{lerma}.



\subsubsection{Usual assumptions}


The following condition allows for the vanishing discount approach to solve the long-run average problem; that is to say, the function $(1-\alpha)V_{\alpha}^*(\cdot)$ converges pointwise to the long-run average value function $J^*(\cdot)$ when $\alpha \to 1$.

\begin{assumption}[\cite{lerma}, Assumption 5.4.1]\label{assumption541}
\phantom{aaa}
\begin{itemize}
\item[(a)] There exists a state $z \in \mS$ and numbers $\beta \in (0,1)$ and $M \geq 0$ such that 
$$(1-\alpha)V_\alpha^*(z) \leq M \qquad \forall \alpha \in [\beta,1).$$

\item[(b)] Moreover, there is a constant $N \geq 0$ and a nonnegative function $b(\cdot)<\infty$ on $\mS$ such that, with $h_\alpha(s):=V_\alpha^*(s)-V_\alpha^*(z)$
$$-b(s) \leq h_\alpha(s)\leq N \qquad \forall s \in \mS \mbox{ and } \alpha \in [\beta,1) $$
\end{itemize}

\end{assumption}

Under the additional assumption we introduce now, \cite{lerma} prove that the constant $\rho^*$ can be retrieved instead by solving the equations with equality \eqref{AROE}. 

\begin{assumption}[\cite{lerma}, Assumption 5.5.1]\label{assumption551}
Assumption \ref{assumption541} holds and, in addition,
\begin{enumerate}
\item[(a)] The function $b(\cdot)$ in Assumption \ref{assumption541}(b) is measurable and such that

$$ \int_{\mS} b(s') Q(ds'|s,u)<\infty \qquad \forall s\in \mS, u \in \mU(s)$$
\item[(b)] Recall the definition $h_\alpha(s):=V^*_\alpha(s)-V^*_\alpha(z)$ for all $s \in \mS$ for some state $z \in \mS$ given in Assumption \ref{assumption541}, and let $\{\alpha(n)\}$ be sequence of discount factors in Lemma \ref{lemma:vanishing_discount}. Then we assume that the sequence $\{h_{\alpha(n)}\}$ is equicontinuous.
\end{enumerate}
\end{assumption}

To establish the connection between the equilibrium problem and the vanishing discount in the measurized framework we adopt part (a) of the following condition, which is usually made to guarantee that the Policy Iteration Algorithm for the AR problem is {\it well defined} \citep{PIA}. 

\begin{assumption} \label{assumption:welldefined}
\phantom{aa}
\begin{enumerate}
\item[(a)] For every stationary policy $f \in \mathbb{F}$ there exist an invariant measure $\nu_f$; i.e.
\begin{equation}\label{fixedpoint}
\nu_f(\mA)=\int_\mS Q(\mA|s,f(s))\nu_f(ds) \qquad \forall \mA \in \mBS
\end{equation}
and the reward function $r(\cdot,f(\cdot))$ is $\nu_f$-integrable.
\item[(b)] For every $f \in \mathbb{F}$ there is a solution $h_f(\cdot)$ to the second equality equation of the \eqref{AROE} 
\end{enumerate}
\end{assumption}

\subsubsection{Theoretical results for the application of T1-T3.\label{sec:AR}}

In this section, we state the most relevant results pertaining the resolution of stochastic MDPs under the AR criterion. 

\paragraph{ (T1) The AR optimality equations.}

Let $\rho^*$ be a constant, $h(\cdot)$ a real-valued measurable function on $\mS$ and $f \in \mathbb{F}$ a selector. We say that the triplet $(\rho^*,h,f)$  is a solution to \eqref{AROE}. Since it is, in general, difficult to obtain a solution to the equations above, sometimes the following relaxed surrogate, called the AR Optimality Inequality (AROI) is used

\begin{align}\label{AROI}\tag{AROI}
\rho^* +h(s) & \leq \sup_{u \in \mU(s)} \left[r(s,u) + \int_\mS h(s')Q(ds'|s,u) \right] \qquad \forall s \in \mS \\
\rho^* +h(s) & \leq r(s,f(s)) + \int_\mS h(s')Q(ds'|s,f(s))  \qquad \ \forall s \in \mS.  \nonumber
\end{align}

The next theorem connects the constant $\rho^*$ solving the optimality equation \eqref{AROI} to the long-run average value function satisfying \eqref{AR}. This result hence bridges the vanishing discount and average reward approaches.  

\begin{theorem}[\cite{lerma}, Theorem 5.4.3]\label{theorem543}
Suppose Assumption \ref{assumption421} holds. Then
\begin{itemize}
\item[(i)] Under Assumption \ref{assumption541}, there exist a constant $\rho^* \geq 0$ and a measurable function $h: \mS \rightarrow \mathbb{R}$ with
\begin{equation*}
-b(s)  \leq h(s) \leq N \qquad \forall s \in \mS, \mbox{  and  } h(z)=0
\end{equation*}
and a selector $f \in \mathbb{F}$ such that the triplet $(\rho^*,h,f)$ satisfies \eqref{AROI}, and, moreover $\pi_f=\{\varphi_f\}_{t\geq 0}$ such that $\varphi_f(f(s)|s)=1$ for all $s \in \mS$ is AR-optimal and $\rho^*$ is the AR-value function, i.e.
\begin{equation}\label{rhooptimal}
J^*(s)=J(\pi_f,s)=\rho^* \qquad \forall s \in \mS
\end{equation}
so that $\rho^*=\sup_\mS \ J^*(s)=\sup_\mS \sup_\Pi \ J(\pi,s)$; in fact, any selector $f \in \mathbb{F}$ that satisfies the second inequality of \eqref{AROI}, also satisfies \eqref{rhooptimal}.
\item[(ii)] Conversely, if $\pi_f \in \Pi^{DS}$ is AR-optimal and satisfies \eqref{rhooptimal}, then there exists a measurable function $\hat{h}(\cdot)$ on $\mS$ such that the triplet $(\rho^*,\hat{h},f)$ satisfies \eqref{AROI}.
\end{itemize}
\end{theorem}

By also imposing Assumption \ref{assumption551} one can claim the existence of a triplet $(\rho^*,h^*,f^*)$ solving \eqref{AROE}. The following result also states that the stationary deterministic policy $\pi_{f^*}=\{f^*\}_{t \geq 0}$ is AR-optimal and that the AR-value function is constant and equal to $\rho^*$.

\begin{theorem}[\cite{lerma}, Theorem 5.5.4]\label{theorem554}
If Assumptions \ref{assumption421} and \ref{assumption551} hold, then there exist $\rho^*\in \mathbb{R}$, $h^*\in \mC(\mS)$ and $f^* \in \mathbb{F}$ satisfying \eqref{AROE}. Moreover, $\pi_{f^*}=\{f^*\}_{t \geq 0}$ is AR-optimal and $J^*(s)=\rho^*$ for all $s \in \mS$. In fact, any deterministic stationary policy $\pi_{f}$ for which $f$ satisfies the second equation of \eqref{AROE} is AR-optimal.
\end{theorem}

\paragraph{(T2) The linear programming formulation.}

In certain cases the solutions to the optimality equations \eqref{AROE} can also be obtained by solving an LP. In this section we loosely derive it; for details on the necessary assumptions for the formal derivation of this formulation and the derivation process itself, the reader is once again referred to Chapter 6 in \cite{lerma}. 

Recall the first \eqref{AROE} equation

\begin{align*}
\rho^* +h(s) & = \sup_{\mU(s)} \left[r(s,u) + \int_\mS h(s')Q(ds'|s,u) \right] \qquad \forall s \in \mS. 
\end{align*}
Let $(\rho^*,h^*)$ be a solution to above-written equation; then the pair verifies the inequality

\begin{align}\label{rho_ineq}
\rho  & \geq  r(s,u) + \int_\mS h(s')Q(ds'|s,u) -h(s) \qquad \forall s \in \mS
\end{align}
We reformulate inequality \eqref{rho_ineq} into an LP (see Eq. (6.4.12) in \cite{lerma}) as

\begin{align}
\inf_{\rho,h(\cdot)} & \ \rho \label{LP_rho} \tag{LP$_\rho$}\\
\mbox{s.t.} &\  \rho \geq  r(s,u) + \int_\mS h(s')Q(ds'|s,u) -h(s) \qquad \forall s \in \mS, u \in \mU(s),   \nonumber
\end{align}
where $h(\cdot)$ belongs to some suitable set of real-valued functions on $\mS$.

\paragraph{(T3) The vanishing discount approach.}
The following lemma establishes a connection between the vanishing discount value function and the constant $\rho^*$. Conveniently, the notation used for this constant coincides with that of the constant term in \eqref{AROE} and \eqref{AROI}. This is so because these constants indeed coincide, as we will state in a later theorem.

\begin{lemma}\label{lemma:vanishing_discount}
Under Assumption \ref{assumption541}, there exist a constant $\rho^*$ with $0 \leq \rho^* \leq$ M, and a sequence of discount factors $\alpha(n) \uparrow 1$ satisfying
\begin{equation} \label{vanishing_discount}
\lim_{n \rightarrow \infty} (1-\alpha(n))V^*_{\alpha(n)}(s)=\rho^* \qquad \forall s \in \mS
\end{equation}
\end{lemma}


\section{Lifting the stochastic AR optimality equations \label{sec:measurizedAROE}}\phantom{aa}\\

To lift the AR optimality equations, we follow the procedure outlined in Section \ref{sec:measurizing}. We start with the LP formulation \eqref{LP_rho} of the stochastic optimality equations \eqref{AROE}. Then we aggregate the constraints over every Markov decision rule $\varphi$:

\begin{align}
\inf_{\rho,h} & \ \rho \label{measurized_AR0} \\
\mbox{s.t.} & \ \rho+h(s)  \geq \int_{\mathcal{U}}  r(s,u) \varphi(du|s)+ \int_{\mathcal{U}} \int_\mS h(s') Q(ds'|s,u)  \varphi(du|s) \qquad \qquad \forall  s \in \mathcal{S}, \varphi \in \Phi. \nonumber
\end{align}

Then, we aggregate the constraints with respect to state distributions $\nu \in \mMpS$:

\begin{align}
\inf_{\rho,h} & \ \rho \label{measurized_AR0} \\
\mbox{s.t.} & \ \medmath{\rho+ \int_{\mS} h(s) d\nu(s) \geq  \int_{\mS} \int_{\mathcal{U}}  r(s,u) \varphi(du|s) d\nu(s)+ \int_{\mS} \int_{\mathcal{U}} \int_\mS h(s') Q(ds'|s,u)  \varphi(du|s) d\nu(s)}, \quad \forall  \nu \in \mMpS, \varphi \in \Phi. \nonumber
\end{align}
To ease notation, define the measurized version of the bias function $\ovh(\cdot)$ as $\ovh: \mMpS \to \mathbb{R}$, given by

\begin{equation} \label{H} 
\ovh(\nu):=\int_\mS h(s) d\nu(s) \qquad \forall \nu \in \mMpS.
\end{equation}
This is consistent with previous definitions and results in this paper (see the proof of Lemma \ref{lemma:assumption541} and Theorem \ref{th:measurized_AROE} for more details). By substituting the definitions of the measurized reward $\ovr$ and the deterministic transition $F$ into \eqref{measurized_AR0}, we obtain \eqref{measurized_LPrho}.

\section{Measure-valued MDPs \label{sec:measure-valued}}

In this section, we formally introduce measure-valued MDPs. Essentially, these are MDPs whose states are measures, although not necessarily probability measures. These MDPs may often be a lifted version of a standard MDP, \textcolor{black}{as we showed in Section \ref{sec:measurized}}. \\

\begin{definition}\label{def:MvMDP}
A measure-valued Markov Decision Model ($\mM v$-MDP) is a five-tuple
\begin{equation}
(\omS,\omU,\{\omU(\nu)| \ \nu \in \omS\},\oQ,\ovr)
\end{equation}
where

\begin{itemize}
\item[(i)] $\omS$ is a space of measures defined over a Borel sample space $\mS$
\item[(ii)] $\omU$ is the set of actions, assumed to be a Borel space
\item[(iii)] $\{\omU(\nu)| \ \nu \in \omS\}$ is the family of nonempty measurable subsets $\omU(\nu)$ of $\omU$, where $\omU(\nu)$ denotes the set of feasible actions when the system is in state $\nu$, and with the property that the set 
$$\overline{\mathbb{K}}=\{ (\nu,\overline{u})| \ \nu\in \omS,\ \overline{u}\in \omU(\nu)\}$$
 of feasible state-action pairs is a measurable subset of $\omS \times \omU$
\item[(iv)] $\oQ$ is a stochastic kernel on $\omS$ given $\overline{\mathbb{K}}$
\item[(v)] $\ovr:\ \overline{\mathbb{K}} \rightarrow \mathbb{R}$ is the reward-per-stage function, assumed to be measurable\\
\end{itemize}
\end{definition}

In our definition of measure-valued MDPs, we specified that $\omU$ and $\mS$ must be Borel spaces. If we limit the state space $\omS$ to the set of probability measures on $\mS$, denoted as $\mMpS$, and equip it with the weak convergence topology, then $\mMpS$ becomes a Borel space. By imposing the standard assumptions (Assumptions \ref{assumption421} and \ref{assumption422}), we can apply Theorem \ref{theorem423}. This implies that we can retrieve the optimal value function and controls from the measure-valued optimality equations

\begin{equation}
\overline{V}^*(\nu)=\sup_{\overline{u} \in \overline{\mU}(\nu)} \ \left\{ \ovr(\nu,\varphi) + \alpha   \mathbb{E}_{\overline{Q}} \left[ \overline{V}^*(\nu')| \nu,\varphi \right] \right\}. \qquad \forall \nu \in \overline{S}, \label{MV-DCOE} \tag{$\mM v$-$\alpha$-DCOE}
\end{equation}

In addition, whenever $\mS$ is the state space of a standard MDP, we can think of a measure-valued MDP as controlling the distribution of states rather than their realizations. This is particularly useful when we the states are not completely observable (as in POMDPs) or when we are managing large populations (as in certain MFMDPs). The following examples formulate partially-observable and mean-field MDPs as measure-valued MDPs.\\

\begin{example}[Mean-field MDPs] \label{ex:MFMDP}
Consider $I$ i.i.d. individuals with states $s_1,...,s_I$ coming from distribution $\nu\in \mMpS$. In a mean-field control problem, these individuals are cooperative agents aiming to maximize the overall social benefit of the system. In this context, the reward function of each agent needs to take into account not only her current state but also the empirical distribution of the states of all the other agents; i.e.,

\begin{equation} \label{empirical_nu}
\hat{\nu}(\mA)=\frac{1}{I} \sum_{i=1}^I \mathbf{1}_{\left\{s_i \in \mA\right\}} \quad \forall \mA \in \mBS.
\end{equation}

Typically, one assumes that the action set $\mU(\bs)$ can be decomposed as $\Pi_{i=1}^I \mU(s_i)$; i.e., there are no linking constraints for the controls. With this notation, the reward function can thus be expressed as 

\begin{equation}\label{MFMDP_reward}
R(\bs,\bu)=\frac{1}{I} \sum_{i=1}^I r(s_i,u_i,\hat{\nu}),
\end{equation}
where $\bs=(s_1,...,s_I)$ and $\bu=(u_i,...,u_i)\in \mU(\bs)$. The transitions are also performed independently and identically, albeit these are random and depend on the i.i.d. individual random shocks $Z_i$, $i=1,...,I$, and the common shock $Z_0$. More specifically, we can decompose the transition function $\mathbf{T}=\Pi_{i=1}^I T(s_i,u_i,\hat{\nu},Z_0,Z_i)$, where

\begin{equation} \label{MFMDP_transition}
\mathbb{P}_{q,Z_0,Z_i}(s_i'\in \mA|s_i,u_i,\hat{\nu},Z_0,Z_i)=T(s_i,u_i,\hat{\nu},Z_0,Z_i)(\mA) \qquad \forall \mA \in \mB(\mS).
\end{equation}
Note that $\mathbb{P}_{q,Z_0,Z_i}$ depends not only on transition kernel $q$ governing the transition of states but also on the probability distribution of the noises. 

Under certain assumptions\footnote{Continuity of the reward function and compactness of the set $\mathcal{Z}$, where the shocks belong to, plus the usual assumptions (Assumptions \ref{assumption421} and \ref{assumption422}).}, \cite{bauerle} shows that when the number of individuals $I$ goes to infinity, this MDP is equivalent to an MDP in the space of probability measures. More specifically, states $\nu \in \mMpS$ are probability measures on $\mS$ and actions $\gamma\in \omU$ are joint probability measures on $\mS \times\mU$. The set of admissible actions from state $\nu$ is

\begin{equation*}
\omU(\nu)=\left\{ \gamma \in \mathcal{M}_{\mathbb{P}}(\mS\times\mU): \ \gamma(\mA\times\mU)=\nu(\mA) \mbox{ and } \int_{s \in \mA} \int_{\mU(s)^c} \gamma(ds,du)=0 \ \ \forall \mA \in \mB(\mS) \right\},
\end{equation*}
where $\mU(s)^c$ is the set of inadmissible actions from state $s$ in the original MDP. The reward function of the measure-valued MDP is the expected reward of the original MDP; i.e., 

\begin{equation} \label{r_bauerle}
\ovr(\nu,\gamma)= \mathbb{E}_\nu \mathbb{E}_\gamma [r(s,u,\hat{\nu})].
\end{equation} 
The next state is a random measure that depends on the common noise $Z_0$ and the individual noise $Z_i$. \cite{bauerle} characterizes this transition through mapping 

\begin{align*}
\overline{T}: \mMpS \times \omU \times \mathcal{Z} & \rightarrow \mMpS\\
(\nu,\gamma,Z_0) & \mapsto \nu'(\cdot)=\int_\mS \int_\mU p(\cdot|s,u,\nu,Z_0) \gamma(du,s) \nu(ds)
\end{align*}
where $p(\mA|s,u,\nu,Z_0)=\mathbb{P}_{q,Z_i}(T(s,u,\hat{\nu},Z_0,Z_i)\in \mA | Z_0)$ inherits the randomness from individual transition $q$ and random noise $Z_i$. The transition function $\overline{T}$ relates to the transition kernel $\oQ$ of the measured-value MDP through $\oQ(\mathcal{P}|\nu,\gamma):=\mathbb{P}(\nu' \in \mathcal{P})=\mathbb{P}_{q,Z_0,Z_i}(\overline{T}(\nu,\gamma,Z_0) \in \mathcal{P})$ for all $\mathcal{P} \in \mB(\mMpS)$. Therefore, one could express an MFMDP as the measure-valued MDP $(\mMpS,\overline{U},\{\omU(\nu)| \ \nu \in \omS\},\oQ,\ovr)$ with the specifications given above. Here \cite{bauerle} showed that the measure-valued MDP inherits Assumptions \ref{assumption421} and \ref{assumption422} from the original MDP, albeit also assuming that the reward function $r$ is continuous (an assumption we do not make). Therefore, the optimal control $\bu^*$ and the optimal value function $\oV^*(\cdot)$ can be retrieved through optimality equations \eqref{MV-DCOE}. More details on \cite{bauerle} and how it relates to the measurized MDPs proposed later in this paper can be found in Appendix \ref{appendix:bauerle}. \hfill $\blacksquare$\\

\end{example}

\begin{example}[Partially Observable MDPs]
Other practical examples of MDPs that can be framed within the measure-valued framework are uncertain MDPs. In a standard MDP, the agent has complete information about the current state of the environment, allowing it to make optimal decisions based on that information. In contrast, in a POMDP the agent does not have direct access to the true state of the environment. Instead, it observes partial, noisy, or incomplete information about the state through observations. This lack of complete information introduces uncertainty and makes decision-making more challenging.

Mathematically, a POMDP can be defined as the six-tuple $(\mS,\mU,Q,r,\mathcal{O},\omega)$, where the four first elements are as in a standard MDP, $\mathcal{O}$ denotes the space of observations, and $\omega(\mA|s',u)$ is the probability of perceiving an observation in $\mA \in \mB(\mathcal{O})$ when action $u$ has been implemented and the process has transitioned to $s'$. Note that we do not specify the set of admissible actions because these are independent of the state; i.e. $\mU(s)=\mU$ for all $s\in \mS$.

It is well known that a POMDP can be modelled as a BMDP, in which the agent keeps track of a belief state $\nu\in \mMpS$ representing the current prior distribution over states $s\in\mS$. The action space $\omU$ coincides with the action space $\mU$ of the POMDP and is also independent of the state. Given the decision maker takes action $u \in \mU$ while in state $\nu \in \mMpS$, the state transitions to $\nu'\in \mMpS$ according to Bayes rules. More specifically, after observing $o$ the new prior becomes

\begin{equation} \label{nu_o}
\nu'_{o}(\mA)=\int_{s \in \mS} \int_{s' \in \mA} \eta(o,u) \omega(o|s',u) Q(ds'|s,u) d\nu(s),
\end{equation}
where $\eta(o,u)=1/\int_\mS \int_\mS \omega(o|s',u) Q(ds'|s,u) d\nu(s)$ is a normalizing constant. Therefore, the transition kernel $\oQ$ of the measure-valued MDP is characterized by

$$\hspace{6cm} \oQ(\mathcal{P}|\nu,u)=\mathbb{P}_\omega (\nu'_o \in \mathcal{P}). \hspace{6cm} \blacksquare$$
%

\end{example}


Much like the earlier examples and measurized MDPs, measure-valued MDPs can sometimes be linked to an underlying MDP with state space $\mS$ and action space $\mU$. Formulating a measure-valued MDP on the space of probability measures over $\mS$ gives rise to states $\nu \in \mMpS$. Interestingly, the set of admissible actions $\omU(\nu)$ is contingent on the current state distribution $\nu$. This enables the modelling of probabilistic constraints on the states of the original MDP. For instance, one can set constraints on various risk measures related to the agent's future perceived costs or limit moments of future distributions over original states. Furthermore, if $\omU$ is the set of Markov decision rules in the original MDP, i.e.
 
 \begin{equation}\label{action=varphi}
 \omU=\{ \varphi \in \mathcal{K}(\mU|\mS): \ \varphi(\mU(s)|s)=1 \ \forall s \in \mS\},
 \end{equation}
 one could also impose restrictions on the distribution of actions taken in the original MDP. These constraints do not arise naturally outside the measure-valued framework. The following examples aim to illustrate how such requirements could easily be added in the proposed framework.\\

\setcounter{example}{0}

\begin{example}[MFMDPs Revisited]
As mentioned in \citet[Remark 3.2]{bauerle}, instead of considering actions as joint probability measures $\gamma \in \mM_\mathbb{P}(\mS\times \mU)$, one could consider actions $\varphi$ to be stochastic kernels belonging to the space $\omU$ defined in \eqref{action=varphi}\footnote{In Section \ref{sec:change_variables} we rigorously explore under which circumstances such a change of variables can be performed without loss of generality. In Section \ref{sec:dualizing} we outline how these joint measures are related to the dual variables of the LP formulation \eqref{LP}.} Now assume that we want to bound the variance of the actions taken by the pool of cooperative agents. Therefore, the set of admissible controls is


$$\omU(\nu)=\{ \varphi \in \omU: \ \mathbb{E}_\nu \mathbb{E}_\varphi[ (u-\mu_u)^2|s]\leq \theta\}, \qquad \forall \nu \in \mMpS,$$
where $\mu_u$ is the mean value of the actions and $\theta\geq 0$ is a parameter. Then it suffices to add the constraint

$$\int_\mS \int_\mU u^2 \varphi(du|s) d\nu(s) - \left[ \int_\mS \int_\mU u \varphi(du|s) d\nu(s) \right]^2 \leq \theta$$
to the optimality equations \eqref{MV-DCOE}. \hfill $\blacksquare$\\

\end{example}


\begin{example}[POMDPs Revisited]
In this example, we assume that, although the transited states cannot be known with certainty, the agent wants to bound the expected probability of landing in states belonging to the set $\mA \in \mB(\mS)$. That is to say, the agent wants to consider the following set of admissible actions 

$$\omU(\nu)=\{ u \in \mU: \ \mathbb{E}_\omega[\mathbb{P}_{\nu_o'}(s' \in \mA))| o] \leq \theta\}, \qquad \forall \nu \in \mMpS.$$
Here, $\nu_o'$ represents the subsequent distribution of states as defined in \eqref{nu_o} and $\theta\geq 0$ is a parameter. This requirement is seamlessly expressed by adding constraint

$$ \int_{s \in \mS} \int_{s' \in \mA} \int_\mathcal{O} \eta(o,u) \omega(do|s',u) Q(ds'|s,u) d\nu(s) \leq \theta $$
to the optimality equations \eqref{MV-DCOE}.  \hfill $\blacksquare$\\
\end{example}




%

%

\section{Facilitating the derivation of measure-valued MFMDPs \label{appendix:bauerle}} 

Consider the MFMDP introduced in Example \ref{ex:MFMDP}. Since the set of feasible actions can be decomposed as $\mU(\bs)=\Pi_{i=1}^I \mU(s_i)$, there are no linking constraints across components $i$ and the set of feasible actions for each agent is the same. In addition, the transitions are i.i.d. over time $\mathbf{T}=\Pi_{i=1}^I T(s_i,u_i,\hat{\nu},Z_0,Z_i)$, thus having that the MFMDP kernel $Q$ can also be somewhat  decomposed as $Q(\cdot|\bs,\bu,Z_0,Z_i)=\Pi_{i=1}^I q(\cdot|s_i,u_i,\hat{\nu},Z_0,Z_i)$, where

\begin{equation*} 
q(s_i'\in \mA|s_i,u_i,\hat{\nu},Z_0,Z_i)=T(s_i,u_i,\hat{\nu},Z_0,Z_i)(\mA) \qquad \forall \mA \in \mB(\mS).
\end{equation*}

We can conceptualize the MFMDP as consisting of $I$ MDPs, one for each agent. Although each agent shares the same reward function and transition kernel, and there are no constraints directly linking their decisions, decomposing the MFMDP into $I$ separate MDPs is not straightforward. This complexity arises because both the reward function and the transitions depend on the distribution $\hat{\nu}$, preventing a straightforward decomposition by state. Nonetheless, \cite{bauerle} builds on the special structure of MFMDPs to show that the problem can be equivalently solved by a unidimensional MDP in the space of measures.  The connection in that paper is built through the empirical MDP, defined over empirical measures. More specifically, \cite{bauerle} defines the empirical value function as

\begin{equation}\label{bauerle_empirical}
\hat{V}_I^*(\hat{\nu}):=\sup_{\hat{\Gamma}} \sum_{t=0}^\infty \alpha^t \mathbb{E}_{\hat{\nu}}^{\hat{\pi}} [\hat{r}(\hat{\nu}_t,\hat{\gamma}_t)]
\end{equation}
for any empirical measure $\hat{\nu}\in \mMpS$ and a policy $\hat{\Gamma}:=(\hat{\gamma}_1,\hat{\gamma}_2,....)$ composed of empirical joint distributions $\hat{\gamma}_t \in \mM_{\mathbb{P}}(\mS\times\mU)$. The reward function $\hat{r}$ of the empirical process is similar to our measurized reward but evaluated on empirical measures

\begin{equation}\label{R_empirical}
\hat{r}(\hat{\nu},\hat{\gamma}):=\int_\mS \int_\mU r(s,u,\hat{\nu}) d\hat{\gamma}(s,u).
\end{equation}
Subsequently, \cite{bauerle} uses the following lemma to show that the optimal value functions of the original and the empirical MDPs coincide: i.e., $V^*(\bs)=V^*_I(\hat{\nu})$ for all $\bs \in \mbS$ and $\hat{\nu}$ as in \eqref{empirical_nu}. 

\begin{lemma}[Lemma 3.1, \cite{bauerle}] \label{lemma:bauerle}
For any feasible state-action pair $(\bs,\bu)\in \mathbb{K}$, there exists an empirical joint distribution $\hat{\gamma} \in \mM_{\mathbb{P}}(\mS\times\mU)$ such that

\begin{equation}\label{bauerle_Rr}
R(\bs,\bu)=\hat{r}(\hat{\nu},\hat{\gamma}),
\end{equation}
where $\hat{\gamma}(\mA,\mU)=\hat{\nu}(\mA)$ for all $\mA \in \mBS$. 
\end{lemma}

Note that \eqref{bauerle_Rr} already implies a state-space collapse: a unique empirical measure suffices to solve the problem with $I$ i.i.d. agents. Similarly to \eqref{bauerle_empirical}, in the next step \cite{bauerle} defines a unidimensional measure-valued process, with value function

\begin{equation}\label{bauerle_measurized}
\oV^*(\nu_0):=\sup_{\Gamma} \sum_{t=0}^\infty \alpha^t \mathbb{E}_{\nu_0}^{\pi} [\overline{r}(\nu,\gamma)], 
\end{equation}
where the reward is defined as in \eqref{R_empirical} but can be evaluated over any probability measures $\nu \in \mMpS$ and $\gamma \in \mM_\mathbb{P}(\mS\times\mU)$, rather than solely on empirical distributions.  In the limit, \cite{bauerle} shows that this measure-valued process retrieves the original value function; i.e., if $\bs=(s_1,....,s_I)\in \mbS$ is such that $s_i \sim \nu$ for all $i=1,...,I$, and $\hat{\nu} \to \nu$ as $I \to \infty$, then $\lim_{I \to \infty} V^*(\bs)=\lim_{I \to \infty} \hat{V}_I^*(\hat{\nu})=\oV^*(\nu).$

As noted in \cite{bauerle}, in the absence of random shocks, the lifted MDP is a deterministic process, much like our measurized MDP. Therefore, \cite{bauerle} is able to leverage the structure of the MFMDP and employ sophisticated mathematical machinery to obtain an alternative derivation of our measurized MDP for MFMDPs. Although similar in spirit, the approach presented in this paper is simpler, more general and encompasses any kind of standard MDP. Here we show how the measurized theory herein can facilitate the analysis of MFMDPs, arriving to the unidimensional measurized MDP \eqref{bauerle_measurized} in a simple and intuitive manner. We just use the assumption on continuity and boundedness of reward $r$ considered in \cite{bauerle}, and Lemma \ref{lemma:bauerle} above. The measurized reward is

\begin{equation} \label{MFMDP-R}
\overline{R}(\bnu,\bvarphi):=\int_\mbS \int_\mbU R(\bs,\bu) \bvarphi(d\bu|\bs)d\bnu(\bs),
\end{equation}
where $\bnu$ and $\bvarphi$ are i.i.d. according to the definition of MFMDP. In addition $\lim_{I \to \infty} \hat{\nu}=\nu$ by the law of large numbers. 
Our Lemma \ref{lemma:change_variables} on the decomposition of $\gamma$ yield

\begin{align*}
R(\bs,\bu;\hat{\nu}) &=\int_\mS \int_\mU r(s,u,\hat{\nu}) d\hat{\gamma}(s,u)\nonumber\\
&= \int_\mS \int_\mU r(s,u,\hat{\nu}) \varphi(du|s) d\hat{\nu}(s)\nonumber
\end{align*}
Now we show that 

\begin{equation}\label{MFMDP-r} 
 \lim_{I \to \infty} R(\bs,\bu;\hat{\nu})= \ovr(\nu,\varphi):=\int_\mS \int_\mU r(s,u,\nu) \varphi(du|s) d\nu(s).
\end{equation}
We follow the easy proof in \cite{bauerle}:
\begin{align*}
&\left|  \int_\mS \int_\mU r(s,u,\hat{\nu}) \varphi(du|s) d\hat{\nu}(s) -  \int_\mS \int_\mU r(s,u,\nu) \varphi(du|s) d\nu(s)\right|& \nonumber\\
&\qquad \leq  \int_\mS \int_\mU \left| r(s,u,\hat{\nu})-r(s,u,\nu)\right| \varphi(du|s) d\hat{\nu}(s) + \left| \int_\mS \int_\mU r(s,u,\nu) \varphi(du|s) (d\hat{\nu}(s)-d\nu(s))\right| &\nonumber\\
&\qquad \leq \sup_{(s,u) \in \mathbb{K}} \left| r(s,u,\hat{\nu})-r(s,u,\nu)\right|  + \left| \int_\mS \int_\mU r(s,u,\nu) \varphi(du|s) (d\hat{\nu}(s)-d\nu(s))\right|& \nonumber
\end{align*}
The first term converges to zero because of the assumption on continuity of $r$ and Lemma 8.1 in \cite{bauerle}. The author claims that the second term converges to zero because of the continuity and boundedness of $r$, and because $\hat{\nu} \to \nu$ when $I \to \infty$. Putting together \eqref{MFMDP-R} and \eqref{MFMDP-r} gives

\begin{align*}
\lim_{I \to \infty} \overline{R}(\bnu,\bvarphi)&= \lim_{I \to \infty} \int_\mbS \int_\mbU R(\bs,\bu)  \bvarphi(d\bu|\bs)d\bnu(\bs)\\
&= \int_\mbS \int_\mbU  \ovr(\nu,\varphi)  \bvarphi(d\bu|\bs)d\bnu(\bs)=\ovr(\nu,\varphi).
\end{align*}
where the second equality comes from the MCT (Assumptions \ref{assumption421} and \ref{assumption:MCT} are necessary). Similarly, the measurized transition in the absence of random shocks is defined as

\begin{align*}
F(\nu,\varphi)(\cdot)=\int_\mS \int_\mU q(\cdot|s,u,\hat{\nu}) \varphi(du|s) d\nu(s).
\end{align*}
Since $q$ is strongly continuous, the function $(s,u,\hat{\nu}) \mapsto q(\cdot|s,u,\hat{\nu})$ is continuous, thus having

\begin{align*}
\lim_{I \to \infty} F(\nu,\varphi)(\cdot)& = \lim_{I \to \infty}\int_\mS \int_\mU q(\cdot|s,u,\hat{\nu}) \varphi(du|s) d\nu(s)\\
& = \int_\mS \int_\mU q(\cdot|s,u,\lim_{I \to \infty}\hat{\nu}) \varphi(du|s) d\nu(s)\\
&=\int_\mS \int_\mU q(\cdot|s,u,\nu) \varphi(du|s) d\nu(s).
\end{align*}

The unidimensional measurized MDP already possesses information on the state distribution of all other agents, condensed in measure $\nu$, because agents are assumed i.i.d. Therefore, the measurized transition of a single agent already contains all the necessary information so that the transition is performed independently using solely the measurized information of one agent. 
This allows us to decompose the measurized MFMDP into a unidimensional problem on measures, yielding \eqref{bauerle_measurized}.



%
\bibliographystyle{jf}
\bibliography{references} 
\end{document}